\documentclass[12pt,reqno]{amsart}

\title[]{Unique Ergodicity in Stochastic Electroconvection}

\author{Elie Abdo}
\address{Department of Mathematics, Temple University, Philadelphia, PA 19122}
\email{abdo@temple.edu}
\author{Nathan Glatt-Holtz}
\address{Department of Mathematics, Tulane University, New Orleans, LA 70118}
\email{negh@tulane.edu}
\author{Mihaela Ignatova}
\address{Department of Mathematics, Temple University, Philadelphia, PA 19122}
\email{ignatova@temple.edu}

\usepackage{mathtools}
\usepackage[margin=1in]{geometry}
\usepackage{amsmath, amsthm, amssymb}
\usepackage{times}
\usepackage{color}
\usepackage{hyperref}
\newcommand{\pa}{\partial}
\newcommand{\la}{\label}
\newcommand{\fr}{\frac}
\newcommand{\na}{\nabla}
\newcommand{\be}{\begin{equation}}
\newcommand{\ee}{\end{equation}}
\newcommand{\ba}{\begin{array}{l}}
\newcommand{\ea}{\end{array}}

\newtheorem{prop}{Proposition}

\newcommand{\beg}{\begin}

\renewcommand{\l}{\Lambda}

\usepackage{cleveref}

\numberwithin{equation}{section}

\usepackage{MnSymbol,wasysym}

\newcommand{\R}{\mathbb R}

\def\ZZ{{\mathbb Z}}
\def\RR{{\mathbb R}}
\def\TT{{\mathbb T}}

\def\PP{\mathbb P}
\def\d{\mathrm{\textbf{d}}}
\def\E{\mathrm{\textbf{E}}}

\begin{document}
\begin{abstract} 
We consider a stochastic electroconvection model describing the nonlinear evolution of a surface charge density in a two-dimensional fluid with additive stochastic forcing.  We prove the existence and uniqueness of solutions, we define the corresponding Markov semigroup, and we study its Feller properties.  When the noise forces enough modes in phase space, we obtain the uniqueness of the smooth invariant measure for the Markov transition kernels associated with the model.    
\end{abstract} 

\vspace{.5cm}

\maketitle

\tableofcontents

\section{Introduction} \la{intro}

Electroconvection refers to the dynamics of electrically conducting fluids under the influence of electrical charges. There are many instances of electroconvection in non-Newtonian and Newtonian fluids, including the flow of nematic and smectic suspensions subject to applied voltage. A particularly interesting example \cite{DDMB,DDMT} considers the dynamics of a thin smectic film in an annular region, driven by an imposed voltage at the boundary. In \cite{ceiv} the behavior of the system was investigated mathematically in the absence of stochastic forcing. The model was described in terms of a surface charge density $q$,
an electric field $E$ and a fluid velocity $u$. The dynamics were confined to a two dimensional domain ($\mathbb T^2$ in the present paper).
The electric field $E$ was derived from a time independent potential $\Phi$ representing the voltage applied at the boundary and a dynamic potential $\l^{-1}q$ due to the charge density $q$, via the relation
\be \la{eq1}
E = - \na \Phi - \na \l^{-1} q,
\ee 
where $\l^{-1}$ denotes the  inverse of the square root of the two-dimensional periodic Laplacian $\l~=~\sqrt{-\Delta}$. 
The current density due to the fluid and the electric field $E$ is  
\be \la{eq2}
J = E + qu, 
\ee
and the charge density obeys the continuity equation
\be \la{eq3}
\pa_t q + \nabla \cdot J = 0.
\ee
The fluid velocity $u$ obeys the incompresible Navier-Stokes equation forced by the electrical forces $qE$ and time independent body forces $f$,  
\be \la{eq4}
\pa_t u + u \cdot \nabla u - \nu \Delta u  + \nabla p  = qE  + f, \quad \nabla \cdot u=0
\ee 
where $p$ is the fluid pressure and $\nu$ is the kinematic viscosity. 

The well-posedness and global regularity of the deterministic model \eqref{eq1}--\eqref{eq4} were obtained in \cite{ceiv} in bounded domains with homogeneous boundary conditions, and the long-time dynamics were investigated in \cite{AI} in the two-dimensional torus $\TT^2$. 

In this paper we consider the stochastic electroconvection model corresponding to the deterministic model \eqref{eq1}--\eqref{eq4},  
\be \la{intro1}
\d q + \na \cdot (qu - \na \l^{-1}q - \na \Phi)dt  = \tilde{g}dW,
\ee 
\be \la{intro11}
\d u + u \cdot \na u dt + \na p dt - \nu \Delta u dt = -q (\na \l^{-1}q + \na \Phi) dt + f dt + g dW,
\ee 
\be \la{intro2}
\na \cdot u = 0
\ee 
forced by time independent noise processes $\tilde{g}dW$ and $gdW$ on  $\TT^2$. For simplicity, we assume that $\nu = 1$. We address the global well-posedness of \eqref{intro1}--\eqref{intro2}, the Feller properties of the Markov semigroup associated with \eqref{intro1}--\eqref{intro2}, and the existence, uniqueness and regularity of the invariant measures for the Markov transition kernels associated with the model \eqref{intro1}--\eqref{intro2}. 

A vast literature treats the well-posedness of stochastic partial differntial equations. Martingale type approaches \cite{B, BT, DPZ, FG, S} were established to prove the existence and uniqueness of solutions to the two-dimensional stochastic Navier-Stokes equations (NSE). In \cite{MS}, the authors use a different approach, independent of the pathwise solutions, based on a generalization of the classical Minty-Browder local monotonicity argument \cite{M,P}, to establish the well-posedness to the stochastic NSE in bounded and unbounded domains. Global existence and uniqueness of strong pathwise solutions were obtained as well for the two-dimensional \cite{EPT, GT, GZ} and three-dimensional \cite{DGTZ, GH} stochastic primitive equations. 

The stochastic electroconvection model \eqref{intro1}--\eqref{intro2} is nonlocal, nonlinear, with critical dissipation in one equation, and consequently the proof of its global well-posedness is rather technical. Under low regularity assumptions imposed on the noises (namely $L^4$ for $\tilde{g}$ and $H^1$ for $g$), we prove that the system \eqref{intro1}--\eqref{intro2} has unique global solutions when the deterministic initial charge density is $L^4$ regular and the deterministic initial  velocity is $H^1$ regular. The existence of solutions is obtained by taking a viscous approximation of \eqref{intro1}--\eqref{intro2}, establishing uniform bounds for the viscous solutions, and using weak convergence. The identification of the drift 
is highly challenging. The reason is that the nonlinearity $qRq$ is not weakly continuous in the spaces we have control in. The remedy is a coercive estimate \eqref{contidea} and use of ideas from \cite{MS}. As a consequence of the existence result, we define the Markov transition kernels on $L^4 \times H^1$ and we show that they are Feller in the norm of $H^{-\fr{1}{2}} \times L^2$. If the noises have higher regularity (namely $\na \tilde{g} \in L^8$ and $\Delta g \in L^2$), then the Markov kernels become Feller in the stronger norm of $H^1 \times H^1$.

We also study the ergodicity of the electroconvection model \eqref{intro1}--\eqref{intro2}. The existence of an invariant measure for the stochastic NSE system was obtained in \cite{DPZ, F, FM}, and the ergodic theory for the stochastic NSE became the center of interest of many subsequent papers (cf.~\cite{BKL2, DX, EH, Fe, M3, M, M2, DX} and references therein). Existence and regularity of invariant measures were obtained in \cite{GKVZ} for the three-dimenional stochastic primitive equations.  In \cite{CGV}, existence and uniqueness of an ergodic invariant measure was established for the 2D fractionally dissipated periodic stochastic Euler equation.

The dissipative term $\l q$ in \eqref{intro1} is critical, and this is a source of technical difficulty. When the potential $\Phi$ vanishes, and with a low regular noise process, we use the Krylov Bogoliubov averaging procedure to prove that the stochastic model \eqref{intro1}--\eqref{intro2} has an invariant measure supported on $H^{\fr{1}{2}} \times H^2$. 
If the noise processes are smooth then the invariant measures are smooth. This follows from bounds of the form
\be \la{ii1}
\E \int_{0}^{T} (\|q\|_{H^{\fr{3}{2}}}^2 + \|u\|_{H^2}^2) dt \le \Gamma_1(u_0, q_0, f, g, \tilde{g}) + \Gamma_2(f,g,\tilde{g})T  
\ee and 
\be \la{ii2}
\E \int_{0}^{T} \log(1 + \|q\|_{H^{k + \fr{1}{2}}}^2 + \|u\|_{H^{k+2}}^2)  dt\le  \Gamma_{1,k}(u_0, q_0, f, g, \tilde{g}) + \Gamma_{2,k}(f,g,\tilde{g})T  
\ee for $k \ge 1$. 
These bounds are obtained by taking advantage of the smoothing
properties of the Stokes operator and the nonlinear coupling, and employing the logarthmic strategy introduced in \cite{GKVZ}. 

The question of uniqueness of invariant measures requires a deeper structural understanding of the interplay of the dynamics and stochastic perturbation. A number of approaches have been used in the recent literature (\cite{BKL2,EMS,HM,HMS, KS,KS2,KS3,M,M2} and references therein).  In this paper we use the asymptotic coupling approach introduced in \cite{H} and  \cite{HMS}. The asymptotic coupling framework was used in \cite{GMR} to obtain uniqueness of the invariant measures of  stochastically forced Navier-Stokes equations, fractionally dissipative Euler equations and damped nonlinear wave equations. In order to show that a stochastic differential equation
\be \la{C1}
dy = F(y) dt + \sum\limits_{l=1}^{d} \sigma_l dW_l   
\ee with initial data $y(0) = y_0$ has only one ergodic measure, the idea is to build a  copy
\be \la{C2}
d\tilde{y} = F(\tilde{y}) dt + G(y, \tilde{y}) \textbf{1}_{t \le \tau} dt + \sum\limits_{l=1}^{d} \sigma_l dW_l   
\ee 
where the feedback control $G$ is such that $y$ and $\tilde{y}$ are forced to approach each other, $y(t) - \tilde{y}(t) \rightarrow 0$ in an appropriate norm, on the event $\left\{\tau = \infty \right\}$  where  $\tau$ is a stopping time such that the coupled system \eqref{C1}--\eqref{C2} has global solutions
with initial data $\tilde{y}(0) = \tilde{y}_0$, and  $\PP(\tau = \infty) > 0$ . Moreover, it is required that 
\be 
\int_{0}^{\infty} |\sigma^{-1}G(y(t), \tilde{y}(t))|^2 \textbf{1}_{t \le \tau} dt < C  
\ee holds (for a.e. $w \in \Omega$) for some deterministic constant $C>0$. If such a construction can be done, then \eqref{C1} has a unique ergodic invariant measure. Finding an appropriate feedback $G$ is typically based on splitting a Hilbert space $X$ into the direct sum of a finite-dimensional space $X_{low}$ and an infinite-dimensional space $X_{high}$ 
\be 
X = X_{low} \oplus X_{high}
\ee 
in such a way that the long time dynamics are controlled by the low frequency part in $X_{low}$. More precisely, the property used is that if the low frequency parts of two solutions are asymptotically the same, then the
high frequency parts in $X_{high}$ are also asymptotically the same. Accordingly, two realizations of \eqref{C1} are coupled in such a way that that their low frequency parts coincide for large time $t > \tau$ provided that they meet at time $t = \tau$. 

 The uniqueness of the invariant measure of the electroconvection model \eqref{intro1}--\eqref{intro2} is obtained by constructing an appropriate feedback control and stopping time. The construction requires  $L^4$ bounds for $q$ and $H^2$ bounds for $u$, and exponential martingale estimates and the Burkholder-Davis-Gundy inequality. The main difficulty is due to the weaker dissipation of the charge densities, and here we use ideas from \cite{GMaR}  to estimate the feedback control.

This paper is organized as follows. In section \ref{EUS}, we show that the system \eqref{intro1}--\eqref{intro2} has a unique global solution provided that the initial charge density has a zero spatial average and is $L^4$ integrable, the initial velocity is divergence-free and is weakly differentiable, and the noise is sufficiently regular. The proof is based on uniform estimates in Lebesgue spaces which are established in Appendix A. In section \ref{sec3}, we define the semigroup associated with \eqref{intro1}--\eqref{intro2} and we prove that it is weak Feller. In the absence of potential $(\Phi = 0)$, we show in section \ref{sec5} the existence of an invariant measure for the Markov transition kernels associated with the electroconvection model \eqref{intro1}--\eqref{intro2} based on the Krylov-Bogoliubov averaging procedure under low regularity assumptions imposed on the noises. In section \ref{solreg}, we prove that any invariant measure of \eqref{intro1}--\eqref{intro2} is smooth provided that the model is forced by smooth noises. Using asymptotic coupling techniques, we prove in section \ref{sec7} the uniqueness of the invariant measure.  In section \ref{SFP}, we address Feller properties in Sobolev norms when the noise processes are sufficiently regular. This uses uniform bounds for the pathwise solution, and these are presented in Appendix B.

\section{Basic Functional Spaces and Notations}

For $1 \le p \le \infty$, we denote by $L^p(\TT^2)$ the Lebesgue spaces of measurable periodic functions $f$ from $\TT^2$ to $\R$ (or $\RR^2)$ that are $p$-integrable on $\TT^2$, that is
\be 
\|f\|_{L^p} = \left(\int_{\TT^2} \|f\|^p \right)^{1/p} <\infty
\ee if $p \in [1, \infty)$ and 
\be 
\|f\|_{L^{\infty}} = {\mathrm{esssup}}_{\TT^2}  |f| < \infty
\ee if $p = \infty$. The $L^2(\TT^2)$ inner product is denoted by $(\cdot,\cdot)_{L^2}$.

For $s > 0$, we denote by $H^s(\TT^2)$ the Sobolev spaces of measurable periodic functions $f$ from $\TT^2$ to $\R$ (or $\RR^2)$ obeying 
\be 
\|f\|_{H^s}^2 = \sum\limits_{k \in \ZZ^2} (1 + |k|^s)^2|q_k|^2 < \infty.
\ee

For a Banach space $(X, \|\cdot\|_{X})$ and $p,q \in [1,\infty]$, we consider the Lebesgue Banach spaces $L^p(\Omega; L_{loc}^q(0,\infty; X))$ of functions $f$ from $X$ to $\R$ (or $\RR^2)$ satisfying 
\be 
\E \left(\int_{0}^{T} \|f\|_{X}^q dt \right)^{\fr{p}{q}} <\infty
\ee for any $T > 0$, with the usual convention when $p = \infty$ or $q = \infty$. The spaces $L_{loc}^q (0,\infty; X)$ and $L^p(\Omega; C^0(0,\infty; X))$ are defined similarly. Here $C^0(0, \infty; X)$ is the space of functions $f$ with the property that the map
\be 
t \mapsto \|f(t)\|_{X}
\ee is continuous for any $f \in X$.

For $s\in\R$, the fractional Laplacian $\l^s$ applied to a mean zero scalar function $f$ is defined as a Fourier multiplier with symbol $|k|^s$, that is, for $f$ given by
\be
f = \sum\limits_{k \in \mathbb{Z}^2 \setminus \left\{0\right\}} f_k e^{ik \cdot x},
\ee
we have that
\be
\l^s f = \sum\limits_{k \in \mathbb{Z}^2 \setminus \left\{0\right\}} |k|^s f_k e^{ik \cdot x}.
\ee

Finally, the periodic Riesz transforms $R = (R_1, R_2)$ applied to scalar functions $f$ are defined as Fourier multipliers 
\be
(R_jf)_k = i\frac{k_j}{|k|} f_k, \quad k\in\ZZ^2\setminus\{0\},\quad j = 1, 2,
\ee
and they are bounded operators on $L^p(\TT^2)$, $1<p<\infty$. We write $R= \na \l^{-1}$.   

Throughout the paper, $C$ denotes a positive universal constant,  and $C(a,b,c,...)$ denotes a positive constant depending on $a$, $b$, $c$, ...

\section{Existence and Uniqueness of Solutions} \la{EUS}

Let $(\Omega, \mathcal{F}, P)$ be a probability space, $\left\{\mathcal{F}_s \right\}_{s \ge 0}$ be a filtration on $(\Omega, \mathcal{F}, P)$, and $\left\{W_k \right\}_{k \ge 1}$ be a collection of independent, identically distributed, real-valued, standard Brownian motions relative to the filtered probability space.

We consider the stochastic electroconvection model 
\be \begin{cases} \la{stochastic}
\d q + u \cdot \na q dt + \l q dt = \Delta \Phi dt +  \tilde{g} dW 
\\\d u + u \cdot \na u dt - \Delta u dt + \na p dt = - q Rq dt - q \na \Phi dt + fdt +  g dW
\\ \na \cdot u = 0
\end{cases}
\ee on $\mathbb{T}^2$ with initial data $q(x,0) = q_0$ and $u(x,0) = u_0$. The unknowns $q(x, t, w)$, $u(x,t,w) = (u_1(x,t,w), u_2 (x,t,w))$, and $p (x,t,w)$ depend on three different variables: position $x \in \mathbb{T}^2$, time  $t \in [0, \infty)$, and outcome  $w \in \Omega$. The body forces $f$ and the potential $\Phi$ depend only on the position variable $x$. The forces $f$ are smooth, divergence-free and have a zero space average. The potential $\Phi$ is assumed to be smooth.  The noise terms $\tilde{g} dW$ and $gdW$ are given by 
\be 
\tilde{g} dW = \sum\limits_{l=1}^{\infty} \tilde{g}_l (x) dW^l (t)
\ee and 
\be 
g dW = \sum\limits_{l=1}^{\infty} g_l(x) dW^l(t).
\ee We assume that the functions $\tilde{g}_l$ are mean-zero and the vector fields $g_l$ are divergence-free for all $l \in \mathbb{N}$. For $k \ge 0$ and $p > 0$, we denote
\be \la{norm1}
\|g\|_{H^k}^2 = \sum\limits_{l=1}^{\infty} \| g_l \|_{H^k}^2,
\ee
\be \la{norm2}
\|\tilde{g}\|_{H^k}^2 = \sum\limits_{l=1}^{\infty} \| \tilde{g}_l \|_{H^k}^2,
\ee and
\be \la{norm3}
\|\tilde{g}\|_{L^p}^p = \int_{\TT^2} \left(\sum\limits_{l=1}^{\infty} |\tilde{g}_l(x)|^2 \right)^{\fr{p}{2}} dx,
\ee
and $g \in H^k$, $\tilde{g} \in H^k$, or $\tilde{g} \in L^p$ if the quantities \eqref{norm1}, \eqref{norm2}, or \eqref{norm3} are finite respectively.

In this section, we prove the existence and uniqueness of solutions of the stochastic model \eqref{stochastic}:

\beg{thm} \la{thm2} Fix a stochastic basis $(\Omega, \mathcal{F}, \mathbb{P}, \left\{\mathcal{F}_t \right\}_{t \ge 0}, W)$. Let $q_0 \in L^4$ have mean zero over $\TT^2$, and let $u_0 \in H^1$ be divergence-free. Suppose $\tilde{g} \in L^4$, $g \in H^1$, $f \in L^2$, and $\Delta \Phi \in L^4$. 
Then there exists a unique pair $(q, u)$ such that $q$ is mean-free, $u$ is divergence-free,  
\be 
u \in L^{2} (\Omega; C^0(0,\infty; L^2) \cap L_{loc}^{\infty}(0,\infty ; H^1) \cap  L_{loc}^2(0,\infty;H^2)),
\ee
\be 
q \in L^{2} (\Omega; C^0(0,\infty; H^{-\fr{1}{2}})  \cap L_{loc}^2(0,\infty; H^{\fr{1}{2}})) \cap L^4(\Omega; L_{loc}^{\infty} (0,\infty; L^4)).
\ee 
Moreover, the elements $(q,u)$ are $\mathcal{F}_t$ adapted and obey
\be 
\d(q, \xi)_{L^2}  + (u \cdot \na q,\xi)_{L^2} dt + (\l q,\xi)_{L^2} dt = (\Delta \Phi,\xi)_{L^2} dt +  (\tilde{g}, \xi)_{L^2} dW
\ee for any $\xi \in H^1$ and a.e. $w \in \Omega$, and 
\be 
\d(u,v)_{L^2} + (u \cdot \na u + qRq, v)_{L^2} dt - (\Delta u,v)_{L^2} dt = (- q \na \Phi,v)_{L^2} dt + (f,v)_{L^2} dt +  (g,v)_{L^2} dW
\ee  for any $v \in H^1$ and a.e. $w \in \Omega$. 
\end{thm}


For each $\epsilon \in (0,1]$, we let $J_{\epsilon}$ be the standard mollifier operator and we consider the viscous approximation
\be 
\begin{cases} \la{nonlinear}
\d q^{\epsilon} + u^{\epsilon} \cdot \na q^{\epsilon} dt + \l q^{\epsilon} dt -\epsilon \Delta q^{\epsilon} dt
= \Delta \Phi dt +   J_{\epsilon} \tilde{g} dW
\\ \d u^{\epsilon} +  u^{\epsilon} \cdot \na u^{\epsilon} dt - \Delta u^{\epsilon} dt + \na p^{\epsilon} dt
= - q^{\epsilon} Rq^{\epsilon}  dt - q^{\epsilon} \na \Phi dt + f dt 
+  J_{\epsilon} g dW
\\\na \cdot u^{\epsilon} = 0
\end{cases}
\ee
with smoothed out initial data $q_0^{\epsilon} = J_{\epsilon}q_0, u_0^{\epsilon} = J_{\epsilon}u_0.$ For each $\epsilon \in (0,1]$, the viscous system \eqref{nonlinear} is forced by smooth noise processes and has local smooth solutions, a fact that can be shown using a fixed point iteration technique. These local solutions extend to global smooth solutions as they remain uniformly bounded in all Sobolev norms, a result that follows from energy-type arguments (see for instance Appendix \ref{B}). In Proposition~\ref{prop1} below, we establish bounds, uniform in time and $\epsilon$, for the solutions of \eqref{nonlinear} in Lebesgue spaces. These estimates are needed to apply the drift identification argument of \cite{MS} and prove Theorem \ref{thm2}.

\beg{prop} \la{prop1} Let $q_0 \in L^4$ have mean zero over $\TT^2$. Let $u_0 \in H^1$ be divergence-free. Suppose $\tilde{g} \in L^4$ and $g \in H^1$. Then the solution $(q^{\epsilon}, u^{\epsilon})$ of \eqref{nonlinear} satisfies 
\beg{align} \la{SS2}
\E \left(\sup\limits_{0 \le t \le T} \|q^{\epsilon}\|_{L^2}^p\right) 
&+ \fr{p^2}{2} \E \left(\int_{0}^{T} \|q^{\epsilon}\|_{L^2}^{p-2} \|\l^{\fr{1}{2}} q^{\epsilon}\|_{L^2}^2 ds \right) \nonumber
\\&\quad\quad \le 2p\|q_0\|_{L^2}^p  
+ C(p) \left(\|\Delta \Phi\|_{L^2}^p  + \|\tilde{g}\|_{L^2}^p \right) T
+ C(p) \|\tilde{g}\|_{L^2}^p T^{\fr{p}{2}}
\end{align} for any $p \ge 2$, 
\beg{align}  \la{SS5}
\E \left\{\sup\limits_{0 \le t \le T} \|q^{\epsilon}\|_{L^4}^p \right\}
+ C(p) \left\{\int_{0}^{T} \|q^{\epsilon}\|_{L^4}^p \right\}
&\le 2p\|q_0\|_{L^4}^p + C(p) \|\Delta \Phi\|_{L^4}^pT \nonumber\\
&\quad\quad+ C(p) \|\tilde{g}\|_{L^4}^p T 
+ C(p)  \|\tilde{g}\|_{L^4}^p T^{\fr{p}{2}}
\end{align} for any $p \ge 4$, 
\beg{align} \la{SS6} 
&\E \left\{\sup\limits_{0 \le t \le T} \|u^{\epsilon}\|_{L^2}^p \right\} 
+ \E \left\{\int_{0}^{T} \|u^{\epsilon} \|_{L^2}^{p-2} \|\na u^{\epsilon}\|_{L^2}^2  dt \right\}  
\le C(p, \|q_0\|_{L^4}, \|u_0\|_{L^2}, f, \Phi, g, \tilde{g})e^{pT}
\end{align} for any $p \ge 2$, and
\beg{align} \la{SS7} 
\E  \left\{\sup\limits_{0 \le t \le T} \|\na u^{\epsilon}(t)\|_{L^2}^2 \right\}
&+ E\left\{\int_{0}^{T} \|\Delta u^{\epsilon} (s)\|_{L^2}^{2}  ds \right\} \nonumber
\\&\quad \le C(\|\na u_0\|_{L^2}, \|q_0\|_{L^4}) + C(\Phi, f, g, \tilde{g})T + C(\tilde{g})T^2.
\end{align} 
\end{prop}

The proof of Proposition \ref{prop1} is based on several applications of It\^o's lemma and is presented in Appendix A.

\beg{prop} \la{mono} Suppose $f \in L^2$ and $\Delta \Phi \in L^4$. Let 
\be \la{operator}
\mathcal{F}(\xi, v) = (v \cdot \na \xi+ \l \xi - \Delta \Phi, v \cdot \na v - \Delta v + \xi R\xi + \xi \na \Phi - f).
\ee 
Let $q_1 \in L^4, q_2 \in L^2, u_1 \in H^2$ and $u_2 \in H^1$. Then there is a positive universal constant $C_0$ such that 
\beg{align} \la{contidea} 
&(\mathcal{F}(q_1, u_1) - \mathcal{F} (q_2, u_2), (\l^{-1}(q_1 - q_2) , u_1 - u_2))_{L^2} \nonumber
\\&\quad+ C_0 K(\Phi, u_1, q_1) \left(\|\l^{-\fr{1}{2}} (q_1 - q_2) \|_{L^2}^2 + \|u_1 - u_2 \|_{L^2}^2 \right) \ge 0
\end{align} holds, where
\be  \la{rcond}
K(\Phi, u_1, q_1) =  \|\na \Phi\|_{L^{\infty}}^2 + \|\na u_1\|_{L^2}^2 + \|\na u_1\|_{L^2} + \|q_1\|_{L^4}^2 + \|q_1\|_{L^4}^4  +  \|\Delta u_1\|_{L^2}^2.
\ee 
\end{prop}

\textbf{Proof:} We have
\beg{align} \la{mon13}
&(\mathcal{F}(q_1, u_1) - \mathcal{F} (q_2, u_2), (\l^{-1}(q_1 - q_2) , u_1 - u_2))_{L^2} \nonumber
\\&= \int_{\TT^2} (u_1 \cdot \na q_1 - u_2 \cdot \na q_2) \l^{-1} (q_1 - q_2) 
+ \int_{\TT^2} \l (q_1 - q_2) \l^{-1} (q_1 - q_2) \nonumber
\\&\quad+ \int_{\TT^2} (u_1 \cdot \na u_1 - u_2 \cdot \na u_2) \cdot (u_1 - u_2) 
- \int_{\TT^2} \Delta (u_1 - u_2) \cdot (u_1 - u_2) \nonumber
\\&\quad\quad+ \int_{\TT^2} (q_1 Rq_1 - q_2 Rq_2) \cdot (u_1 - u_2)
+ \int_{\TT^2} (q_1 - q_2) \na \Phi \cdot (u_1 - u_2).
\end{align}
Integrating by parts, we have 
\be \la{mon11}
\int_{\TT^2} \l (q_1 - q_2) \l^{-1} (q_1 - q_2) 
-  \int_{\TT^2} \Delta (u_1 - u_2) \cdot (u_1 - u_2)
= \|q_1 - q_2 \|_{L^2}^2 + \|\na (u_1 - u_2)\|_{L^2}^2.
\ee
By H\"older and Young inequalities, we have
\be 
\left|\int_{\TT^2} (q_1 - q_2) \na \Phi \cdot (u_1 - u_2) \right|
\le C\|\na \Phi\|_{L^{\infty}}^2\|u_1 - u_2\|_{L^2}^2 + \fr{1}{4} \|q_1 - q_2\|_{L^2}^2.
\ee We note that
\beg{align}
&\int_{\TT^2} (u_1 \cdot \na u_1 - u_2 \cdot \na u_2) \cdot (u_1 - u_2)  \nonumber
\\&= \int_{\TT^2} ((u_1 - u_2) \cdot \na u_1 ) \cdot (u_1 - u_2) + \int_{\TT^2} (u_2 \cdot \na (u_1 - u_2)) \cdot (u_1 - u_2)  \nonumber
\\&=  \int_{\TT^2} ((u_1 - u_2) \cdot \na u_1 ) \cdot (u_1 - u_2)
\end{align} in view of the divergence-free condition satisfied by $u_2$, and hence 
\beg{align}
&\left| \int_{\TT^2} (u_1 \cdot \na u_1 - u_2 \cdot \na u_2) \cdot (u_1 - u_2) \right|
\le \|\na u_1\|_{L^2} \|u_1 - u_2\|_{L^4}^2 \nonumber
\\&\quad\le C\|\na u_1\|_{L^2} \|u_1 - u_2\|_{L^2} \|\na (u_1 - u_2)\|_{L^2}  + C\|\na u_1\|_{L^2} \|u_1 - u_2\|_{L^2}^2 \nonumber
\\&\quad\le C\left(\|\na u_1\|_{L^2}^2 + \|\na u_1\|_{L^2} \right) \|u_1 - u_2\|_{L^2}^2 
+ \fr{1}{4} \|\na (u_1 - u_2)\|_{L^2}^2
\end{align} where we used Ladyzhenskaya's interpolation inequality applied to $u_1 - u_2$. 
Now, we write
\beg{align}\la{EQ1}
&\int_{\TT^2} (u_1 \cdot \na q_1 - u_2 \cdot \na q_2) \l^{-1} (q_1 - q_2) 
= \int_{\TT^2} ((u_1 - u_2) \cdot \na q_1) \l^{-1} (q_1 - q_2) \nonumber
\\&\quad\quad+ \int_{\TT^2} ((u_2 - u_1) \cdot \na (q_1 - q_2)) \l^{-1} (q_1 - q_2)
+ \int_{\TT^2} (u_1 \cdot \na (q_1 - q_2)) \l^{-1}(q_1 - q_2)
\end{align} and 
\beg{align}\la{EQ2}
&\int_{\TT^2} (q_1 Rq_1 - q_2 Rq_2) \cdot (u_1 - u_2)
= \int_{\TT^2} (q_1 - q_2) R q_1 \cdot (u_1 - u_2) \nonumber
\\&\quad\quad+ \int_{\TT^2} (q_2 - q_1)R(q_1 - q_2) \cdot (u_1 - u_2)
+ \int_{\TT^2} q_1 R (q_1 - q_2) \cdot (u_1 - u_2). 
\end{align} Adding \eqref{EQ1} and \eqref{EQ2}, four terms cancel out, namely
\be 
\int_{\TT^2} ((u_2 - u_1) \cdot \na (q_1 - q_2)) \l^{-1} (q_1 - q_2)
= - \int_{\TT^2} (q_2 - q_1)R(q_1 - q_2) \cdot (u_1 - u_2)
\ee and 
\be 
\int_{\TT^2} ((u_1 - u_2) \cdot \na q_1) \l^{-1} (q_1 - q_2)
= - \int_{\TT^2} q_1 R (q_1 - q_2) \cdot (u_1 - u_2), 
\ee
due to the divergence-free condition satisfied by $u_2 - u_1$. 
We estimate 
\beg{align}
&\left| \int_{\TT^2} (q_1 - q_2) R q_1 \cdot (u_1 - u_2) \right|
\le \|Rq_1\|_{L^4} \|q_1 - q_2\|_{L^2} \|u_1 - u_2\|_{L^4} \nonumber
\\&\quad\le C \|q_1\|_{L^4} \|q_1 - q_2\|_{L^2} \left(\|u_1 - u_2\|_{L^2} + \|u_1 - u_2\|_{L^2}^{\fr{1}{2}} \|\na (u_1 - u_2) \|_{L^2}^{\fr{1}{2}} \right) \nonumber
\\&\quad\le C\left(\|q_1\|_{L^4}^2 + \|q_1\|_{L^4}^4 \right) \|u_1 - u_2\|_{L^2}^2 + \fr{1}{4} \|q_1 - q_2\|_{L^2}^2 + \fr{1}{4} \|\na (u_1 - u_2) \|_{L^2}^{2}
\end{align} using H\"older's inequality, the boundedness of the Riesz transforms in $L^4$, Ladyzhenskaya's inequality, and Young's inequality. In view of the commutator estimate (see \cite[Proposition 3]{AI}) 
\be \la{AIcom}
\|\l^{-\fr{1}{2}} (v \cdot \na \rho) - v \cdot \na \l^{-\fr{1}{2}}\rho \|_{L^2} \le C \|\Delta v\|_{L^2} \|\rho\|_{L^2}
\ee
that holds for any divergence-free $v \in H^2$ and mean-zero $\rho \in L^2$,  we have
\beg{align} \la{mon12}
&\left|\int_{\TT^2} u_1 \cdot \na (q_1 - q_2) \l^{-1}(q_1 - q_2)\right| \nonumber
\\&= \left|\int_{\TT^2} \left[\l^{-\fr{1}{2}} (u_1 \cdot \na (q_1 - q_2)) - u_1 \cdot \na \l^{-\fr{1}{2}}(q_1 - q_2)  \right] \l^{-\fr{1}{2}}(q_1 - q_2)\right| \nonumber
\\&\quad\le C\|\Delta u_1\|_{L^2} \|\l^{-\fr{1}{2}} (q_1 - q_2)\|_{L^2} \|q_1 - q_2\|_{L^2}  \nonumber
\\&\quad \le C\|\Delta u_1\|_{L^2}^2 \|\l^{-\fr{1}{2}} (q_1 - q_2)\|_{L^2}^2 + \fr{1}{4} \|q_1 - q_2\|_{L^2}^2.
\end{align} 
Here we also used that $u_1$ is divergence-free. 
Collecting the bounds \eqref{mon11}--\eqref{mon12} and applying them to \eqref{mon13}, we obtain
\beg{align}\la{EQ3}
&(\mathcal{F}(q_1, u_1) - \mathcal{F} (q_2, u_2), (\l^{-1}(q_1 - q_2) , u_1 - u_2))_{L^2}  \nonumber
\\&\quad+ C_0 K(\Phi, u_1, q_1) \left(\|u_1 - u_2\|_{L^2}^2 + \|\l^{-\fr{1}{2}} (q_1 - q_2) \|_{L^2}^2 \right) \nonumber
\\&\ge \fr{1}{4}\left(\|\na (u_1 - u_2)\|_{L^2}^2 + \|q_1 - q_2\|_{L^2}^2 \right)
\ge 0
\end{align} where $K(\Phi, u_1, q_1)$ is given by \eqref{rcond}. This finishes the proof of Proposition \ref{mono}.

Now, we prove Theorem \ref{thm2}. 

\textbf{Proof of Theorem \ref{thm2}:}  Let 
\be 
\mathcal{F}_1 (q^{\epsilon}, u^{\epsilon}) = u^{\epsilon} \cdot \na q^{\epsilon}
\ee and 
\be 
\mathcal{F}_2 (q^{\epsilon}, u^{\epsilon}) = u^{\epsilon} \cdot \na u^{\epsilon} + q^{\epsilon}Rq^{\epsilon}.
\ee We note that 
\beg{align} 
\|\mathcal{F}_1 \|_{H^{-1}}^2 
&\le \|u^{\epsilon}\|_{L^4}^2 \|q^{\epsilon}\|_{L^4}^2
\le C\left(\|u^{\epsilon}\|_{L^2}^2 + \|u^{\epsilon}\|_{L^2} \|\na u^{\epsilon}\|_{L^2}\right) \|q^{\epsilon}\|_{L^4}^2 \nonumber
\\&\le C \|u^{\epsilon}\|_{L^2}^4 + C\| q^{\epsilon}\|_{L^4}^4 + C \|u^{\epsilon}\|_{L^2}^2 \|\na u^{\epsilon}\|_{L^2}^2
\end{align} using Ladyzhenskaya's interpolation inequality, and 
\beg{align}
\|\mathcal{F}_2\|_{H^{-1}}^2 
&\le \|u^{\epsilon}\|_{L^4}^4 + \|q^{\epsilon}\|_{L^4}^2\|R q^{\epsilon}\|_{L^2}^2
\le  C\|u^{\epsilon}\|_{L^2}^4 + C\|u^{\epsilon}\|_{L^2}^2 \|\na u^{\epsilon}\|_{L^2}^2 + C\|\l^{\fr{1}{2}} q^{\epsilon}\|_{L^2}^2\|q^{\epsilon}\|_{L^2}^2
\end{align} using the boundedness of the Riesz transforms in $L^2$. 
In view of the bounds \eqref{SS2}, \eqref{SS5} with $p=4$, and \eqref{SS6}, we deduce that $\mathcal{F}_1$ and $\mathcal{F}_2$ are uniformly bounded in
\be 
L^2(\Omega; L_{loc}^2(0,\infty; H^{-1}(\TT^2))).
\ee 
Therefore, up to subsequences $\mathcal{F}_1 (q^{\epsilon}, u^{\epsilon})$ and  $\mathcal{F}_2 (q^{\epsilon}, u^{\epsilon})$ converge weakly to some functions $F_1$ and $F_2$, respectively, in 
\be
 L^2(\Omega; L_{loc}^2(0,\infty; H^{-1}(\TT^2))).
\ee 
 Moreover, up to subsequences, $u^{\epsilon}$ converges weakly to some function $u$ in
 \be 
L^{2} (\Omega; L_{loc}^{\infty}(0, \infty; H^1(\TT^2))) \cap L^2(\Omega; L_{loc}^2(0,\infty; H^2(\TT^2))),
\ee  in view of the bound \eqref{SS7}, and $q^{\epsilon}$ converges weakly to some function $q$ in 
\be L^{4} (\Omega; L_{loc}^{\infty} (0,\infty; L^4(\TT^2))) \cap L_{loc}^2(\Omega; L^2(0,\infty; H^{1/2}(\TT^2))),
\ee in view of the bounds \eqref{SS2} with $p=2$ and \eqref{SS5} with $p=4$. 

Now we write the equations satisfied by $(q^{\epsilon}, u^{\epsilon})$ and $(q,u)$ as
\be  \la{thm21}
\d (q^{\epsilon}, u^{\epsilon}) + \mathcal{F}(q^{\epsilon}, u^{\epsilon}) dt  + (0,  \na p^{\epsilon} ) dt= (J_{\epsilon} \tilde{g}, J_{\epsilon}g) dW
\ee 
where $\mathcal{F}$ is as in \eqref{operator}, and 
\be  \la{thm22}
\d (q, u) + \mathcal{F}_0 dt = (\tilde{g}, g ) dW
\ee in $L^2(\Omega; L_{loc}^2(0,\infty; H^{-1}(\TT^2)))$, where
\be 
\mathcal{F}_0 = (F_1 + \l q - \Delta \Phi, F_2 - \Delta u + q \nabla \Phi - f).
\ee We show that for almost every $w \in \Omega$ and $t \in [0,\infty)$, we have 
\be \la{dr3}
\mathcal{F}(q,u) = \mathcal{F}_0
\ee in the sense of distributions. 

We note that 
\be 
(q,u) \in L^{2} (\Omega; C^0(0,\infty; H^{-\fr{1}{2}}(\TT^2)))  \times L^{2} (\Omega; C^0(0,\infty; L^2(\TT^2))) 
\ee and $(\l^{-1}q,u)$ obeys the energy equality 
\beg{align} \la{1111}
&\d \left(\|\l^{-\fr{1}{2}} q\|_{L^2}^2 + \|u\|_{L^2}^2 \right)
+ 2 (\mathcal{F}_0, (\l^{-1}q,u))_{L^2} dt \nonumber
\\&\quad= (\|\l^{-\fr{1}{2}}\tilde{g} \|_{L^2}^2 + \|g\|_{L^2}^2) dt
+ 2((\tilde{g}, g), (\l^{-1}q, u))_{L^2} dW
\end{align} (see \cite{GK}).
For a pair
\be 
(\tilde{q}, \tilde{u}) \in L^4(\Omega; L_{loc}^4(0,\infty; L^4(\TT^2))) \times L^2(\Omega; L_{loc}^2(0,\infty; H^2(\TT^2))),
\ee such that $\tilde{q}$ has mean zero and $\tilde{u}$ is divergence-free, we define 
\be  \la{rcond}
r(t,\tilde{q},\tilde{u}) = C_0 \int_{0}^{t} \left[\|\na \Phi\|_{L^{\infty}}^2 + \|\na \tilde{u}\|_{L^2}^2 + \|\na \tilde{u}\|_{L^2} + \|\tilde{q}\|_{L^4}^2 + \|\tilde{q}\|_{L^4}^4  +  \|\Delta \tilde{u}\|_{L^2}^2  \right] ds
\ee where $C_0$ is the constant in \eqref{contidea}.

In order to show the drift identification claim \eqref{dr3}, it is sufficient to show that
\beg{align} \la{dr4}
\E \left\{\int_{0}^{T} 2e^{-r(t)} (\mathcal{F}(q,u) -  \mathcal{F}_0, (\l^{-1}\Psi_1, \Psi_2))_{L^2} dt \right\} \ge 0
\end{align} for all $(\Psi_1, \Psi_2) \in L^4(\Omega; L_{loc}^4(0,\infty; L^4(\TT^2))) \times L^2(\Omega; L_{loc}^2(0,\infty; H^2(\TT^2)))$ such that $\Psi_1$ has mean zero and $\Psi_2$ is divergence-free. Indeed, \eqref{dr4} implies that 
\be 
\E \left\{\int_{0}^{T} 2e^{-r(t)} \|\mathcal{F}(q,u) -  \mathcal{F}_0 \|_{H^{-1} \times H^{-1}}^2 dt \right\} = 0
\ee
from which we conclude that $\mathcal{F}(q,u) = \mathcal{F}_0$ in $H^{-1} \times H^{-1}$ a.e. on $\Omega \times [0,T]$. Accordingly, we proceed to prove \eqref{dr4}.

Denoting $\d r(t)$ by $\dot{r}(t)$, we have 
\beg{align}
&\d \left[e^{-r(t)} \left(\|\l^{-\fr{1}{2}} q\|_{L^2}^2 + \|u\|_{L^2}^2 \right)  \right] + e^{-r(t)} (2\mathcal{F}_0 + \dot{r} (q, u)  , (\l^{-1}q,u) )_{L^2} dt \nonumber
\\&\quad= e^{-r(t)} \left(\|\l^{-\fr{1}{2}}\tilde{g}\|_{L^2}^2 +  \|g\|_{L^2}^2 \right) 
+  e^{-r(t)} ((\tilde{g}, g), ((\l^{-1}q, u))_{L^2} dW 
\end{align} in view of \eqref{1111}. 
Consequently, and using the analogous It\^o stochastic equation obeyed by $e^{-r(t)} \left(\|\l^{-\fr{1}{2}}q^{\epsilon}\|_{L^2}^2 + \|u^{\epsilon}\|_{L^2}^2 \right)$ and the weak lower semi-continuity, we obtain
\beg{align}
&\E \left\{-\int_{0}^{T} e^{-r(t)}(2\mathcal{F}_0 + \dot{r} (q,u), (\l^{-1}q, u))_{L^2} dt \right\} \nonumber
\\&\quad= \E \left\{e^{-r(T)}  \left(\|\l^{-\fr{1}{2}} q(T)\|_{L^2}^2 + \|u(T)\|_{L^2}^2 \right) -  \left(\|\l^{-\fr{1}{2}} q_0\|_{L^2}^2 + \|u_0\|_{L^2}^2 \right)\right\} \nonumber
\\&\quad\quad+ \E\left\{ - \int_{0}^{T} e^{-r(t)}  \left(\|\l^{-\fr{1}{2}}\tilde{g}\|_{L^2}^2 + \|g\|_{L^2}^2 \right) dt  \right\} \nonumber
\\&\quad\le \liminf_{\epsilon \to 0} \E \left\{e^{-r(T)}  \left(\|\l^{-\fr{1}{2}} q^{\epsilon}(T)\|_{L^2}^2 + \|u^{\epsilon}(T)\|_{L^2}^2 \right)\right\} + \lim_{\epsilon \to 0} \E \left\{-  \left(\|\l^{-\fr{1}{2}} J_{\epsilon}q_0\|_{L^2}^2 + \|J_{\epsilon}u_0\|_{L^2}^2 \right)  \right\} \nonumber
\\&\quad\quad+ \lim_{\epsilon \to 0} \E \left\{ - \int_{0}^{T} e^{-r(t)}  \left(\|\l^{-\fr{1}{2}}J_{\epsilon}\tilde{g}\|_{L^2}^2 + \|J_{\epsilon}g\|_{L^2}^2\right) dt   \right\} \nonumber
\\&\quad= \liminf\limits_{\epsilon \to 0} \E \left\{-\int_{0}^{T} e^{-r(t)}(2\mathcal{F}(q^{\epsilon}, u^{\epsilon}) + \dot{r} (q^{\epsilon},u^{\epsilon}), (\l^{-1}q^{\epsilon}, u^{\epsilon}))_{L^2} dt \right\}, 
\end{align} which implies that
\beg{align} \la{thm23}
&\E \left\{\int_{0}^{T} e^{-r(t)}(2\mathcal{F}_0 + \dot{r} (q,u), (\l^{-1}q, u))_{L^2} dt \right\} \nonumber
\\&\quad\ge \limsup_{\epsilon \to 0} \E \left\{\int_{0}^{T} e^{-r(t)}(2\mathcal{F}(q^{\epsilon}, u^{\epsilon}) + \dot{r} (q^{\epsilon},u^{\epsilon}), (\l^{-1}q^{\epsilon}, u^{\epsilon}))_{L^2} dt \right\}. 
\end{align}
In view of \eqref{contidea}, we have 
\beg{align} \la{thm24}
&\E \left\{\int_{0}^{T} e^{-r(t)} (2\mathcal{F}(\tilde{q}, \tilde{u}) + \dot{r} (\tilde{q}, \tilde{u}), (\l^{-1}\tilde{q}, \tilde{u}) - (\l^{-1}q^{\epsilon}, u^{\epsilon}))_{L^2} dt \right\} \nonumber
\\&\quad\ge \E \left\{\int_{0}^{T} e^{-r(t)} (2\mathcal{F}(q^{\epsilon}, u^{\epsilon}) + \dot{r} (q^{\epsilon}, u^{\epsilon}), (\l^{-1} \tilde{q}, \tilde{u}) - (\l^{-1}q^{\epsilon}, u^{\epsilon}) )_{L^2} dt \right\}
\end{align}  for any $(\tilde{q}, \tilde{u}) \in L^4(\Omega; L_{loc}^4(0,\infty; L^4)) \times L^2(\Omega; L_{loc}^2(0,\infty; H^2))$ such that $\tilde{q}$ has mean zero and $\tilde{u}$ is divergence-free.

Putting \eqref{thm23} and \eqref{thm24} together, we obtain 
\beg{align}
&\E \left\{\int_{0}^{T} e^{-r(t)} (2\mathcal{F} (\tilde{q}, \tilde{u}) + \dot{r} (\tilde{q}, \tilde{u}), (\l^{-1}\tilde{q}, \tilde{u})- (\l^{-1}q, u))_{L^2} dt \right\} \nonumber
\\&\quad= \lim_{\epsilon \to 0} \E\left\{\int_{0}^{T} e^{-r(t)} (2\mathcal{F} (\tilde{q}, \tilde{u}) + \dot{r} (\tilde{q}, \tilde{u}), (\l^{-1}\tilde{q}, \tilde{u})- (\l^{-1}q^{\epsilon}, u^{\epsilon}))_{L^2} dt \right\} \nonumber
\\&\quad\ge \liminf_{\epsilon \to 0} \E \left\{\int_{0}^{T} e^{-r(t)} (2\mathcal{F}(q^{\epsilon}, u^{\epsilon}) + \dot{r} (q^{\epsilon}, u^{\epsilon}), (\l^{-1} \tilde{q}, \tilde{u}) - (\l^{-1}q^{\epsilon}, u^{\epsilon}) )_{L^2} dt \right\} \nonumber
\\&\quad=  \E \left\{\int_{0}^{T} e^{-r(t)} (2\mathcal{F}_0 + \dot{r} (q, u), (\l^{-1} \tilde{q}, \tilde{u}) )_{L^2} dt \right\} \nonumber
\\&\quad\quad- \limsup_{\epsilon \to 0} \E \left\{\int_{0}^{T} e^{-r(t)} (2\mathcal{F}(q^{\epsilon}, u^{\epsilon}) + \dot{r} (q^{\epsilon}, u^{\epsilon}),  (\l^{-1}q^{\epsilon}, u^{\epsilon}) )_{L^2} dt \right\} \nonumber
\\&\quad\ge \E \left\{\int_{0}^{T} e^{-r(t)} (2\mathcal{F}_0 + \dot{r} (q,u), (\l^{-1} \tilde{q}, \tilde{u}) - (\l^{-1}q, u))_{L^2} dt \right\}
\end{align} for any $(\tilde{q}, \tilde{u}) \in L^4(\Omega; L_{loc}^4(0,\infty; L^4)) \times L^2(\Omega; L_{loc}^2(0,\infty; H^2))$ such that $\tilde{q}$ has mean zero and $\tilde{u}$ is divergence-free. 
Letting
\be 
(\tilde{q}, \tilde{u}) = (q, u) + \lambda \Psi
\ee where $\lambda > 0$ and $\Psi = (\Psi_1, \Psi_2) \in L^4(\Omega; L_{loc}^4(0,\infty; L^4)) \times L^2(\Omega; L_{loc}^2(0,\infty; H^2))$, $\Psi_1$ having mean zero and $\Psi_2$ being divergence-free, we obtain 
\beg{align}
&\E \left\{\int_{0}^{T} e^{-r(t)} (2\mathcal{F} ((q,u) + \lambda \Psi) + \dot{r} ((q,u) + \lambda \Psi ), \lambda (\l^{-1} \Psi_1, \Psi_2))_{L^2} dt \right\} \nonumber
\\&\quad\ge \E \left\{\int_{0}^{T} e^{-r(t)} (2\mathcal{F}_0 + \dot{r} (q,u),\lambda (\l^{-1} \Psi_1, \Psi_2))_{L^2} dt \right\}.
\end{align} We divide by $\lambda$, and then take the limit as $\lambda$ goes to zero. We obtain \eqref{dr4} from which we conclude that $\mathcal{F}_0 = \mathcal{F} (q,u)$. 

Uniqueness of solutions is obtained as for the deterministic system \cite[Theorem 2]{AI}. Indeed, if we suppose the existence of two different solutions, and we write the equations obeyed by their difference, then we obtain deterministic equations which are independent of the noise. We omit further details. 

\beg{rem}The existence of unique pathwise solutions can be obtained by setting
\be 
Q^{\epsilon} = q^{\epsilon} - \int_{0}^{t} e^{(t-s)\l} \tilde{g}(x) dW
\ee and 
\be 
U^{\epsilon} = u^{\epsilon} - \int_{0}^{t} e^{-(t-s)\Delta} g(x)dW,
\ee writing the determinitic system obeyed by $(Q^{\epsilon}, U^{\epsilon})$, establishing pointwise in $w$ bounds for $(Q^{\epsilon}, U^{\epsilon})$ in 
\be 
(L_{loc}^{\infty}(0,\infty; L^4) \cap L_{loc}^2(0,\infty; H^{\fr{1}{2}})) \times (L_{loc}^{\infty}(0,\infty; H^1) \cap L_{loc}^2(0,\infty; H^2)),
\ee and passing to the limit using the Aubin-Lions lemma. However, this requires higher regularity assumptions on the noise processes forcing the system (as shown in Proposition \ref{pathneeded} below). Consequently, the identification of drift technique minimizes the regularity conditions imposed on the noises $g$ and $\tilde{g}$.  
\end{rem}

\beg{rem} If the ranges of $\tilde{g}$ and $g$ are infinite countable and their components are time-dependent, then the existence and uniqueness of solutions to the corresponding stochastic electroconvection model are obtained on the time interval $[0,T]$ provided that the following regularity condition
\be 
\int_{0}^{T} \left[ \|\tilde{g}(t)\|_{L^4}^4 + \|g(t)\|_{H^1}^2  \right] dt < \infty
\ee holds.
\end{rem}

\section{Electroconvection Semigroup and Weak Feller Properties} \la{sec3}

We consider the space
\be 
\mathcal{H} = \dot{H}^{-\fr{1}{2}} \times L_{\sigma}^2
\ee consisting of vectors $(\xi, v)$ where $\xi \in H^{-\fr{1}{2}}$ has mean zero and $v \in L^2$ is divergence-free, and we consider the space
\be 
\mathcal{V} = \dot{L}^4 \times H_{\sigma}^1
\ee
consisting of vectors $(\xi, v)$ where $\xi \in L^4$ has mean zero and $v \in H^1$ is divergence-free. We define the norms $\|\cdot\|_{\mathcal{H}}$ and $\|\cdot\|_{\mathcal{V}}$  by 
\be 
\|(\xi, v)\|_{\mathcal{H}}^2 = \|\l^{-\fr{1}{2}} \xi\|_{L^2}^2 + \|v\|_{L^2}^2.
\ee and 
\be 
\|(\xi, v)\|_{\mathcal{V}}^2 = \| \xi\|_{L^4}^2 + \|v\|_{H^1}^2
\ee respectively. 
For a time $t \ge 0$ and a Borel set $A \in \mathcal{B}(\mathcal{V})$,  we define the Markov transition kernels associated with \eqref{stochastic} by 
\be 
P_t ((q_0, u_0), A) = \PP((q,u) (t, (q_0, u_0)) \in A)
\ee where $(q,u)(t, (q_0, u_0))$ denotes the solution of the stochastic model \eqref{stochastic} with initial data $(q_0, u_0$) at time $t$.  

Let $\mathcal{M}_b(V)$ be the collection of bounded real-valued Borel measurable functions on $\mathcal{V}$. For each $t \ge 0$ and $\varphi \in \mathcal{M}_b(V)$, we define the Markovian semigroup (which will also be denoted by $\left\{P_t\right\}_{t \ge 0})$ by 
\be 
P_t \varphi (\cdot) = \E \varphi ((q,u) (t, \cdot)) = \int_{V} \varphi (q_0, u_0) P_t(\cdot, d(q_0,u_0)).
\ee
Let $C_{b}(\mathcal{V}, \|\cdot\|_{\mathcal{H}})$ be the space of continuous bounded real-valued functions on the space $(\mathcal{V}, \|\cdot\|_{\mathcal{H}}) $, and $C_{g}(\mathcal{V}, \|\cdot\|_{\mathcal{H}})$ be the space of real continuous functions $\varphi$ on the space $(\mathcal{V}, \|\cdot\|_{\mathcal{H}}) $, with growth 
\be  \la{grow}
|\varphi(\xi, v)| \le C(1 + \|\l^{-\fr{1}{2}} \xi\|_{L^2}^2 + \|v\|_{L^2}^2).
\ee  We point out that continuity of $\varphi$ on the space $(\mathcal{V}, \|\cdot\|_{\mathcal{H}})$ means that if $(\xi_n ,v_n) \in \mathcal{V}$ converges to $(\xi, v)$ in the norm $\|\cdot\|_{\mathcal{H}}$, then $\varphi(\xi_n, v_n)$ converges to $\varphi(\xi, v)$.
The Markovian semigroup $\left\{P_t \right\}_{t \ge 0}$ has the following weak Feller properties:

\beg{thm} \la{MFeller} The semigroup $\left\{P_t \right\}_{t \ge 0}$ is Markov-Feller on $C_b (\mathcal{V},\|\cdot\|_{\mathcal{H}})$ and $C_{g}(\mathcal{V}, \|\cdot\|_{\mathcal{H}})$, that is if $\varphi \in C_b (\mathcal{V},\|\cdot\|_{\mathcal{H}})$, then $P_t \varphi \in C_b (\mathcal{V},\|\cdot\|_{\mathcal{H}})$ and if $\varphi \in C_{g}(\mathcal{V}, \|\cdot\|_{\mathcal{H}})$, then $P_t \varphi \in C_{g}(\mathcal{V}, \|\cdot\|_{\mathcal{H}}).$
\end{thm}

In the proof of Theorem \ref{MFeller} presented below, we use Propositions \ref{Continuity} and \ref{weakboundd}.

\begin{prop} \la{Continuity} (Continuity) Let $(q_0^1, u_0^1)$ and $(q_0^2, u_0^2)$ be in $\mathcal{V}$. Suppose $\tilde{g} \in L^4$ and $g \in H^1$. Then the corresponding solutions $(q_1, u_1)$ and $(q_2, u_2)$ obey 
\begin{align} \la{conttt}
&\|u_1(t) - u_2(t)\|_{L^2}^2 + \|\l^{-\fr{1}{2}} q_1(t) - \l^{-\fr{1}{2}}q_2(t)\|_{L^2}^2\nonumber\\
&\quad\le \exp \left\{r(t)\right\} \left[\|u_{0}^1 - u_{0}^2\|_{L^2}^2 + \|\l^{-\fr{1}{2}} q_{0}^1 - \l^{-\fr{1}{2}} q_{0}^2 \|_{L^2}^2 \right]
\end{align} with probability 1, where 
\be 
r(t) = C_0 \int_{0}^{t} \left[\|\na \Phi\|_{L^{\infty}}^2 + \|\na u_1\|_{L^2}^2 + \|\na u_1\|_{L^2} + \|q_1\|_{L^4}^2 + \|q_1\|_{L^4}^4  +  \|\Delta u_1\|_{L^2}^2  \right] ds
\ee is well-defined and finite almost surely. 
\end{prop}

\textbf{Proof:} We write the equations obeyed by the differences $q_1 - q_2$ and $u_1 - u_2$, and we take their $L^2$ inner product with $\l^{-1}(q_1-q_2)$ and $u_1 - u_2$ respectively. We add the resulting energy equalities and we obtain 
\beg{align}
&\fr{1}{2} \fr{d}{dt} \left[\|\l^{-\fr{1}{2}} (q_1 - q_2) \|_{L^2}^2 + \|u_1 - u_2 \|_{L^2}^2 \right]  \nonumber
\\&\quad+ (\mathcal{F}(q_1, u_1) - \mathcal{F} (q_2, u_2), (\l^{-1}(q_1 - q_2) , u_1 - u_2))_{L^2} = 0
\end{align} where $\mathcal{F}$ is given by \eqref{operator}.
In view of \eqref{contidea}, we have 
\beg{align}
&\fr{1}{2} \fr{d}{dt} \left[\|\l^{-\fr{1}{2}} (q_1 - q_2) \|_{L^2}^2 + \|u_1 - u_2 \|_{L^2}^2 \right] \nonumber 
\\&\quad- r(t, q_1, u_1) \left[\|\l^{-\fr{1}{2}} (q_1 - q_2) \|_{L^2}^2 + \|u_1 - u_2 \|_{L^2}^2\right]
\le 0
\end{align} where $r(t, q_1, u_1)$ is given by \eqref{rcond}. Multiplying by the integrating factor $e^{-\int_{0}^{t} r(s) ds}$ and integrating in time from $0$ to $t$ give \eqref{conttt}.

\beg{prop} \la{weakboundd} Let $(q_0, u_0) \in \mathcal{V}$.  Suppose $\tilde{g} \in L^4$ and $g \in H^1$. Then the unique solution $(q, u)$ of \eqref{stochastic}  obeys
\be \la{weakbound}
\E \left\{ \sup\limits_{0 \le t \le T} (\|\l^{-\fr{1}{2}} q\|_{L^2}^2 + \|u\|_{L^2}^2) \right\}
\le  \E \left\{\|\l^{-\fr{1}{2}} q_{0}\|_{L^2}^2 + \|u_{0}\|_{L^2}^2 + C(\Phi,f, g, \tilde{g}) \right\}e^{C(\Phi)T}. 
\ee  
\end{prop}

\textbf{Proof:} By It\^o's lemma, we have 
\beg{align} \la{EQ31}
\d \|\l^{-\fr{1}{2}}q\|_{L^2}^2 + 2 \|q\|_{L^2}^2 dt
&= -2 (u \cdot \na q, \l^{-1} q)_{L^2} dt + 2 (\Delta \Phi, \l^{-1}q)_{L^2} dt\nonumber
\\\quad\quad&+ \|\l^{-\fr{1}{2}} \tilde{g} \|_{L^2}^2  dt
+ 2  (\l^{-\fr{1}{2}}\tilde{g}, \l^{-\fr{1}{2}}q )_{L^2} dW
\end{align} and 
\beg{align}\la{EQ32}
\d \|u\|_{L^2}^2 + 2\|\na u \|_{L^2}^2 dt
&= - 2(u \cdot \na u, u)_{L^2} - 2(qRq, u)_{L^2} dt- 2(q\na \Phi, u)_{L^2}dt + 2(f, u)_{L^2}dt \nonumber
\\\quad\quad&+ \|g\|_{L^2}^2 dt + 2  (g, u)_{L^2} dW.
\end{align}
We add the equations \eqref{EQ31} and \eqref{EQ32}. Integrating by parts, we have
\be 
(u \cdot \na q, \l^{-1} q)_{L^2} = - (u \cdot Rq, q)_{L^2} = - (qRq, u)_{L^2},
\ee and using the cancellation 
\be 
(u \cdot \na u, u)_{L^2} = 0, 
\ee we obtain the differential equation 
\beg{align}\la{EQ33}
&\d \left\{\|\l^{-\fr{1}{2}}q\|_{L^2}^2 + \|u\|_{L^2}^2 \right\} + 2(\|q\|_{L^2}^2 + \|\na u\|_{L^2}^2) dt
= 2 (\Delta \Phi, \l^{-1}q)_{L^2} dt - 2(q\na \Phi, u)_{L^2} dt + 2(f, u)_{L^2} dt \nonumber
\\&\quad\quad+  \|\l^{-\fr{1}{2}} \tilde{g} \|_{L^2}^2 dt + \|g\|_{L^2}^2 dt
+ 2  (\l^{-\fr{1}{2}}\tilde{g}, \l^{-\fr{1}{2}}q )_{L^2} dW 
+ 2  (g, u)_{L^2} dW.
\end{align} From \eqref{EQ33}, we arrive at the differential inequality 
\beg{align}
&\d \left\{\|\l^{-\fr{1}{2}}q\|_{L^2}^2 + \|u\|_{L^2}^2 \right\} + (\|q\|_{L^2}^2 + \|\na u\|_{L^2}^2) dt
\le C(\|\l \Phi\|_{L^2}^2 + \|f\|_{L^2}^2) dt + C(\|\na \Phi\|_{L^{\infty}}^2 + 1)\|u\|_{L^2}^2 dt \nonumber
\\&\quad\quad+  \|\l^{-\fr{1}{2}}\tilde{g}\|_{L^2}^2 dt + \|g\|_{L^2}^2 dt
+ 2  (\l^{-\fr{1}{2}}\tilde{g}, \l^{-\fr{1}{2}}q )_{L^2} dW 
+ 2  (g, u)_{L^2} dW.
\end{align}
Letting 
\be 
\rho = \|\na \Phi\|_{L^{\infty}}^2 + 1,
\ee we obtain 
\beg{align}
\d \left\{e^{-C\rho t}(\|\l^{-\fr{1}{2}}q\|_{L^2}^2 + \|u\|_{L^2}^2) \right\} 
&\le C(\|\l \Phi\|_{L^2}^2 + \|f\|_{L^2}^2)e^{-C\rho t} dt   
+  \|\l^{-\fr{1}{2}} \tilde{g} \|_{L^2}^2 dt   + \|g\|_{L^2}^2 dt \nonumber
\\&\quad\quad+ 2  (\l^{-\fr{1}{2}}\tilde{g}, \l^{-\fr{1}{2}}q )_{L^2} dW 
+ 2  (g, u)_{L^2} dW.
\end{align}
Integrating in time from $0$ to $t$, taking the supremum over $[0,T]$, applying the expectation $\E$ in $w$, and using suitable martingale estimates, we obtain \eqref{weakbound}. This completes the proof of Proposition~\ref{weakboundd}. 

Now we prove Theorem \ref{MFeller}:

\textbf{Proof of Theorem \ref{MFeller}:} Fix $\varphi \in C_g(\mathcal{V}, \|\cdot\|_{\mathcal{H}})$. Suppose $(\xi_n, v_n)$ converges to $(\xi, v)$ in $(\mathcal{V}, \|\cdot\|_{\mathcal{H}})$, that is 
\be  \la{semi4}
\|\l^{-\fr{1}{2}} (\xi_n - \xi)\|_{L^2}^2 + \|v_n - v\|_{L^2}^2 \rightarrow 0.
\ee In view of the continuity property given in Proposition~\ref{Continuity}, we have 
\be 
\|q(t, \xi_n) - q(t,\xi)\|_{H^{-\fr{1}{2}}} \rightarrow 0
\ee and 
\be 
\|u(t, v_n) - u(t, v)\|_{L^2} \rightarrow 0.
\ee Since $\varphi$ is continuous on $(\mathcal{V}, \|\cdot\|_{\mathcal{H}})$, we conclude that 
\be 
\varphi ((q,u)(t, (\xi_n, v_n))) \rightarrow \varphi ((q,u)(t, (\xi, v))) 
\ee and hence 
\be 
\E \varphi ((q,u)(t, (\xi_n, v_n)))\rightarrow \E \varphi ((q,u)(t, (\xi, v)))
\ee by the Lebesgue Dominated Convergence Theorem, which can be applied due to the growth condition \eqref{grow}, the bound \eqref{weakbound}, and the convergence \eqref{semi4} yielding the boundedness of the sequence of initial datum $(\xi_n, v_n)$ in the $\mathcal{H}$-norm.  This shows that $\left\{P_t\right\}_{t \ge 0}$ is Feller on $C_{g}(\mathcal{V}, \|\cdot\|_{\mathcal{H}})$. Similarly, $\left\{P_t\right\}_{t \ge 0}$ is Feller on $C_b(\mathcal{V}, \|\cdot\|_{\mathcal{H}})$. This ends the proof of Theorem \ref{MFeller}.

\section{Existence and Regularity of Invariant Measures in the Absence of Potential} \la{sec5}

In this section, we consider the electroconvection system 
\be \begin{cases} \la{invsys}
\d q + u \cdot \na q dt + \l q dt =  \tilde{g} dW
\\\d u + u \cdot \na u dt - \Delta u dt + \na p dt = - q Rq dt + fdt +  g dW
\\ \na \cdot u = 0
\end{cases}
\ee in  $\mathbb{T}^2 \times [0, \infty) \times \Omega$ where the potential $\Phi=0$. We note that the system \eqref{invsys} is in the mean-zero frame: if the initial charge density and velocity are assumed to have a zero spatial average, then the solution $(q, u)$ will have mean zero over $\TT^2$ for all positive times $t \ge 0$.

Let $\dot{L}^p$ and $\dot{H}^s$ be the spaces of $L^p$ and $H^s$ functions with zero spatial averages respectively.
Let $H$ and $V$ be the spaces of $L^2$ and $H^1$ functions that are divergence-free and mean zero respectively.  
Let 
\be 
\dot{\mathcal{H}} = \dot{H}^{-\fr{1}{2}} \times H
\ee and 
\be 
\dot{\mathcal{V}} = \dot{L}^4 \times V
\ee  with 
\be 
\|(q, u)\|_{\dot{\mathcal{H}}}^2 = \|\l^{-\fr{1}{2}} q\|_{L^2}^2 + \|u\|_{L^2}^2
\ee and 
\be 
\|(q, u)\|_{\dot{\mathcal{V}} }^2 = \|q\|_{L^4}^2 + \|\na u\|_{L^2}^2
\ee respectively. 
We note that $\dot{\mathcal{V}}$ is compactly embedded in $\dot{\mathcal{H}}$. We define the operator $\mathcal{A}$ on $\mathcal{D}(\mathcal{A}) = \dot{H}^2 \times (H^2 \cap H)$ by 
\be 
\mathcal{A}(\rho, v) = (-\Delta \rho, - P \Delta v)
\ee
where $P$ is the Leray-Hodge projector.  
There is an orthonormal basis of $L^2 \times  H$ consisting of eigenfunctions $\left\{(e_k, w_k) \right\}_{k=1}^{\infty}$ of $\mathcal A$, such that
\be 
(-\Delta e_k, -P\Delta w_k) = \lambda_k (e_k, w_k)
\ee 
where the sequence of eigenvalues $\left\{\lambda_k\right\}_{k=1}^{\infty}$ of $\mathcal A$  counted with multiplicity is nondecreasing and diverges to $\infty$.
Asymptotically, $\lambda_k \ge ck$ for  $k \ge 1$. Let $P_N$ and $Q_N$ be the orthogonal projections of $\mathcal{\dot{H}}$ onto the space spanned by the first $N$ eigenfunctions of $\mathcal{A}$, $(e_k, w_k)$ corresponding to eigenvalues $\lambda_k$, and its orthogonal complement respectively. We have the inequality 
\be \la{gP}
\|Q_N (\l^{-\fr{1}{2}}\rho, v)\|_{L^2}^2 \le \fr{1}{\sqrt{\lambda_{N+1}}} \|(\rho, \na v)\|_{L^2}^2
\ee  
which holds for all $N\ge 1$.

The Markov transition kernels $\left\{P_t \right\}_{t \ge 0}$ associated with the electroconvection model \eqref{invsys}, 
\be 
P_t((q_0, u_0), A) = \mathbb{P} ((q,u)(t, (q_0,u_0)) \in A),
\ee
are defined on $\dot{\mathcal{V}}$ and are $\dot{\mathcal{H}}$-Feller as shown in Theorem \ref{MFeller}.
Here we establish the existence of invariant measures for the Markov transition kernels $\left\{P_t\right\}_{t \ge 0}$.

\beg{thm} \la{weakinvariance} Suppose that $g \in V$ and $\tilde{g} \in \dot{L}^4$. There exists an invariant measure $\mu$ for the Markov transition kernels associated with \eqref{invsys}. Moreover 
\be \la{measurereg11}
\int_{\dot{\mathcal{V}}} \left[\|\l^{\fr{1}{2}} q \|_{L^2}^2 + \|\Delta u\|_{L^2}^2 \right] d\mu((q,u)) \le C < \infty
\ee  
for any invariant measure $\mu$ of \eqref{invsys}, where $C$ is positive constant depending only on $\|f\|_{L^2}$, $\|g\|_{H^1}$, and $\|\tilde{g}\|_{L^4}$. 
\end{thm}

The proof of Theorem \ref{weakinvariance} uses the following auxiliary propositions and is presented at the end of this section. All the estimates can be done rigorously by taking a viscous system approximating~\eqref{invsys}, deriving the bounds for the mollified solution, and then inheriting them to the solution of~\eqref{invsys} using the lower semi-continuity of the norms. We present formal proofs, omitting the approximation. 

\beg{prop} \la{sec51} Let $q_0 \in \dot{H}^{-\fr{1}{2}}$ and $u_0 \in H$. Suppose $g \in L^2$ and $\tilde{g} \in \dot{H}^{-\fr{1}{2}}$. Then
\be \la{invariant0}
\int_{0}^{t} \E \left[\|q(s)\|_{L^2}^2 + \|\na u(s)\|_{L^2}^2 \right]  ds
\le \|\l^{-\fr{1}{2}}q_0\|_{L^2}^2 + \|u_0\|_{L^2}^2 + \left[\|\l^{-\fr{1}{2}}\tilde{g}\|_{L^2}^2 + \|g\|_{L^2}^2 + \|f\|_{L^2}^2 \right]t
\ee holds for all $t \ge 0$. 
\end{prop}

\textbf{Proof:} The sum of the $H^{-\fr{1}{2}}$ norm of $q$ and $L^2$ norm of $u$ obeys the energy equality
\beg{align}
&\d \left\{\|\l^{-\fr{1}{2}}q\|_{L^2}^2 + \|u\|_{L^2}^2 \right\} + 2(\|q\|_{L^2}^2 + \|\na u\|_{L^2}^2) dt \nonumber
\\&=  2(f, u)_{L^2} dt 
+  \|\l^{-\fr{1}{2}} \tilde{g} \|_{L^2}^2 dt 
+ \|g\|_{L^2}^2 dt
+ 2  (\l^{-\fr{1}{2}}\tilde{g}, \l^{-\fr{1}{2}}q )_{L^2} dW 
+ 2  (g, u)_{L^2} dW
\end{align}
(cf. \eqref{EQ31}--\eqref{EQ33} above) which gives the differential inequality
\beg{align}
&\d \left\{\|\l^{-\fr{1}{2}}q\|_{L^2}^2 + \|u\|_{L^2}^2 \right\} + (\|q\|_{L^2}^2 + \|\na u\|_{L^2}^2) dt \nonumber
\\&\le \|f\|_{L^2}^2 dt 
+  \|\l^{-\fr{1}{2}} \tilde{g} \|_{L^2}^2 dt 
+ \|g\|_{L^2}^2 dt
+ 2  (\l^{-\fr{1}{2}}\tilde{g}, \l^{-\fr{1}{2}}q )_{L^2} dW 
+ 2  (g, u)_{L^2} dW
\end{align} where we used the Poincar\'e inequality to bound $L^2$ norm of the mean-free vector $u$ by the $L^2$ norm of its first order derivative. We integrate in time from $0$ to $t$ and we apply $\E$. We obtain the desired bound \eqref{invariant0}.

\beg{prop} Let $q_0 \in \dot{L}^2$. Suppose $\tilde{g} \in \dot{L}^2$. Then 
\be \la{invariant00}
\int_{0}^{t} \E \|\l^{\fr{1}{2}} q(s)\|_{L^2}^2  ds
\le \|q_0\|_{L^2}^2 + \|\tilde{g}\|_{L^2}^2 t
\ee holds for all $t \ge 0$. 
\end{prop}

\textbf{Proof:} The $L^2$ norm of $q$ evolves according to 
\be 
\d \|q\|_{L^2}^2 + 2\|\l^{\fr{1}{2}} q\|_{L^2}^2
=   \| \tilde{g}\|_{L^2}^2 dt 
+ 2  ( \tilde{g}, q)_{L^2} d W
\ee 
where we used the cancellation $(u \cdot \na q, q)_{L^2} = 0$.  We integrate in time from $0$ to $t$ and we apply $\E$. We obtain \eqref{invariant00}.

\beg{prop} Let $p \ge 4$. Let $q_0 \in \dot{L}^4$. Suppose $\tilde{g} \in \dot{L}^4$. Then 
\be \la{invariant1}
\int_{0}^{t} \E \|q(s)\|_{L^4}^p  ds
\le C(p) \left[\|q_0\|_{L^4}^p + \|\tilde{g}\|_{L^4}^pt \right]
\ee holds for all $t \ge 0$. 
\end{prop}

\textbf{Proof:} The $p$-th power of the $L^4$ norm of $q$ obeys the energy inequality
\be \la{pcharge}
\d \|q\|_{L^4}^p  + \fr{cp}{2} \|q\|_{L^4}^p  
\le C \|\tilde{g} \|_{L^4}^p  dt
+ p \|q\|_{L^4}^{p-4}  (\tilde{g}, q^3)_{L^2} dW.
\ee Integrating in time from $0$ to $t$ and applying $\E$, we obtain the desired bound \eqref{invariant1}.

\beg{prop} Let $u_0 \in V$ and $q_0 \in \dot{L}^4$. Suppose $g \in V$ and $\tilde{g} \in \dot{L}^4$. Then 
\beg{align} \la{invariant2}
&\E \|\na u(t)\|_{L^2}^2 + \E \left\{\int_{0}^{t} \|\Delta u(s)\|_{L^2}^2 ds \right\} \nonumber
\\&\quad\quad\le C \left[\|q_0\|_{L^4}^4 + \|\na u_0\|_{L^2}^2 + \left(\|f\|_{L^2}^2 + \|\na g\|_{L^2}^2 + \|\tilde{g}\|_{L^4}^4 \right)t \right]
\end{align}  holds for all $t \ge 0$. 
\end{prop}

\textbf{Proof:} The $L^2$ norm of $\na u$ obeys
\beg{align} \la{gradvelenergy}
\d \|\na u \|_{L^2}^2  + 2\|\Delta u \|_{L^2}^2
= 2 (q Rq, \Delta u)_{L^2} dt 
-2 (f, \Delta u)_{L^2} dt
+ \|\na g \|_{L^2}^2 dt
- 2 (g, \Delta u)_{L^2} dW.
\end{align} 
Here we used the identity
\be 
(u \cdot \na u, \Delta u)_{L^2} = 0
\ee that holds in the two-dimensional periodic setting on $\TT^2$. In view of the boundedness of the Riesz transforms on $L^4$, we have
\be 
|(qRq, \Delta u)_{L^2}| \le \|q\|_{L^4} \|Rq\|_{L^4} \|\Delta u\|_{L^2}
\le C\|q\|_{L^4}^2 \|\Delta u\|_{L^2}.
\ee  Consequently, an application of Young's inequality yields 
\beg{align}
\d \|\na u \|_{L^2}^2  + \|\Delta u \|_{L^2}^2 dt
\le C\|q\|_{L^4}^4 dt + C\|f\|_{L^2}^2 dt + \|\na g\|_{L^2}^2 dt - 2 (g, \Delta u)_{L^2} dW. 
\end{align} Integrating in time from $0$ to $t$ and applying $\E$, we obtain 
\beg{align} 
\E \|\na u(t)\|_{L^2}^2 
&+ \int_{0}^{t} \E \|\Delta u(s)\|_{L^2}^2 ds
\le \|\na u_0\|_{L^2}^2 \nonumber
\\&\quad + C \left(\|f\|_{L^2}^2 + \|\na g\|_{L^2}^2 \right)t
+ C \E \left\{\int_{0}^{t} \|q(s)\|_{L^4}^4 ds \right\}. 
\end{align}
In view of the bound \eqref{invariant1} applied with $p=4$, we obtain \eqref{invariant2}. 

\beg{prop} \la{prop4} Suppose $g \in V$, $\tilde{g} \in \dot{L}^{4}$, and $f \in \dot{L^2}$. Let 
\be 
\nu_T (A) = \fr{1}{T} \int_{0}^{T} \PP ((q(s), u(s)) \in A ) ds.
\ee Then $\left\{\nu_T \right\}$ is tight in $\dot{\mathcal{H}}$ for $u_0 = q_0 = 0$. 
\end{prop}

\textbf{Proof:} Suppose $u_0 = q_0 = 0$. Let $R > 0$, and let $B_R$ be the ball of radius $R$ in $\dot{L}^2 \times V$ (which is compact in $\dot{\mathcal{H}}$). By Chebyshev's inequality, 
\beg{align} 
\sup\limits_{T > 0} \nu_{T} (B_R^c) 
&= \sup\limits_{T > 0} \fr{1}{T} \int_{0}^{T} \PP (\|(q, u)\|_{\dot{L}^2 \times V} \ge R) dt \nonumber
\\&\le \fr{1}{R^2} \sup\limits_{T > 0} \fr{1}{T} \int_{0}^{T} \E \|(q, u)\|_{\dot{L}^2 \times V}^2) dt \rightarrow 0
\end{align} as $R \rightarrow \infty$ in view of the bound \eqref{invariant0} that is linear in $T$. Therefore, the family $\left\{\nu_T \right\}$ is tight in $\dot{\mathcal{H}}$, ending the proof of Proposition~\ref{prop4}.

Now we prove Theorem \ref{weakinvariance}. 

\textbf{Proof of Theorem \ref{weakinvariance}:} We adapt the notation $w  = (q,u)$ and write solutions as $w (t, w_0)$. From the weak Feller property obtained in Theorem \ref{MFeller}, the tightness of the time-averaged measures obtained in Proposition \ref{prop4}, and the Krylov-Bogoliubov averaging procedure, we conclude that there exists a probability measure $\mu$ satisfying
\be 
\int_{\dot{\mathcal{H}}} \varphi (w_0) d\mu(w_0) = \int_{\dot{\mathcal{H}}} \int_{\dot{\mathcal{H}}} \fr{1}{T} \int_{0}^{T} P_t (w_0, dw) \varphi (w) dt d\mu(w_0)
\ee for any $T > 0$ and any $\varphi \in C_b(\mathcal{\dot{H}})$. Now we study the regularity of $\mu$ and we prove \eqref{measurereg11}. 
For $n \ge 1$, we let $P_n$ be the projection onto the space spanned by the first $n$ eigenfunctions of $-\Delta$. 
For $n \ge 1, R>0, w = (q,u) \in \dot{\mathcal{H}}$, we let
\be 
\Psi_{n,R} (w) = \left[\|P_n q \|_{L^2}^2 + \|\na P_n u\|_{L^2}^2 \right] \wedge R
\ee and we note that $\Psi_{n,R} \in C_b(\dot{\mathcal{H}})$. In view of \eqref{invariant0}, we estimate
\beg{align} 
&\left|\fr{1}{T} \int_{0}^{T} \int_{\dot{\mathcal{H}}} P_t(w_0, dw) \Psi_{n,R}(w) dt \right|
= \left|\fr{1}{T} \int_{0}^{T} \E \Psi_{n,R}(w(t, w_0)) dt \right| \nonumber
\\&\le \fr{1}{T} \E \int_{0}^{T} \left[\|q\|_{L^2}^2 + \|\na u \|_{L^2}^2 \right] \le \left(\|\l^{-\fr{1}{2}}q_0\|_{L^2}^2 + \|u_0\|_{L^2}^2\right) T^{-1} + \|\l^{-\fr{1}{2}}\tilde{g}\|_{L^2}^2 + \|g\|_{L^2}^2 + \|f\|_{L^2}^2
\end{align} 
for any $T > 0$. Let $B_{\dot{\mathcal{H}}}(\rho)$ be the ball
\be 
B_{\dot{\mathcal{H}}}(\rho) = \left\{w \in \dot{\mathcal{H}}: \|w\|_{\dot{\mathcal{H}}}^2 \le \rho^2 \right\}.
\ee Then, using invariance, we have
\beg{align}
&\int_{\dot{\mathcal{H}}} \Psi_{n,R} (w_0) d\mu(w_0)
\le \int_{B_{\dot{\mathcal{H}}}(\rho)} \left|\fr{1}{T} \int_{0}^{T} \int_{\dot{\mathcal{H}}} P_t(w_0, dw) \Psi_{n,R}(w) dt \right| d\mu(w_0) \nonumber
\\&\quad+ \int_{\dot{\mathcal{V}} \setminus B_{\dot{\mathcal{H}}}(\rho)} \left|\fr{1}{T} \int_{0}^{T} \int_{\dot{\mathcal{H}}} P_t(w_0, dw) \Psi_{n,R}(w) dt \right| d\mu(w_0) \nonumber
\\&\le \left[\rho^2 T^{-1}+ \|\l^{-\fr{1}{2}}\tilde{g}\|_{L^2}^2 + \|g\|_{L^2}^2 + \|f\|_{L^2}^2 \right] \mu(B_{\dot{\mathcal{H}}}(\rho)) + R\mu(\dot{\mathcal{H}} \setminus B_{\dot{\mathcal{H}}}(\rho)).
\end{align}
We choose $\rho$ large enough so that 
\be 
R\mu(\dot{\mathcal{H}} \setminus B_{\dot{\mathcal{H}}}(\rho)) \le 1
\ee and then we choose $T$ large enough so that 
\be 
\rho^2 T^{-1} \le 1
\ee and we get
\be 
\int_{\dot{\mathcal{H}}} \Psi_{n,R} (w_0) d\mu(w_0) \le  2 +  \|\l^{-\fr{1}{2}}\tilde{g}\|_{L^2}^2 + \|g\|_{L^2}^2 + \|f\|_{L^2}^2 .
\ee 
By Fatou's lemma, we have
\be 
\int_{\dot{\mathcal{H}}} \left\{\left[\|q_0 \|_{L^2}^2 + \|\na u_0\|_{L^2}^2 \right] \wedge R \right\} d\mu (w_0) \le  2 +  \|\l^{-\fr{1}{2}}\tilde{g}\|_{L^2}^2 + \|g\|_{L^2}^2 + \|f\|_{L^2}^2 
\ee 
and by the Monotone Convergence Theorem, we obtain
\be 
\int_{\dot{\mathcal{H}}} \left[\|q_0 \|_{L^2}^2 + \|\na u_0\|_{L^2}^2 \right] d\mu (w_0) \le   2 +  \|\l^{-\fr{1}{2}}\tilde{g}\|_{L^2}^2 + \|g\|_{L^2}^2 + \|f\|_{L^2}^2 .
\ee 
Therefore, the invariant measure $\mu$ is supported on $\mathcal{X}_2 = \dot{L}^2 \times V$.
Next we upgrade the regularity of the measure $\mu$. For $w = (q, u) \in \mathcal{X}_2$, we define
\be 
\Psi_{n,R}^2 (w) = \left[\|\l^{\fr{1}{2}} P_n q \|_{L^2}^2 +  \|\na P_n u \|_{L^2}^2 \right] \wedge R.
\ee 
 In view of the bounds \eqref{invariant0} and \eqref{invariant00}, we have
\beg{align} \la{000}
&\left|\fr{1}{T} \int_{0}^{T} \E \Psi_{n,R}^2(w(t, w_0)) dt \right| \nonumber
\le \fr{1}{T} \E \int_{0}^{T} \left[\|\l^{\fr{1}{2}}q\|_{L^2}^2 + \|\na u \|_{L^2}^2 \right] dt
\\&\le \left(2\|q_0\|_{L^2}^2 +  \|u_0\|_{L^2}^2 \right)T^{-1} + 2\|\tilde{g}\|_{L^2}^2 + \|g\|_{L^2}^2 + \|f\|_{L^2}^2
\end{align} for any $T > 0$. Letting $B_{\mathcal{X}_2}(\rho)$ be the ball
\be 
B_{\mathcal{X}_2}(\rho) = \left\{w = (q,u) \in \mathcal{X}_2: \|q\|_{L^2}^2 + \|\na u\|_{L^2}^2 \le \rho^2 \right\}, 
\ee 
we use \eqref{000} and invariance to obtain 
\beg{align}
&\int_{\mathcal{X}_2} \Psi_{n,R}^2 (w_0) d\mu(w_0) \nonumber
= \int_{\mathcal{X}_2} \fr{1}{T} \int_{0}^{T} \E \Psi_{n,R}^2(w(t, w_0)) dt d\mu(w_0)
\\&\le \left[2 \rho^2 T^{-1} + 2\|\tilde{g}\|_{L^2}^2 + \|g\|_{L^2}^2 + \|f\|_{L^2}^2 \right] \mu(B_{\mathcal{X}_2}(\rho)) + R\mu(\mathcal{X}_2 \setminus B_{\mathcal{X}_2}(\rho)).
\end{align} 
We choose $\rho$ large enough and $T$ large enough so that 
\be 
\int_{\mathcal{X}_2} \Psi_{n,R}^2 (w_0) d\mu(w_0) \le 2 + 2\|\tilde{g}\|_{L^2}^2 + \|g\|_{L^2}^2 + \|f\|_{L^2}^2.
\ee 
By Fatou's lemma and the Monotone Convergence Theorem, we obtain
\be 
\int_{\mathcal{X}_2} \left[\|\l^{\fr{1}{2}}q_0 \|_{L^2}^2 + \|\Delta u_0\|_{L^2}^2 \right] d\mu (w_0) \le  2 + 2\|\tilde{g}\|_{L^2}^2 + \|g\|_{L^2}^2 + \|f\|_{L^2}^2.
\ee Therefore, the invariant measure $\mu$ is supported on $\mathcal{X}_3 = \dot{H}^{\fr{1}{2}} \times  V$.  
Finally, for $w = (q, u) \in \mathcal{X}_3$, we define
\be 
\Psi_{n,R}^3 (w) = \left[\|\l^{\fr{1}{2}} P_n q \|_{L^2}^2 +  \|\Delta P_n u \|_{L^2}^2 \right] \wedge R.
\ee 
 In view of the bounds \eqref{invariant00} and \eqref{invariant2}, we have
\beg{align} \la{0000}
&\left|\fr{1}{T} \int_{0}^{T} \E \Psi_{n,R}^3(w(t, w_0)) dt \right| \nonumber
\le \fr{1}{T} \E \int_{0}^{T} \left[\|\l^{\fr{1}{2}}q\|_{L^2}^2 + \|\Delta u \|_{L^2}^2 \right] dt
\\&\le \left(\|q_0\|_{L^2}^2 + C\|q_0\|_{L^4}^4 + C\|\na u_0\|_{L^2}^2 \right)T^{-1} + \|\tilde{g}\|_{L^2}^2 + C\|f\|_{L^2}^2 + C\|\na g\|_{L^2}^2 + C\|\tilde{g}\|_{L^4}^4
\end{align} for any $T > 0$. We let $B_{\mathcal{X}_3}(\rho)$ be the ball
\be 
B_{\mathcal{X}_3}(\rho) = \left\{w = (q,u) \in \mathcal{X}_3: \|\l^{\fr{1}{2}} q\|_{L^2}^2 + \|\na u\|_{L^2}^2 \le \rho^2 \right\}. 
\ee  Using the bound \eqref{000}, invariance, and the continuous embedding of $H^{\fr{1}{2}}$ in $L^4$, we obtain
\beg{align}
&\int_{\mathcal{X}_3} \Psi_{n,R}^3 (w_0) d\mu(w_0) \nonumber
= \int_{\mathcal{X}_3} \fr{1}{T} \int_{0}^{T} \E \Psi_{n,R}^3(w(t, w_0)) dt d\mu(w_0)
\\&\le C\left[(\rho^2 + \rho^4)T^{-1} + \|\tilde{g}\|_{L^2}^2 + \|f\|_{L^2}^2 + \|\na g\|_{L^2}^2 + \|\tilde{g}\|_{L^4}^4 \right] \mu(B_{\mathcal{X}_3}(\rho)) + R\mu(\mathcal{X}_3 \setminus B_{\mathcal{X}_3}(\rho)).
\end{align} 
We choose $\rho$ large enough and $T$ large enough so that
\be 
\int_{\mathcal{X}_3} \left[\|\l^{\fr{1}{2}}q_0 \|_{L^2}^2 + \|\Delta u_0\|_{L^2}^2 \right] d\mu (w_0) \le   C\left(1 + \|\tilde{g}\|_{L^2}^2 + \|f\|_{L^2}^2 + \|\na g\|_{L^2}^2 + \|\tilde{g}\|_{L^4}^4 \right).
\ee Therefore, the invariant measure $\mu$ is supported on $\dot{H}^{\fr{1}{2}} \times  (H^2 \cap V)$.
This ends the proof of Theorem~\ref{weakinvariance}.

\section{Higher Regularity of Invariant Measures} \la{solreg}

In this section, we prove that any invariant measure of \eqref{invsys} is more regular than $\dot{H}^{\fr{1}{2}} \times (H^2 \cap V)$.

\beg{thm} \la{higherreg}
Suppose $g$ and $\tilde{g}$ are smooth. If $\mu$ is an invariant measure of \eqref{invsys}, then $\mu$ is smooth and satisfies
\be 
\int_{\dot{\mathcal{V}}} \log \left[1 + \|u\|_{H^k}^2 + \|q\|_{H^k}^2 \right] d\mu((q,u))  \le C(k, f, g, \tilde{g}) < \infty.
\ee for any $k \ge 0$. 
\end{thm}

The proof of Theorem \ref{higherreg} is based on the following auxilliary propositions 
 and is presented at the end of this section.

\beg{prop} \la{sec61} Let $u_0 \in V$ and $q_0 \in \dot{L}^4$. Suppose $g \in V$ and $\tilde{g} \in \dot{L}^4$. Let $p \ge 4$. Then
\beg{align} \la{invariant3}
&\E \left\{\int_{0}^{t} \|\na u(s)\|_{L^2}^{p-2} \|\Delta u(s)\|_{L^2}^2 ds \right\} \nonumber
\\&\quad\quad\le C(p)\left[\|q_0\|_{L^4}^{2p} + \|\na u_0\|_{L^2}^{p} + \left(\|f\|_{L^2}^{p} + \|\na g\|_{L^2}^{p} + \|\tilde{g}\|_{L^4}^{2p} \right) t\right]
\end{align}  holds for all $t \ge 0$. 
\end{prop}

\textbf{Proof:} The $L^2$ norm of $\na u$ evolves according to the stochastic energy equality
\be 
\d \|\na u\|_{L^2}^{2} + 2 \|\Delta u\|_{L^2}^2 dt
= 2(qRq - f, \Delta u)_{L^2} dt
+ \|\na g\|_{L^2}^2 dt
+ 2 (\na g, \na u)_{L^2} dW.
\ee Consequently, the $p$-th power of $\|\na u\|_{L^2}$ obeys
\beg{align}
&\d \|\na u\|_{L^2}^p 
+ p\|\na u\|_{L^2}^{p-2}\|\Delta u\|_{L^2}^2 dt
= p \|\na u\|_{L^2}^{p-2} (qRq - f, \Delta u)_{L^2} dt
+ \fr{p}{2} \|\na g\|_{L^2}^2 \|\na u\|_{L^2}^{p-2} dt \nonumber
\\&\quad\quad+ p \left(\fr{p}{2} - 1 \right) \|\na u\|_{L^2}^{p-4} (\na g, \na u)_{L^2}^2 dt
+ p\|\na u\|_{L^2}^{p-2} (\na g, \na u)_{L^2} dW \nonumber
\\&\quad\le \fr{p}{4} \|\na u\|_{L^2}^{p-2} \|\Delta u\|_{L^2}^2 dt 
+ \fr{p}{4} \|\na u\|_{L^2}^p dt
+ C(p) \|q\|_{L^4}^{2p} dt
+ C(p) \left[\|f\|_{L^2}^p + \|\na g\|_{L^2}^p\right]dt \nonumber
\\&\quad\quad\quad+ p\|\na u\|_{L^2}^{p-2} (\na g, \na u)_{L^2} dW.
\end{align} In view of the Poincar\'e inequality, we obtain 
\beg{align} 
&\d \|\na u\|_{L^2}^p 
+ \fr{p}{2} \|\na u\|_{L^2}^{p-2}\|\Delta u\|_{L^2}^2 dt \nonumber
\\&\quad\le C(p) \|q\|_{L^4}^{2p} dt
+ C(p) \left[\|f\|_{L^2}^p + \|\na g\|_{L^2}^p\right]dt 
+ p\|\na u\|_{L^2}^{p-2} (\na g, \na u)_{L^2} dW.
\end{align}
We integrate in time from $0$ to $t$ and we apply $\E$. In view of the bound \eqref{invariant1}, we obtain \eqref{invariant3}.

\beg{prop} Let $u_0 \in V$ and $q_0 \in \dot{L}^4$. Suppose $g \in V$ and $\tilde{g} \in \dot{L}^4$. Then 
\beg{align} \la{invariant4}
&\E \left\{\int_{0}^{t} \|\na u(s)\|_{L^2}^2 \|\Delta u(s)\|_{L^2}^2 \|q(s)\|_{L^4}^4 ds \right\} \nonumber
\\&\quad\le C(f,g,\tilde{g}) \left[\|q_0\|_{L^4}^{12} + \|q_0\|_{L^4}^8 + \|q_0\|_{L^4}^4 + \|\na u_0\|_{L^2}^4 + \|q_0\|_{L^4}^4\|\na u_0\|_{L^2}^4 + t\right]
\end{align} holds for all $t \ge 0$.
\end{prop}

\textbf{Proof:} The stochastic process $\|\na u\|_{L^2}^4 \|q\|_{L^4}^4$ obeys
\be 
\d \left[\|\na u\|_{L^2}^4 \|q\|_{L^4}^4 \right]
= \|\na u\|_{L^2}^4 \d \|q\|_{L^4}^4 
+ \|q\|_{L^4}^4 \d \|\na u\|_{L^2}^4
+ \d \|\na u\|_{L^2}^4 \cdot \d \|q\|_{L^4}^4.
\ee The $4$-th power of the $L^2$ norm of $\na u$ evolves according to 
\beg{align} \la{evolprod}
&\d \|\na u\|_{L^2}^4
= -4 \|\na u \|_{L^2}^2 \|\Delta u\|_{L^2}^2 dt
+ 4\|\na u\|_{L^2}^2 (qRq  - f, \Delta u)_{L^2} dt \nonumber
\\&\quad\quad+2 \|\na u\|_{L^2}^2 \|\na g\|_{L^2} ^2dt
+ 4  |(g, \Delta u)_{L^2}|^2 dt
-4 \|\na u\|_{L^2}^2  (g, \Delta u)_{L^2} dW
\end{align} whereas the $4$-th power of the $L^4$ norm of $q$ evolves according to 
\beg{align}  \la{invariant551}
\d \|q\|_{L^4}^4 
= - 4 (\l q, q^3)_{L^2} dt
+ 6   (\tilde{g}^2, q^2)_{L^2} dt
+ 4  (\tilde{g}, q^3)_{L^2} dW.
\end{align}
Consequently, the product $\|\na u\|_{L^2}^4 \|q\|_{L^4}^4$ satisfies the energy equality
\beg{align}
&\d \left[\|q\|_{L^4}^4 \|\na u\|_{L^2}^4 \right]
= -4 \|\na u\|_{L^2}^4 (\l q, q^3)_{L^2} dt
+ 6\|\na u\|_{L^2}^4  (\tilde{g}^2, q^2)_{L^2} dt \nonumber
\\&\quad + 4\|\na u\|_{L^2}^4  (\tilde{g}, q^3)_{L^2} dW 
-4\|q\|_{L^4}^4 \|\na u\|_{L^2}^2\|\Delta u\|_{L^2}^2 dt 
+ 4\|q\|_{L^4}^4 \|\na u\|_{L^2}^2 (qRq - f, \Delta u)_{L^2} dt \nonumber
\\&\quad\quad +2\|q\|_{L^4}^4 \|\na u\|_{L^2}^2 \|\na g\|_{L^2}^2 dt
+ 4\|q\|_{L^4}^4  (g, \Delta u)_{L^2}^2 dt
- 4\|q\|_{L^4}^4 \|\na u\|_{L^2}^2  (g, \Delta u)_{L^2} dW \nonumber
\\&\quad\quad\quad- 16\|\na u\|_{L^2}^2  (\tilde{g}, q^3)_{L^2} (g, \Delta u)_{L^2} dt
\end{align} which yields the energy inequality
\beg{align} \la{proof1}
&\d \left[\|q\|_{L^4}^4 \|\na u\|_{L^2}^4 \right]
+4c \|\na u\|_{L^2}^4 \|q\|_{L^4}^4 dt
+4\|q\|_{L^4}^4 \|\na u\|_{L^2}^2\|\Delta u\|_{L^2}^2 dt \nonumber
\\&\quad \le 6\|\na u\|_{L^2}^4  (\tilde{g}^2, q^2)_{L^2} dt
+ 4\|q\|_{L^4}^4 \|\na u\|_{L^2}^2 (qRq - f, \Delta u)_{L^2} dt \nonumber
\\&\quad\quad +2\|q\|_{L^4}^4 \|\na u\|_{L^2}^2 \|\na g\|_{L^2}^2 dt
+ 4\|q\|_{L^4}^4  (\na g, \na u)_{L^2}^2 dt
- 16\|\na u\|_{L^2}^2  (\tilde{g}, q^3)_{L^2} (g, \Delta u)_{L^2} dt \nonumber
\\&\quad\quad\quad - 4\|q\|_{L^4}^4 \|\na u\|_{L^2}^2  (g, \Delta u)_{L^2} dW 
+ 4\|\na u\|_{L^2}^4  (\tilde{g}, q^3)_{L^2} dW.
\end{align} 
Here, we used the nonlinear Poincar\'e inequality for the fractional Laplacian in $L^4$ applied to the mean zero function $q$ (see \cite{AI,CGV})
\be  \la{3}
\int_{\TT^2} q^3 \l q dx \ge c \|q\|_{L^4}^4.
\ee 
By the Cauchy-Schwartz inequality, Young's inequality and the Poincar\'e inequality applied to the mean zero function $\na u$, we estimate
\beg{align} 
&\left|6\|\na u\|_{L^2}^4  (\tilde{g}, q^2)_{L^2} \right|
\le 6\|\na u\|_{L^2}^4 \|\tilde{g}\|_{L^4}^2  \|q\|_{L^4}^2  \nonumber
\\&\quad \le \fr{c}{8} \|\na u\|_{L^2}^4 \|q\|_{L^4}^4 
+ C \|\tilde{g}\|_{L^4}^4 \|\na u\|_{L^2}^2 \|\Delta u\|_{L^2}^2.
\end{align}
The boundedness of the Riesz transforms on $L^4$ yields
\beg{align}
&\left|4\|q\|_{L^4}^4 \|\na u\|_{L^2}^2 (qRq - f, \Delta u)_{L^2} \right|
\le C\|q\|_{L^4}^6 \|\na u\|_{L^2}^2 \|\Delta u\|_{L^2} 
+ C\|q\|_{L^4}^4 \|\na u\|_{L^2}^2 \|\Delta u\|_{L^2} \|f\|_{L^2} \nonumber
\\&\quad \le \fr{1}{8} \|q\|_{L^4}^4\|\na u\|_{L^2}^2 \|\Delta u\|_{L^2}^2 
+ \fr{c}{8} \|q\|_{L^4}^4 \|\na u\|_{L^2}^4 
+ C\|q\|_{L^4}^{12}
+ C\|q\|_{L^4}^4 \|f\|_{L^2}^4. 
\end{align}
We bound 
\be 
2\|q\|_{L^4}^4 \|\na u\|_{L^2}^2 \|\na g\|_{L^2}^2 
\le \fr{c}{8} \|q\|_{L^4}^4 \|\na u\|_{L^2}^4 
+ C\|\na g\|_{L^2}^4 \|q\|_{L^4}^4
\ee and 
\be 
4\|q\|_{L^4}^4  (\na g, \na u)_{L^2}^2 
\le 4\|q\|_{L^4}^4\|\na u\|_{L^2}^2 \|\na g\|_{L^2}^2
\le \fr{c}{8} \|q\|_{L^4}^4 \|\na u\|_{L^2}^4 
+ C\|\na g\|_{L^2}^4 \|q\|_{L^4}^4
\ee using Young's inequality. 
Finally, we estimate
\beg{align} \la{proof2}
&\left|16\|\na u\|_{L^2}^2  (\tilde{g}, q^3)_{L^2} (g, \Delta u)_{L^2}\right|
\le 16\|\na u\|_{L^2}^3 \|q\|_{L^4}^3  \|\tilde{g}\|_{L^4} \|\na g\|_{L^2} \nonumber
\\&\quad \le \fr{c}{8} \|\na u\|_{L^2}^4 \|q\|_{L^4}^4 
+ C\|\tilde{g}\|_{L^4}^4 \left( \|\na g\|_{L^2} \right)^4.
\end{align}
Putting \eqref{proof1}--\eqref{proof2} together, we obtain the differential inequality 
\beg{align}
&\d \left[\|q\|_{L^4}^4 \|\na u\|_{L^2}^4 \right]
+ c\|\na u\|_{L^2}^4 \|q\|_{L^4}^4 dt
+ \|q\|_{L^4}^4 \|\na u\|_{L^2}^2 \|\Delta u\|_{L^2}^2 dt \nonumber 
\\&\quad \le C(\tilde{g}) \|\na u\|_{L^2}^2 \|\Delta u\|_{L^2}^2 dt
+ C(f, g) \|q\|_{L^4}^4 dt
+ C(g, \tilde{g}) dt
+ C\|q\|_{L^4}^{12} dt \nonumber
\\&\quad\quad - 4\|q\|_{L^4}^4 \|\na u\|_{L^2}^2  (g, \Delta u)_{L^2} dW 
+ 4\|\na u\|_{L^2}^4  (\tilde{g}, q^3)_{L^2} dW.
\end{align}
We integrate in time from $0$ to $t$ and we apply $\E$. The bound \eqref{invariant1} applied with $p=4$ and $p=12$ together with the bound \eqref{invariant3} gives the desired estimate \eqref{invariant4}. 

\beg{prop} Let $u_0 \in V$ and $q_0 \in \dot{H}^{\fr{1}{2}}$. Suppose $g \in V$ and $\tilde{g} \in \dot{H}^{\fr{1}{2}}$. Then 
\beg{align} \la{invariant6}
\E \left\{\int_{0}^{t} \log (1 + \|\na q(s)\|_{L^2}^2) ds \right\}
&\le \log(1 + \|\l^{\fr{1}{2}}q_0\|_{L^2}^2) + \|\na u_0\|_{L^2}^2 + C\|q_0\|_{L^4}^4 + \|q_0\|_{L^2}^2 \nonumber
\\&\quad+ C\left(\|f\|_{L^2}^2 + \|\na g\|_{L^2}^2 + \|\tilde{g}\|_{L^4}^4 + \|\l^{\fr{1}{2}} \tilde{g}\|_{L^2}^2 \right)t
\end{align}  holds for all $t \ge 0$. 
\end{prop}

\textbf{Proof:} The $\dot{H}^{\fr{1}{2}}$ norm of $q$ obeys
\be 
\d \|\l^{\fr{1}{2}} q\|_{L^2}^2 
+ 2\|\l q\|_{L^2}^2 dt
= - 2(u \cdot \na q, \l q)_{L^2} dt
+ \|\l^{\fr{1}{2}} \tilde{g}\|_{L^2}^2 dt
+ 2(\tilde{g}, \l q)_{L^2} dW. 
\ee For each $t \ge 0$, let 
\be 
X(t) = \|\l^{\fr{1}{2}}q(t)\|_{L^2}^2
\ee and 
\be
\bar{X}(t) = \|\l q(t)\|_{L^2}^2.
\ee
By It\^o's lemma, we have 
\beg{align} 
&\d \log (1 + X) 
+ \fr{2\bar{X}}{1+X} dt
= -\fr{2}{1+X} (u \cdot \na q, \l q)_{L^2} dt \nonumber
\\&\quad+ \fr{1}{1+X} \|\l^{\fr{1}{2}} \tilde{g}\|_{L^2}^2 dt
- \fr{2}{(1+X)^2} (\tilde{g}, \l q)_{L^2}^2 dt
+ \fr{2}{1+X} (\tilde{g}, \l q)_{L^2} dW.
\end{align}
The nonlinear term is estimated using commutator estimates (see \cite[Proposition 3]{AI})
\beg{align}
\left|\int_{\TT^2} (u \cdot \na q) \l q \right|
= \left|\int_{\TT^2} (\l^{\fr{1}{2}} (u \cdot \na q) - u \cdot \na \l^{\fr{1}{2}}q) \l^{\fr{1}{2}}q  \right|
\le C\|\Delta u\|_{L^2}  \|\l q\|_{L^2} \|\l^{\fr{1}{2}}q\|_{L^2},
\end{align} hence
\be 
\d \log (1+X) + \fr{2\bar{X}}{1+X} dt 
\le \fr{C}{1+X} \|\Delta u\|_{L^2} \sqrt{X} \sqrt{\bar{X}} dt
+ \|\l^{\fr{1}{2}} \tilde{g}\|_{L^2}^2 dt
+ \fr{2}{1+X} (\tilde{g}, \l q)_{L^2} dW.
\ee After applying Young's inequality, we obtain 
\be 
\d \log(1+X) + \fr{\bar{X}}{1 + X} dt
\le C\|\Delta u\|_{L^2}^2 dt
+ \|\l^{\fr{1}{2}} \tilde{g}\|_{L^2}^2 dt
+ \fr{2}{1+X} (\tilde{g}, \l q)_{L^2} dW.
\ee 
Next, we integrate in time from $0$ to $t$, apply $\E$, and obtain
\be 
\E \int_{0}^{t} \fr{\bar{X}}{1+X} ds
\le \log(1 + X(0)) + C\int_{0}^{t} \E \|\Delta u(s)\|_{L^2}^2 ds
+ \|\l^{\fr{1}{2}} \tilde{g}\|_{L^2}^2 t.
\ee Therefore, 
\beg{align} \la{logestimate1}
\E \int_{0}^{t} \log (1 + \bar{X}) ds
&= \E \int_{0}^{t} \log \left(\fr{1 + \bar{X}}{1 + X} \right) ds
+ \E \int_{0}^{t} \log (1+X) ds \nonumber
\\&\quad\le \E \int_{0}^{t} \fr{\bar{X}}{1+X} ds
+ \E \int_{0}^{t} X ds.
\end{align} 
In view of the bounds \eqref{invariant00} and \eqref{invariant2}, we obtain \eqref{invariant6}, completing the proof.

\beg{prop} Let $u_0 \in V$ and $q_0 \in \dot{H}^1$. Suppose $g \in V$ and $\tilde{g} \in \dot{H}^1$. Then 
\beg{align} \la{invariant5}
&\E \left\{\int_{0}^{t} \|\l^{\fr{3}{2}} q(s)\|_{L^2}^2 ds \right\}  \nonumber
\\&\quad\le \|\na q_0\|_{L^2}^2 + C(f,g,\tilde{g}) \left[\|q_0\|_{L^4}^{12} + \|q_0\|_{L^4}^8 + \|q_0\|_{L^4}^4 + \|\na u_0\|_{L^2}^4 + \|q_0\|_{L^4}^4\|\na u_0\|_{L^2}^4 + t\right]
\end{align} holds for any $t \ge 0$. 
\end{prop}

\textbf{Proof:} By It\^o's lemma, we have
\beg{align} \la{invariant550}
&\d \|\na q\|_{L^2}^2
+ 2\|\l^{\fr{3}{2}} q\|_{L^2}^2 dt\nonumber\\
&\quad= 2(u \cdot \na q, \Delta q)_{L^2} dt
+ \|\na \tilde{g}\|_{L^2}^2 dt
- 2  (\tilde{g}, \Delta q)_{L^2} dW.
\end{align}  
In order to estimate the nonlinear term, we integrate by parts, use the divergence-free property $\nabla \cdot u=0$, to obtain 
\be 
(u \cdot \na q, \Delta q)_{L^2} 
= \sum\limits_{k,j \in \left\{1, 2\right\}} \int_{\TT^2} u_j \pa_j q \pa_{kk} q dx
= - \sum\limits_{k,j \in \left\{1, 2\right\}} \int_{\TT^2} \pa_k u_j \pa_j q \pa_{k} q dx.
\ee 
We bound
\beg{align}
&|(u \cdot \na q, \Delta q)_{L^2}| 
\le \|\na u\|_{L^4} \|\na q\|_{L^{\fr{8}{3}}}^2 \le C\|\na u\|_{L^4} \|q\|_{L^4}^{\fr{1}{2}} \|\l^{\fr{3}{2}} q\|_{L^2}^{\fr{3}{2}} \nonumber
\\&\quad\le  C\|\na u\|_{L^2}^{\fr{1}{2}}\|\Delta u\|_{L^2}^{\fr{1}{2}} \|q\|_{L^4}^{\fr{1}{2}} \|\l^{\fr{3}{2}} q\|_{L^2}^{\fr{3}{2}} 
\le \fr{1}{2} \|\l^{\fr{3}{2}} q\|_{L^2}^2 + C\|\na u\|_{L^2}^{2}\|\Delta u\|_{L^2}^{2} \|q\|_{L^4}^{2}
\end{align} in view of H\"older's inequality with exponents $4, 8/3, 8/3$, the interpolation estimate \cite[Proposition~2]{AI}
\be 
\|\l^{\fr{3}{2}} q\|_{L^2}^2 \ge C\|q\|_{L^4}^{-\fr{2}{3}} \|\na q\|_{L^{\fr{8}{3}}}^{\fr{8}{3}},
\ee and Ladyzhenskaya's interpolation inequality. We obtain
\be \la{EQ56}
\d \|\na q\|_{L^2}^2
+ \|\l^{\fr{3}{2}} q\|_{L^2}^2 dt
\le C\|\na u\|_{L^2}^{2} \|\Delta u\|_{L^2}^2 \|q\|_{L^4}^{2} dt
+ \|\na \tilde{g}\|_{L^2}^2 dt
- 2  (\tilde{g}, \Delta q)_{L^2} dW.
\ee 
 Hence, an application of Young's inequality yields
\beg{align} 
\d \|\na q\|_{L^2}^2
&+ \|\l^{\fr{3}{2}} q\|_{L^2}^2 dt
\le C\|\na u\|_{L^2}^{2} \|\Delta u\|_{L^2}^2 \|q\|_{L^4}^{4} dt \nonumber
\\&\quad+ C\|\na u\|_{L^2}^{2} \|\Delta u\|_{L^2}^2 dt
+ \|\na \tilde{g}\|_{L^2}^2 dt
- 2  (\tilde{g}, \Delta q)_{L^2} dW.
\end{align}
We integrate in time from $0$ to $t$ and we apply $\E$. In view of \eqref{invariant3} and \eqref{invariant4}, we obtain \eqref{invariant5}. 

\beg{prop} \la{bootstrap}  Let $k \ge 0$. Let $q_0 \in \dot{H}^{k+1}$ and $u_0 \in H^{k+2} \cap H$. Suppose $\tilde{g} \in \dot{H}^{k+1}$ and $g \in H^{k+2} \cap H$. If the estimate
\beg{align} \la{smoothcond}
&\E \int_{0}^{t} \log (1 + \|(-\Delta)^{\fr{k}{2} + \fr{1}{4}}q(s) \|_{L^2}^2 + \|(-\Delta)^{\fr{k+2}{2}} u(s) \|_{L^2}^2) ds \nonumber 
\\&\quad\quad\le C\log(1 + \|(-\Delta)^{\fr{k}{2}} q_0\|_{L^2}^2 + \|(-\Delta)^{\fr{k+1}{2}} u_0\|_{L^2}^2) \nonumber
\\&\quad\quad\quad\quad+ C(f, g, \tilde{g}, k) \left[\|\na q_0\|_{L^2}^{12} + \|\na u_0\|_{L^2} ^8 + 1 + t\right]
\end{align} holds for all $t \ge 0$, then the following estimate
\beg{align} \la{smoothresult}
&\E \int_{0}^{t} \log (1 + \|(-\Delta)^{\fr{k}{2} + \fr{3}{4}}q(s) \|_{L^2}^2 + \|(-\Delta)^{\fr{k+3}{2}} u(s) \|_{L^2}^2) ds \nonumber
\\&\quad\quad\le C\log(1 + \|(-\Delta)^{\fr{k+1}{2}} q_0\|_{L^2}^2 + \|(-\Delta)^{\fr{k+2}{2}} u_0\|_{L^2}^2) \nonumber
\\&\quad\quad\quad\quad+ C(f, g, \tilde{g}, k) \left[\|\na q_0\|_{L^2}^{12} + \|\na u_0\|_{L^2} ^8 + 1 + t\right]
\end{align} holds for all $t \ge 0$.
\end{prop}

\textbf{Proof:} The It\^o lemma yields
\beg{align}
&\d \|(-\Delta)^{\fr{k+1}{2}}q\|_{L^2}^2 + 2\|(-\Delta)^{\fr{k}{2} + \fr{3}{4}}q\|_{L^2}^2 dt \nonumber
\\&\quad= - 2(u \cdot \na q, (-\Delta)^{k+1}q)_{L^2} dt
+ \|(-\Delta)^{\fr{k+1}{2}}\tilde{g}\|_{L^2}^2 dt
+ 2(\tilde{g}, (-\Delta)^{k+1}q)_{L^2} dW
\end{align} and 
\beg{align}
&\d \|(-\Delta)^{\fr{k+2}{2}} u\|_{L^2}^2 
+2 \|(-\Delta)^{\fr{k+3}{2}} u\|_{L^2}^2 \nonumber
\\&\quad= -2(qRq + u \cdot \na u - f, (-\Delta)^{k+2}u)_{L^2} dt
+ \|(-\Delta)^{\fr{k+2}{2}} g\|_{L^2} dt
+ 2(g, (-\Delta)^{k+2}u)_{L^2} dW.
\end{align} 
Let 
\be 
X = \|(-\Delta)^{\fr{k+1}{2}}q\|_{L^2}^2 + \|(-\Delta)^{\fr{k+2}{2}}u\|_{L^2}^2,
\ee
\be 
\bar{X} = \|(-\Delta)^{\fr{k}{2} + \fr{3}{4}}q\|_{L^2}^2 + \|(-\Delta)^{\fr{k+3}{2}} u\|_{L^2}^2,
\ee
\be 
M = 2(\tilde{g}, (-\Delta)^{k+1}q)_{L^2} + 2(g, (-\Delta)^{k+2}u)_{L^2},
\ee and 
\be 
N = \|(-\Delta)^{\fr{k+1}{2}}\tilde{g}\|_{L^2}^2 + \|(-\Delta)^{\fr{k+2}{2}} g\|_{L^2}^2.
\ee Then the stochastic process $X$ evolves according to 
\be 
\d X + 2 \bar{X} dt
= - 2(u \cdot \na q, (-\Delta)^{k+1}q)_{L^2} dt 
-2(qRq + u \cdot \na u - f, (-\Delta)^{k+2}u)_{L^2} dt
+ N dt + M dW.
\ee An application of It\^o's lemma gives the stochastic energy equality 
\beg{align}
&\d \log(1+X) + \fr{2\bar{X}}{1+X} dt
=-\fr{2}{1+X}(u \cdot \na q, (-\Delta)^{k+1}q)_{L^2} dt  \nonumber
\\&\quad- \fr{2}{1+X} |(qRq + u \cdot \na u - f, (-\Delta)^{k+2}u)_{L^2} dt
+ \fr{N}{1+X} dt
- \fr{M^2}{2(1+X)^2} dt
+ \fr{M}{1+X} dW,
\end{align} from which we obtain the following differential inequality 
\beg{align}
&\d \log(1+X) + \fr{2\bar{X}}{1+X} dt
\le \fr{2}{1+X}|(u \cdot \na q, (-\Delta)^{k+1}q)_{L^2} | dt  \nonumber
\\&\quad+ \fr{2}{1+X} |(qRq + u \cdot \na u - f, (-\Delta)^{k+2}u)_{L^2}| dt
+ N dt 
+ \fr{M}{1+X} dW.
\end{align}
In view of the commutator estimate
\be \la{kato}
\|\l^{s} (FG) - F\l^s G\|_{L^p} \le C \|\na F\|_{L^{p_1}} \|\l^{s-1} G\|_{L^{p_2}} + C \|\l^{s} F\|_{L^{p_3}} \|G\|_{L^{p_4}}
\ee that holds for any $s > 0$, $p \in (1, \infty)$,  $p_2, p_3 \in (1,\infty)$, $\fr{1}{p} = \fr{1}{p_1} + \fr{1}{p_2} = \fr{1}{p_3} + \fr{1}{p_4}$, and all appropriately smooth functions $F$ and $G$ (see \cite[Lemma A.1]{AI}), we estimate 
\beg{align}
&\left|(u \cdot \na q, (-\Delta)^{k+1} q)_{L^2} \right| \nonumber
\le \|(-\Delta)^{\fr{k+1}{2}} q\|_{L^2} \|(-\Delta)^{\fr{k+1}{2}}(u \cdot \na q) - u \cdot \na (-\Delta)^{\fr{k+1}{2}}q \|_{L^2}
\\&\quad\le C\|(-\Delta)^{\fr{k+1}{2}} q\|_{L^2} \left[\|\na u\|_{L^4} \|(-\Delta)^{\fr{k+1}{2}} q \|_{L^4} + \|\na q\|_{L^4} \|(-\Delta)^{\fr{k}{2} + \fr{1}{2}} u\|_{L^4} \right] \nonumber
\\&\quad\le C\|\Delta u\|_{L^2} \sqrt{X} \sqrt{\bar{X}} + C\|\l^{\fr{3}{2}}q\|_{L^2}X.
\end{align}
 Here, we used the continuous Sobolev embedding of $H^{\fr{1}{2}}$ in $L^4$. In view of the fractional product estimate 
\be 
\|\l^{s} (FG)\|_{L^p}
\le C \left[\|F\|_{L^{p_1}} \|\l^{s} G\|_{L^{p_2}} + \|\l^s F\|_{L^{p_3}} \|G\|_{L^{p_4}} \right]
\ee that holds for any $s > 0$, $p \in (1, \infty)$,  $p_2, p_3 \in (1,\infty)$, $\fr{1}{p} = \fr{1}{p_1} + \fr{1}{p_2} = \fr{1}{p_3} + \fr{1}{p_4}$, and all appropriately smooth functions $F$ and $G$ (see \cite[Lemma A.1]{KS3}), we estimate 
\beg{align}
&\left|(qRq, (-\Delta)^{k+2}u)_{L^2}\right| \nonumber
= \left|((-\Delta)^{\fr{k+1}{2}} (qRq), (-\Delta)^{\fr{k+3}{2}} u)_{L^2} \right|
\\&\quad\le C\left[\|Rq\|_{L^{\infty}} \|(-\Delta)^{\fr{k+1}{2}}q\|_{L^2} + \|q\|_{L^{\infty}} \|(-\Delta)^{\fr{k+1}{2}} Rq\|_{L^2} \right] \|(-\Delta)^{\fr{k+3}{2}} u\|_{L^2} \nonumber
\\&\quad\le C\|\l^{\fr{3}{2}}q\|_{L^2} \|(-\Delta)^{\fr{k+1}{2}} q\|_{L^2} \|(-\Delta)^{\fr{k+3}{2}}u\|_{L^2} \nonumber
\\&\quad\le C \|\l^{\fr{3}{2}}q\|_{L^2} \sqrt{X} \sqrt{\bar{X}}
\end{align}
after integrating by parts, using the continuous Sobolev embedding of $H^{\fr{3}{2}}$ in $L^{\infty}$, and using the boundedness of the Riesz transform on $H^{\fr{3}{2}}$. 
As for the nonlinear term in $u$, we integrate by parts, apply the commutator estimate \eqref{kato}, use the continuous embedding of $H^{\fr{1}{2}}$ in $L^4$, and estimate 
\beg{align}
&|(u \cdot \na u, (-\Delta)^{k+2}u)_{L^2}|
= |((-\Delta)^{\fr{k+2}{2}} (u \cdot \na u), (-\Delta)^{\fr{k+2}{2}} u)_{L^2}| \nonumber
\\&\quad= |((-\Delta)^{\fr{k+2}{2}} (u \cdot \na u) - u \cdot \na (-\Delta)^{\fr{k+2}{2}} u, (-\Delta)^{\fr{k+2}{2}} u)_{L^2}| \nonumber
\\&\quad\le C\|\na u\|_{L^4} \|(-\Delta)^{\fr{k+2}{2}} u\|_{L^4}  \|(-\Delta)^{\fr{k+2}{2}} u\|_{L^2}
\\&\quad\le C\|\Delta u\|_{L^2} \sqrt{\bar{X}} \sqrt{X}.
\end{align} 
Therefore, we obtain the inequality 
\beg{align} 
&\d \log(1+X) + \fr{2\bar{X}}{1+X} dt 
\le N dt + \fr{M}{1+X} dW \nonumber
\\&\quad+ \fr{C}{1+X} \left[\|\Delta u\|_{L^2} \sqrt{X} \sqrt{\bar{X}} + \|(-\Delta)^{\fr{k+1}{2}}f\|_{L^2} \sqrt{\bar{X}} + \|\l^{\fr{3}{2}}q\|_{L^2}\sqrt{X} \sqrt{\bar{X}} \right] dt 
\end{align} which boils down to
\beg{align}
\d \log(1+X) + \fr{\bar{X}}{1+X} dt
&\le C\|\Delta u\|_{L^2}^2 dt
+ C\|\l^{\fr{3}{2}}q\|_{L^2}^2 dt \nonumber
\\&\quad\quad+ C\|(-\Delta)^{\fr{k+1}{2}} f\|_{L^2}^2 dt
+N dt + \fr{M}{1+X} dW
\end{align}
after application of Young's inequality. We integrate in time from $0$ to $t$ and we apply $\E$. Using the bounds \eqref{invariant2} and \eqref{invariant5}, and applying Young's inequality, we conclude that 
\be 
\E \int_{0}^{t} \fr{\bar{X}}{1+X} ds \le \log (1 + X(0)) + C(f, g, \tilde{g}, k) (\|\na q_0\|_{L^2}^{12} + \|\na u_0\|_{L^2}^8 + 1 + t)
\ee  for all $t \ge 0$. Bounding similarly to \eqref{logestimate1}, we have
\beg{align}
\E \int_{0}^{t} \log(1+ \bar{X}) ds
&\le  \log (1 + X(0)) + C(f, g, \tilde{g}, k) (\|\na q_0\|_{L^2}^{12} + \|\na u_0\|_{L^2}^8 + 1 + t) \nonumber
\\&\quad\quad+ \E \int_{0}^{t} \log (1 + X) ds.
\end{align}
Since
\be 
X \le \|(-\Delta)^{\fr{k}{2} + \fr{3}{4}} q\|_{L^2} \|(-\Delta)^{\fr{k}{2} + \fr{1}{4}} q\|_{L^2}
+ \|(-\Delta)^{\fr{k+2}{2}} u\|_{L^2} \|(-\Delta)^{\fr{k+3}{2}} u\|_{L^2},
\ee  
we have
\be 
1 + X \le \left[1 + \bar{X}\right]^{\fr{1}{2}} \left[1 + \|(-\Delta)^{\fr{k}{2} + \fr{1}{4}} q\|_{L^2}^2 + \|(-\Delta)^{\fr{k+2}{2}} u\|_{L^2}^2 \right]^{\fr{1}{2}}
\ee and so
\be 
\log(1+X)\le \fr{1}{2} \log(1+ \bar{X}) + \fr{1}{2} \log (1 + \|(-\Delta)^{\fr{k}{2} + \fr{1}{4}} q\|_{L^2}^2 + \|(-\Delta)^{\fr{k+2}{2}} u\|_{L^2}^2).
\ee Therefore,
\beg{align}
\fr{1}{2} \E \int_{0}^{t} \log(1 + \bar{X}) ds
&\le \log (1 + X(0)) + C(f, g, \tilde{g}, k) (\|\na q_0\|_{L^2}^{12} + \|\na u_0\|_{L^2}^8 + 1 + t) \nonumber
\\&\quad\quad+ \fr{1}{2}\E \int_{0}^{t}  \log (1 + \|(-\Delta)^{\fr{k}{2} + \fr{1}{4}} q\|_{L^2}^2 + \|(-\Delta)^{\fr{k+2}{2}} u\|_{L^2}^2) ds.
\end{align} 
In view of \eqref{smoothcond}, we obtain \eqref{smoothresult}.

We end this section by proving Theorem \ref{higherreg}. 

\textbf{Proof of Theorem \ref{higherreg}:} Suppose $\mu$ is an invariant measure of \eqref{invsys}. By Theorem \ref{weakinvariance}, $\mu$ is supported on $H^{\fr{1}{2}} \times (H^2 \cap H)$. In view of the bounds \eqref{invariant6} and \eqref{invariant5}, and repeating the same argument used to prove Theorem \ref{weakinvariance}, we conclude that $\mu$ is supported on $H^{\fr{3}{2}} \times (H^2 \cap H)$. Now we bootstrap using Proposition \ref{bootstrap} and we deduce that $\mu$ is supported on $H^{k + \fr{3}{2}} \times H^{k+3}$ for any $k \ge 0$. This shows that $\mu$ is smooth and completes the proof of Theorem \ref{higherreg}. 
 
\section{Uniqueness of Invariant Measures} \la{sec7}

In this section, we prove that \eqref{invsys} has a unique ergodic invariant measure provided that the ranges of $\tilde{g}$ and $g$ are large enough in phase space. 
Uniqueness is obtained by employing asymptotic coupling arguments from \cite{GMR}.

\beg{thm} \la{t7} Suppose that $g \in V$ and $\tilde{g} \in \dot{L}^4$. There exists $N = N(f, g, \tilde{g})$ such that if $P_N \mathcal{\dot{H}} \subset range(\tilde{g}, g)$, then  \eqref{invsys} has a unique ergodic invariant measure.
\end{thm}

In order to prove Theorem \ref{t7}, we need the following proposition:

\beg{prop} \la{p11} Let $R > 0$. Then there exist positive universal constants $c$ and $C$ such that the estimates
\beg{align}  \la{111}
\MoveEqLeft[4] \PP \left(\sup\limits_{t \ge 0} \left(\|\na u(t)\|_{L^2}^2 + \fr{1}{2} \int_{0}^{t} \|\Delta u(s)\|_{L^2}^2 ds - \|\na u_0\|_{L^2}^2 \right.\right. \nonumber
\\&\left.{} \left.{}- C(\|f\|_{L^2}^2 + \|\na g\|_{L^2}^2) t - C\int_{0}^{t} \|q(s)\|_{L^4}^4 ds \right)  > R\right)  
 \le e^{-\fr{R}{8 \|g\|_{L^2}^2}}
\end{align} 
and
\be \la{112}
\PP \left(\sup\limits_{t \ge 0} \left(\|q(t)\|_{L^4}^4 + c \int_{0}^{t} \|q(s)\|_{L^4}^4 ds - \|q_0\|_{L^4}^4 - 2 - C\|\tilde{g}\|_{L^4}^4 t \right) > R\right) \le \fr{C\left(\|\tilde{g}\|_{L^4}^{16} + \|q_0\|_{L^4}^{16} \right)}{R+1}
\ee hold. 
\end{prop}

\textbf{Proof:}  We integrate in time from $0$ to $t$ the differential inequality  
\beg{align} 
\d \|\na u\|_{L^2}^2 +  \|\Delta u\|_{L^2}^2 dt
\le C\|q\|_{L^4}^4  dt 
+ C\|f\|_{L^2}^2 dt
+ C\|\na g\|_{L^2}^2 dt
- 2 (g, \Delta u)_{L^2} dW
\end{align}  
(see \eqref{gradvelenergy})
and take the supremum over $t \ge 0$ to obtain
\beg{align} 
&\sup\limits_{t \ge 0} \left\{\|\na u(t)\|_{L^2}^2 + \fr{1}{2} \int_{0}^{t} \|\Delta u(s)\|_{L^2}^2 ds
- \|\na u_0\|_{L^2}^2 - C(\|f\|_{L^2}^2 + \|\na g\|_{L^2}^2) t - C\int_{0}^{t} \|q(s)\|_{L^4}^4 ds \right\} \nonumber
\\&\quad\le \sup\limits_{t \ge 0} \left\{\int_{0}^{t} 2(g, - \Delta u)_{L^2} dW(s) - \fr{1}{2} \int_{0}^{t} \|\Delta u\|_{L^2}^2 ds  \right\} \nonumber
\\&\quad=\sup\limits_{t \ge 0} \left\{\int_{0}^{t} 2(g, - \Delta u)_{L^2} dW(s) -  \fr{1}{8 \|g\|_{L^2}^2} \int_{0}^{t} 4 \|g\|_{L^2}^2 \|\Delta u\|_{L^2}^2 ds  \right\}.
\end{align}  
Exponential martingale inequalities \cite[(3.4)]{GMR} imply 
\be 
\PP (\sup\limits_{t \ge 0} \left\{\int_{0}^{t} 2(g, - \Delta u)_{L^2} dW(s) - \fr{1}{8 \|g\|_{L^2}^2} \int_{0}^{t} 4 \|g\|_{L^2}^2 \|\Delta u\|_{L^2}^2 ds  \right\} > R)
\le e^{-\fr{R}{8 \|g\|_{L^2}^2}}.
\ee 
Therefore \eqref{111} is established. The derivation of \eqref{112} is based on ideas from \cite{GMaR}. Indeed, the $L^4$ norm of $q$ evolves according to
\beg{align}  
\d \|q\|_{L^4}^4 
+ 4 (\l q, q^3)_{L^2} dt 
= 6 (\tilde{g}^2, q^2)_{L^2} dt
+ 4  (\tilde{g}, q^3)_{L^2} dW.
\end{align} 
(see \eqref{pcharge}).
By the Poincar\'e inequality for the fractional Laplacian in $L^4$, we have 
\be 
(\l q, q^3)_{L^2} \ge c\|q\|_{L^4}^4 
\ee 
Thus, we obtain the differential inequality
\be 
\d \|q\|_{L^4}^4 + c\|q\|_{L^4}^4 dt  
\le  C \|\tilde{g}\|_{L^4}^4 dt
+ 4 (\tilde{g}, q^3)_{L^2} dW.
\ee
We integrate from $0$ to $t$, and take the supremum over $t\ge 0$. We obtain 
\beg{align} 
&\sup\limits_{t \ge 0} \left\{\|q(t)\|_{L^4}^4 + c\int_{0}^{t} \|q(s)\|_{L^4}^4 ds - \|q_0\|_{L^4}^4 - 2 -C\|\tilde{g}\|_{L^4}^4 t \right\} \nonumber
\\&\quad\quad\quad\quad\le \sup\limits_{t \ge 0} \left\{\int_{0}^{t} 4(\tilde{g}, q^3)_{L^2} dW(s) - t - 2 \right\} 
\end{align} 
which implies
\beg{align} 
&\PP \left(\sup\limits_{t \ge 0} \left(\|q(t)\|_{L^4}^4 + c\int_{0}^{t} \|q(s)\|_{L^4}^4 ds - \|q_0\|_{L^4}^4 - 2  -C\|\tilde{g}\|_{L^4}^4 t \right) \ge R \right) \nonumber
\\&\quad\quad\quad\quad\le \PP \left(\sup\limits_{t \ge 0} (M(t) - t - 2) \ge R \right)
\end{align} 
for any $R > 0$, where $M(t)$ is the martingale term 
\be 
M(t) = 4\int_{0}^{t} (\tilde{g}, q^3)_{L^2} dW(s).
\ee 
We have
\be 
\left\{\sup\limits_{t \ge 0} (M(t) - t - 2) \ge R \right\} \subset \bigcup_{n \ge 0} \left\{\sup\limits_{t \in [n, n+1)} (M(t) - t - 2) \ge R \right\}
\ee
and 
\be 
\left\{\sup\limits_{t \in [n, n+1)} (M(t) - t - 2) \ge R \right\}
\subset \left\{M^*(n+1) \ge R + n + 2 \right\}
\ee 
where 
\be 
M^*(t) = \sup\limits_{s \in [0,t]} |M(s)|.
\ee Using the Burkholder-Davis-Gundy inequality \cite{Kara}
\be 
\E M^*(t)^4 \le C\E ([M](t)^2)
\ee where $[M](t)$ is the quadratic variation 
\be 
[M](t) = 16\int_{0}^{t} (\tilde{g}, q^3)_{L^2}^2 ds,
\ee we obtain 
\beg{align} 
&\E M^*(t)^4 
\le C \E([M](t)^2)
\le C \E \left(\int_{0}^{t} (\tilde{g}, q^3)_{L^2}^2 ds \right)^2 \nonumber
\\&\quad\quad\le C\|\tilde{g}\|_{L^4}^4 \E \left(\int_{0}^{t} \|q\|_{L^4}^6  ds\right)^2 
\le C \|\tilde{g}\|_{L^4}^4 t \E \left(\int_{0}^{t} \|q\|_{L^4}^{12} ds \right) \nonumber
\\&\quad\quad\le C\|\tilde{g}\|_{L^4}^4t \left(\|q_0\|_{L^4}^{12} + \|\tilde{g}\|_{L^4}^{12}t \right)
\le C\left(\|\tilde{g}\|_{L^4}^{16} + \|q_0\|_{L^4}^{16} \right)(1+t)^2.
\end{align} Here we used the estimate \eqref{invariant1} applied for $p=12$.
Therefore, 
\beg{align}
&\PP \left(\sup\limits_{t \ge 0} (M(t) - t -2)\ge R \right)
\le \sum\limits_{n \ge 0} \PP (M^*(n+1) \ge R + n + 2)\nonumber
\\&\quad\quad\le \sum\limits_{n \ge 0} \fr{\E M^*(n+1)^4}{(R+ n +2)^4} 
\le C\left(\|\tilde{g}\|_{L^4}^{16} + \|q_0\|_{L^4}^{16} \right) \sum\limits_{n \ge 0} \fr{(n+2)^2}{(R+n+2)^4} \nonumber
\\&\quad\quad\le C\left(\|\tilde{g}\|_{L^4}^{16} + \|q_0\|_{L^4}^{16} \right)
\sum\limits_{n \ge 0} \fr{1}{(R + n + 2)^2}
\le \fr{C\left(\|\tilde{g}\|_{L^4}^{16} + \|q_0\|_{L^4}^{16} \right)}{R + 2} 
\end{align} in view of the Chebyshev's inequality. This gives \eqref{112} ending the proof of Proposition \ref{p11}.

Finally, we prove the uniqueness result:

\textbf{Proof of Theorem \ref{t7}:} Fix $(q_0, u_0)$ and $(Q_0, U_0)$ in $\dot{\mathcal{V}}$. Our aim is to establish the conditions for the asymptotic coupling framework presented in Section 2.4 of \cite{GMR}. To this end, we consider $(q,u)$ solving \eqref{invsys} with $(q(0), u(0)) = (q_0,u_0)$, and $(Q,U)$ solving 
\be \begin{cases} \la{invSYS}
\d (Q, U) +  (\l Q, -\Delta U) dt + (0, \na P) dt 
= (-U \cdot \na Q, -U \cdot \na U - Q RQ  + f ) dt 
\\\quad\quad+ (\tilde{g}, g)dW + {\bf{1}}_{\tau_K > t} \lambda P_N(q-Q, u-U) dt 
\\ \na \cdot U = 0
\end{cases}
\ee with $(Q(0), U(0)) = (Q_0, U_0)$, where
\be 
\tau_K = \inf\limits_{t \ge 0} \left\{\int_{0}^{t} \|P_N(\l^{-\fr{1}{2}} (q - Q),(u-U))\|_{L^2}^2 ds \ge K \right\}.
\ee and $K$, $N$ and $\lambda$ are positive constants to be determined later.

By Girsanov's theorem \cite[Theorem 2.2]{GMR}, the law of $(Q,U)$ is absolutely continuous with respect to the solution $(q,u)(\cdot, (Q_0,U_0))$ of \eqref{invsys} corresponding to $(Q_0, U_0)$ for any choices of $\lambda>0$ and $K>0$. Consequently, the uniqueness of the invariant measure follows from an application of Corollary 2.1 in \cite{GMR}, provided that we can find some positive constants $\lambda$ and $K$ such that $(q,u) - (Q,u) \rightarrow 0$ in the norm of $\mathcal{\dot{H}}$ on a set of positive measure.

Let 
\be 
v = u - U, \pi = p - P, \xi = q - Q.
\ee
Then $(\xi, v)$ obeys 
\beg{align} \la{weq}
&\partial_t (\xi,v) + (\l \xi, -\Delta v) + {\bf{1}}_{\tau_K > t} \lambda P_N (\xi, v) + (0, \pi) \nonumber
\\&\quad= \left(- u \cdot \na q + U \cdot \na Q,  - u \cdot \na u + U \cdot \na U - qRq + QRQ  \right)
\end{align}
Let $\omega = ( \xi, v)$. 
Taking the $L^2$ inner product of \eqref{weq} with $(\l^{-1} \xi, v)$, we obtain the differential inequality
\beg{align} 
&\fr{1}{2} \fr{d}{dt} \|\omega\|_{\mathcal{\dot{H}}}^2 + \|\xi\|_{L^2}^2 +  \|\na v\|_{L^2}^2 + {\bf{1}}_{\tau_K > t} \lambda \|P_N (\l^{-\fr{1}{2}} \xi, v)\|_{L^2}^2 \nonumber\\
&\quad= (- u \cdot \na q + U \cdot \na Q, \l^{-1}\xi)_{L^2} + ( - u \cdot \na u + U \cdot \na U - qRq + QRQ, v)_{L^2} \nonumber
\\&\quad= (-v \cdot \na Q - u \cdot \na \xi, \l^{-1} \xi)_{L^2}
+ (-v \cdot \na u - U \cdot \na v, v)_{L^2}
+ (-\xi Rq - QR\xi, v)_{L^2} \nonumber
\\&\quad= - (u \cdot \na \xi, \l^{-1}\xi)_{L^2} - (\xi Rq, v)_{L^2} 
- (v \cdot \na u, v)_{L^2}
\end{align}
where we used the cancellations
\be 
(U \cdot \na v, v )_{L^2} = 0
\ee   and 
\be 
(v \cdot \na Q, \l^{-1}\xi)_{L^2} = - (v \cdot R\xi, Q)_{L^2} = - (QR\xi, v)_{L^2}.
\ee
We estimate
\be 
|(v \cdot \na u, v)_{L^2}|
\le \|v\|_{L^4}^2\|\na u\|_{L^2}
\le C\|v\|_{L^2} \|\na v\|_{L^2} \|\na u\|_{L^2}
\le \fr{1}{4} \|\na v\|_{L^2}^2 + C\|\na u\|_{L^2}^{2}\|v\|_{L^2}^2,
\ee 
\be 
|(\xi Rq, v)_{L^2}| \le \|\xi\|_{L^2}\|v\|_{L^4}\|Rq\|_{L^4} 
\le \fr{1}{4} \|\xi\|_{L^2}^2 + \fr{1}{4} \|\na v\|_{L^2}^2 + C\|q\|_{L^4}^4 \|v\|_{L^2}^2,
\ee and
\beg{align} 
&|(u \cdot \na \xi, \l^{-1}\xi)_{L^2}|
= |(\l^{-\fr{1}{2}}(u \cdot \na \xi) - u \cdot \na \l^{-\fr{1}{2}}\xi, \l^{-\fr{1}{2}}\xi)_{L^2}| \nonumber
\\&\le C\|\Delta u\|_{L^2} \|\l^{-\fr{1}{2}} \xi\|_{L^2}\|\xi\|_{L^2}
\le \fr{1}{4} \|\xi\|_{L^2}^2 + C\|\Delta u\|_{L^2}^2 \|\l^{-\fr{1}{2}}\xi\|_{L^2}^2
\end{align} using H\"older's inequality, Ladyzhenskaya's interpolation inequality, Young's inequality, the boundedness of the Riesz transform on $L^4$, and the commutator estimate \eqref{AIcom}. 
This yields the differential inequality
\beg{align} 
& \fr{d}{dt} \|\omega\|_{\mathcal{\dot{H}}}^2 + \|\xi\|_{L^2}^2 +  \|\na v\|_{L^2}^2 + {\bf{1}}_{\tau_K > t} \lambda \|P_N (\l^{-\fr{1}{2}} \xi, v)\|_{L^2}^2 \nonumber\\
&\quad\le (C + \|\Delta u\|_{L^2}^2 + \|q\|_{L^4}^4) \|\omega\|_{\mathcal{\dot{H}}}^2.
\end{align}
For a fixed integer $N$, we have
\beg{align} 
&\|\xi\|_{L^2}^2  + \|\na v\|_{L^2}^2 + {\bf{1}}_{\tau_K > t}\lambda \|P_N (\l^{-\fr{1}{2}} \xi, v)\|_{L^2}^2 \nonumber
\\&\quad\ge {\bf{1}}_{\tau_K > t} (\lambda_{N+1}^{\fr{1}{2}} \|Q_N(\l^{-\fr{1}{2}} \xi, v)\|_{L^2}^2 + \lambda_{N+1}^{\fr{1}{2}} \|P_N(\l^{-\fr{1}{2}}\xi, v)\|_{L^2}^2) \nonumber
\\&\quad\ge {\bf{1}}_{\tau_K > t} \lambda_{N+1}^{\fr{1}{2}} \|\omega\|_{\mathcal{\dot{H}}}^2
\end{align} for $\lambda \ge \lambda_{N+1}^{\fr{1}{2}}$ in view of the inequality \eqref{gP}.
Hence
\be 
\fr{d}{dt} \|\omega\|_{\mathcal{\dot{H}}}^2  + {\bf{1}}_{\tau_K > t} \lambda_{N+1}^{\fr{1}{2}} \|\omega\|_{\mathcal{\dot{H}}}^2 
\le (C +  C\|q\|_{L^4}^4 + C\|\Delta u\|_{L^2}^2)\|\omega\|_{\mathcal{\dot{H}}}^2.  
\ee 
Integrating in time, we obtain
\be 
\|\omega(t)\|_{\mathcal{\dot{H}}}^2 \le \|\omega_0\|_{\mathcal{\dot{H}}}^2 \exp \left\{-\lambda_{N+1}^{\fr{1}{2}}t + \int_{0}^{t} (C +   C\|q\|_{L^4}^4 + C\|\Delta u\|_{L^2}^2) ds  \right\}
\ee for any $t \in [0, \tau_K]$.
For $R \ge 0$, we consider the sets
\begin{align}
\MoveEqLeft[4] E_R = \left\{\sup\limits_{t \ge 0} \left( \|\na u(t)\|_{L^2}^2 + \fr{1}{2} \int_{0}^{t} \|\Delta u(s)\|_{L^2}^2 ds \right.\right. \nonumber\\
&\left.{}\left.{}- \|\na u_0\|_{L^2}^2 - C(\|f\|_{L^2}^2 + \|\na g\|_{L^2}^2) t - C\int_{0}^{t} \|q(s)\|_{L^4}^4 ds \right)   \le R\right\}
\end{align}
and
\be 
F_R = \left\{\sup\limits_{t \ge 0} \left(\|q(t)\|_{L^4}^4 + c \int_{0}^{t} \|q(s)\|_{L^4}^4 ds - \|q_0\|_{L^4}^4 - 2 - C\|\tilde{g}\|_{L^4}^4 t \right) \le R\right\}.
\ee  
By Proposition \ref{p11}, we have $\PP(E_R \cap F_R) > 0$ when $R$ is sufficiently large.
Indeed,
\beg{align} 
\PP(E_R \cap F_R) &= \PP(E_R) + \PP(F_R) - \PP(E_R \cup F_R) 
> 1 -e^{-\fr{R}{8 \|g\|_{L^2}^2}} - \fr{C\left(\|\tilde{g}\|_{L^4}^{16} + \|q_0\|_{L^4}^{16} \right)}{R+1}> 0
\end{align} when $R$ is large. 
Consequently, on $E_R \cap F_R$ and for $t \in [0,\tau_K]$, we have 
\be 
\|\omega(t)\|_{\mathcal{\dot{H}}}^2 \le \|\omega_0\|_{\mathcal{\dot{H}}}^2 e^{-\fr{t}{2} \lambda_{N+1}^{\fr{1}{2}}} e^{\left(-\fr{1}{2} \lambda_{N+1}^{\fr{1}{2}} +  C(f, g, \tilde{g})\right)t  } e^{C(\|\na u_0\|_{L^2}, \|q_0\|_{L^4}, R)}.
\ee
We choose an integer $N = N(f, g, \tilde{g})$ large enough so that 
\be 
-\fr{1}{2} \lambda_{N+1}^{\fr{1}{2}} + C(f, g, \tilde{g}) \le 0
\ee 
yielding
\be 
\|\omega(t)\|_{\mathcal{\dot{H}}}^2 \le \|\omega_0\|_{\mathcal{\dot{H}}}^2 e^{-\fr{t}{2} \lambda_{N+1}^{\fr{1}{2}}}  e^{C(\|\na u_0\|_{L^2}, \|q_0\|_{L^4}, R) }
\ee on $E_R \cap F_R$ and for $t \in [0,\tau_K]$.
Finally, we choose $K$ large enough such that $E_R \cap F_R  \subset \left\{\tau_K = \infty \right\}$ and we conclude that on the nontrivial set $E_R \cap F_R $
\be 
(q(t) - Q(t), u(t) - U(t)) \rightarrow 0
\ee in $\mathcal{\dot{H}}$ as $t \rightarrow \infty$. This completes the proof of Theorem \ref{t7}.

\section{Feller Property in the $H^1$ norm} \la{SFP}

We consider the space
\be  
\tilde{\mathcal{V}} = \dot{H}^{1}(\TT^2) \times V
\ee with norm
\be 
\|(\xi, v)\|_{\mathcal{\tilde{V}}} = \|\na \xi\|_{L^2} + \|\na v\|_{L^2}^2.
\ee

In this section, we show that the transition kernels associated with \eqref{invsys} are Feller in the norm of $\mathcal{\tilde{V}}$. 

\beg{thm} \la{stfel} Suppose that $g \in \cap H^2 \cap H$ and $\tilde{g} \in \dot{H}^1$ such that $\na \tilde{g} \in L^8$. Then the semigroup $\left\{\tilde{P}_t \right\}_{t \ge 0}$ is Markov-Feller on $C_b (\mathcal{\tilde{V}}).$
\end{thm}

We need the following propositions.

\beg{prop} (Continuity in $\tilde{V}$) Let $(q_0^1, u_0^1)$ and $(q_0^2, u_0^2)$ be in $\mathcal{\tilde{V}}$. Suppose $\tilde{g} \in \dot{H}^{1}$ and $g \in V$. Then the corresponding solutions $(q_1, u_1)$ and $(q_2, u_2)$ obey 
\beg{align} \la{strongfeller2}
&\|\na u_1(t) - \na u_2(t)\|_{L^2}^2 + \|\na q_1(t) - \na q_2(t)\|_{L^2}^2 \nonumber\\
&\quad\le \exp \left\{C\tilde{C}(t)\right\} \left[\|\na u_{0}^1 - \na u_{0}^2\|_{L^2}^2 + \|\na q_{0}^1 - \na q_{0}^2 \|_{L^2}^2 \right]
\end{align} with probability 1, where 
\be 
\tilde{C}(t) = \int_{0}^{t} \left[\|\l^{\fr{3}{2}} q_1\|_{L^2}^2 + \|\na u_1\|_{L^2}^2 + \|q_2\|_{L^4}^2 + \|\Delta u_2\|_{L^2}^2 \right] ds
\ee is well-defined and finite almost surely. 
\end{prop}

\textbf{Proof:} Let $q = q_1 - q_2$ and $u = u_1 - u_2$. The norm $\|\na q\|_{L^2}$ satisfies the energy inequality
\beg{align}
&\fr{1}{2} \fr{d}{dt} \|\na q\|_{L^2}^2
+ \|\l^{\fr{3}{2}} q\|_{L^2}^2 
\le \left|\int_{\TT^2} (u \cdot \na q_1) \Delta q  \right| 
+ \left|\int_{\TT^2} (u_2 \cdot \na q) \Delta q \right| \nonumber
\\&\le C\|\na u\|_{L^4} \|\na q\|_{L^2} \|\na q_1\|_{L^4}
+ C\|\na u_2\|_{L^4} \|\na q\|_{L^4} \|\na q\|_{L^2}.
\end{align} where we integrated by parts and used the divergence-free condition of $u_2$ and $u$. 
Applying Young's inequality and using the continuous embedding of $H^{\fr{1}{2}}$ in $L^4$, we obtain 
\be \la{ssss1}
\fr{d}{dt} \|\na q\|_{L^2}^2 + \|\l^{\fr{3}{2}} q\|_{L^2}^2 
\le \fr{1}{8} \|\Delta u\|_{L^2}^2 + C\left[\|\l^{\fr{3}{2}} q_1\|_{L^2}^2 + \|\Delta u_2\|_{L^2}^2 \right]\|\na q\|_{L^2}^2.
\ee On other hand, the norm $\|\na u\|_{L^2}$ obeys 
\beg{align}
\fr{1}{2} \fr{d}{dt} \|\na u\|_{L^2}^2 + \|\Delta u\|_{L^2}^2
&\le \left|\int_{\TT^2} (u \cdot \na u_1) \Delta u \right|
+ \left|\int_{\TT^2} (u_2 \cdot \na u) \Delta u \right|
+ \left|\int_{\TT^2} (q_1 Rq) \Delta u \right|
+ \left|\int_{\TT^2} (qRq_2) \Delta u \right| \nonumber
\\&\le C\|\na u\|_{L^4}^2 \|\na u_1\|_{L^2}
+ C\|\na u\|_{L^4}^2 \|\na u_2\|_{L^2} \nonumber
\\&\quad+ C\|q_1\|_{L^4} \|\l^{\fr{1}{2}}q\|_{L^2}\|\Delta u\|_{L^2}
+ C\|q_2\|_{L^4} \|\l^{\fr{1}{2}}q \|_{L^2} \|\Delta u\|_{L^2},
\end{align} hence
\be \la{ssss2}
\fr{d}{dt} \|\na u\|_{L^2}^2 + \|\Delta u\|_{L^2}^2
\le C\left[\|\na u_1\|_{L^2}^2 + \|\na u_2\|_{L^2}^2 \right]\|\na u\|_{L^2}^2 + C\left[\|q_1\|_{L^4}^2 + \|q_2\|_{L^4}^2 \right]\|\na q\|_{L^2}^2.
\ee Adding \eqref{ssss1} and \eqref{ssss2}, we get 
\be 
\fr{d}{dt} \left[\|\na q\|_{L^2}^2 + \|\na u\|_{L^2}^2 \right] 
\le C \left[\|\l^{\fr{3}{2}} q_1\|_{L^2}^2 + \|\na u_1\|_{L^2}^2 + \|q_2\|_{L^4}^2 + \|\Delta u_2\|_{L^2}^2 \right]\left[\|\na q\|_{L^2}^2 + \|\na u\|_{L^2}^2 \right] 
\ee which gives \eqref{strongfeller2}.

\beg{prop}\la{pathneeded} Suppose $\na \tilde{g} \in L^8$ and $\Delta g \in L^2$. Let $(q_0, u_0) \in \mathcal{\dot{V}}$ and $T > 0$. Then the solution $(q,u)$ to the system \eqref{invsys} is uniformly bounded (almost surely) in 
\be 
L_{loc}^2(0, \infty; L^4(\TT^2)) \times L_{loc}^2(0,\infty; H^2(\TT^2))
\ee 
by some constant depending only on $g, \tilde{g}, f, \|\na u_0\|_{L^2}$ and $\|q_0\|_{L^4}$. 
Consequently, if $(\xi_n, v_n) \in \mathcal{\tilde{V}}$ is a sequence of initial datum such that $\left\{(\xi_n, v_n)\right\}_{n=1}^{\infty}$ converges to $(\xi, v)$ in $\mathcal{\tilde{V}}$, then
\be \la{strongfeller3}
\limsup\limits_{n \to \infty} \int_{0}^{T} \left[\|q(t, (\xi_n, v_n))\|_{L^4}^2 + \|\Delta u(t, (\xi_n, v_n))\|_{L^2}^2 \right] dt < \infty
\ee almost surely. 
\end{prop}

The proof of Proposition \ref{pathneeded} is presented in Appendix B.

Now we prove Theorem \ref{stfel}:

\textbf{Proof of Theorem \ref{stfel}:}  Fix $\varphi \in C_b(\mathcal{\tilde{V}})$. Suppose $(\xi_n, v_n)$ converges to $(\xi, v)$ in $\mathcal{\tilde{V}}$, that is 
\be  \la{semi41}
\|\na \xi_n - \na \xi\|_{L^2}^2 + \|\na v_n - \na v\|_{L^2}^2 \rightarrow 0.
\ee In view of the continuity in $\mathcal{\tilde{V}}$ given by \eqref{strongfeller2}, we have 
\beg{align} 
&\|\na u(t, (\xi_n, v_n)) - \na u(t, (\xi,v))\|_{L^2}^2 
+ \|\na q(t, (\xi_n, v_n)) - \na q(t, (\xi, v))\|_{L^2}^2 \nonumber
\\&\le e^{C(K_n(t) + K(t)} \left[\|\na v_n - \na v\|_{L^2}^2 + \|\na \xi_n - \na \xi\|_{L^2}^2 \right]
\end{align} where
\be 
K_n(t) = \int_{0}^{t} \left[\|q(t, (\xi_n, v_n))\|_{L^4}^2 + \|\Delta u (t, (\xi_n, v_n))\|_{L^2}^2 \right] ds 
\ee and 
\be 
K(t) = \int_{0}^{t} \left[\|\l^{\fr{3}{2}} q(t, (\xi,v))\|_{L^2}^2 + \|\na u(t, (\xi, v))\|_{L^2}^2 \right] ds.
\ee
In view of \eqref{invariant2} and \eqref{invariant5}, we have the finiteness of $K(t)$ for almost every $w \in \Omega$. In view of \eqref{strongfeller3}, we have 
\be 
\limsup\limits_{n \to \infty}  K_n(s) ds < \infty 
\ee for almost every $w \in \Omega$. This implies that 
\be 
\|\na u(t, (\xi_n, v_n)) - \na u(t, (\xi,v))\|_{L^2}^2 
+ \|\na q(t, (\xi_n, v_n)) - \na q(t, (\xi, v))\|_{L^2}^2 \rightarrow 0.
\ee  Since $\varphi$ is continuous on $\mathcal{\tilde{V}}$, we conclude that 
\be 
\varphi ((q,u)(t, (\xi_n, v_n))) \rightarrow \varphi ((q,u)(t, (\xi, v))) 
\ee and hence 
\be 
\E \varphi ((q,u)(t, (\xi_n, v_n)))\rightarrow \E \varphi ((q,u)(t, (\xi, v)))
\ee due to the boundedness of $\varphi$. This completes the proof of Theorem \ref{stfel}.

\appendix

\section{Uniform Bounds in Lebesgue Spaces} \la{A}

In this Appendix, we prove Proposition \ref{prop1}. For simplicity, we ignore the viscous term $-\epsilon \Delta q^{\epsilon}$ in \eqref{stochastic} because it does not have any major contribution in estimating the solutions of the mollified system \eqref{nonlinear} and vanishes as we take the limit $\epsilon\to 0$. We also drop the $\epsilon$ superscript. 

The proof is divided into 7 main steps.

\textbf{Step 1.} We prove that the estimate \eqref{SS2} holds when $p=2$. 

\textbf{Proof of Step 1.} By It\^o's lemma, we have 
\be
\d q^2 
= -2 q (u \cdot \na q) dt 
-2 q \l q dt 
+2  q  \Delta \Phi dt 
+  \tilde{g}^2 dt
+ 2  q  \tilde{g} dW.
\ee 
We integrate in the space variable over $\TT^2$. In view of the divergence-free condition obeyed by $u$, the nonlinear term vanishes, that is 
\be 
(u \cdot \na q, q)_{L^2} = 0,
\ee  which yields the energy equality
\be \la{ql2}
\d \|q\|_{L^2}^2 + 2\|\l^{\fr{1}{2}} q\|_{L^2}^2
= 2 (\Delta \Phi, q)_{L^2} +  \| \tilde{g}\|_{L^2}^2 dt 
+ 2  ( \tilde{g}, q)_{L^2} d W.
\ee We estimate 
\be 
|(\Delta \Phi, q)_{L^2}| = |(\l^{\fr{3}{2}} \Phi, \l^{\fr{1}{2}} q)_{L^2}| \le \fr{1}{2} \|\l^{\fr{3}{2}} \Phi\|_{L^2}^2 + \fr{1}{2} \|\l^{\fr{1}{2}} q\|_{L^2}^2
\ee using the H\"older and Young inequalities. We obtain the differential inequality 
\be 
\d \|q\|_{L^2}^2 + \|\l^{\fr{1}{2}} q\|_{L^2}^2 dt
\le \|\l^{\fr{3}{2}} \Phi\|_{L^2}^2 dt +\|\tilde{g}\|_{L^2}^2 dt + 2 ( \tilde{g}, q)_{L^2} dW.
\ee
Integrating in time from $0$ to $t$, we get 
\beg{align} \la{S1}
&\|q(t,w)\|_{L^2}^2 + \int_{0}^{t} \|\l^{\fr{1}{2}} q (s,w)\|_{L^2}^2 ds \nonumber
\\&\quad\le \|q_0\|_{L^2}^2 + \left(\|\l^{\fr{3}{2}} \Phi\|_{L^2}^2 + \|\tilde{g}\|_{L^2}^2  \right)t 
+ 2 \int_{0}^{t}  ( \tilde{g}, q)_{L^2} dW. 
\end{align}
We take the supremum over all $t \in [0,T]$,
\beg{align} 
&\sup\limits_{0 \le t \le T} \|q(w)\|_{L^2}^2 + \int_{0}^{T} \|\l^{\fr{1}{2}} q (s,w)\|_{L^2}^2 ds  \nonumber
\\&\quad\le 2\|q_0\|_{L^2}^2 + 2\left(\|\l^{\fr{3}{2}} \Phi\|_{L^2}^2 + \|\tilde{g}\|_{L^2}^2 \right)T 
+ 4 \sup\limits_{0 \le t \le T} \left|\int_{0}^{t}  ( \tilde{g}, q)_{L^2} dW\right|. 
\end{align} Now we apply the expectation $\E$. In view of the martingale estimate (see \cite{DPZ}), 
\be 
\E \left\{\sup\limits_{0 \le t \le T} \left|\int_{0}^{t}  ( \tilde{g}, q)_{L^2} dW \right| \right\}
\le C \E \left\{\left(\int_{0}^{T}  ( \tilde{g}, q)_{L^2}^2 dt \right)^{\fr{1}{2}}\right\},
\ee we have 
\beg{align}
&\E \left\{\sup\limits_{0 \le t \le T} \left|\int_{0}^{t}  ( \tilde{g}, q)_{L^2} dW \right| \right\}
\le C \E \left\{\left(\int_{0}^{T} \|q\|_{L^2}^2\|\tilde{g}\|_{L^2}^2 dt \right)^{\fr{1}{2}} \right\} \nonumber
\\&\quad\le \E \left\{\left(\sup\limits_{0 \le t \le T} \|q\|_{L^2} \right) \left(C\int_{0}^{T} \|\tilde{g}\|_{L^2}^2 dt \right)^{\fr{1}{2}} \right\} 
\le \fr{1}{8} \E \left\{\sup\limits_{0 \le t \le T} \|q\|_{L^2}^2 \right\} + C\|\tilde{g}\|_{L^2}^2 T
\end{align} 
This gives \eqref{SS2} when $p = 2$.

\textbf{Step 2.} We prove that the estimate \eqref{SS2} holds for any $p \in [4, \infty)$. 

\textbf{Proof of Step 2.} Applying It\^o's lemma to the process $F(X_t (w))$ where $X_t(w) = \|q(t,w)\|_{L^2}^2$ obeys \eqref{ql2} and $F(\xi) = \xi^{\fr{p}{2}}$, we derive the energy equality 
\beg{align}
\d (\|q\|_{L^2}^2)^{\fr{p}{2}} 
&= - p \|q\|_{L^2}^{p-2} \|\l^{\fr{1}{2}} q\|_{L^2}^2 dt \nonumber
\\&\quad+ p \|q\|_{L^2}^{p-2} (\Delta \Phi, q)_{L^2} dt
+ \fr{p}{2} \|q\|_{L^2}^{p-2}  \| \tilde{g}\|_{L^2}^2 dt \nonumber
\\&\quad +p \left(\fr{p}{2} -1 \right) \|q\|_{L^2}^{p-4}  |( \tilde{g}, q)_{L^2} |^2 dt
+  p\|q\|_{L^2}^{p-2} ( \tilde{g}, q)_{L^2} dW, 
\end{align} which yields the differential inequality
\beg{align}
&\d \|q\|_{L^2}^p + p \|q\|_{L^2}^{p-2} \|\l^{\fr{1}{2}} q\|_{L^2}^2 dt
\le p \|q\|_{L^2}^{p-1} \|\Delta \Phi\|_{L^2}dt \nonumber
\\&\quad\quad+ \fr{p}{2} (p-1) \|q\|_{L^2}^{p-2} \|\tilde{g}\|_{L^2}^2 dt
+   p\|q\|_{L^2}^{p-2} ( \tilde{g}, q)_{L^2} dW. 
\end{align} 
In view of the bound
\be 
\|q\|_{L^2} \le \|\l^{\fr{1}{2}} q\|_{L^2},
\ee we have
\beg{align}  \la{qpeq}
&\d \|q\|_{L^2}^p  + \fr{p}{4} \|q\|_{L^2}^p dt + \fr{p}{2} \|q\|_{L^2}^{p-2} \|\l^{\fr{1}{2}} q\|_{L^2}^2 dt \nonumber
\\&\quad\le C(p) \left(\|\Delta \Phi\|_{L^2}^p  + \|\tilde{g}\|_{L^2}^{p} \right)dt
+  p\|q\|_{L^2}^{p-2} ( \tilde{g}, q)_{L^2} dW
\end{align} where we used Young's inequality to estimate 
\be 
p \|q\|_{L^2}^{p-1} \|\Delta \Phi\|_{L^2} 
\le C(p) \|\Delta \Phi\|_{L^2}^p + \fr{p}{8} \|q\|_{L^2}^p
\ee and 
\be 
\fr{p}{2} (p-1) \|q\|_{L^2}^{p-2} \|\tilde{g}\|_{L^2}^2
\le C(p) \|\tilde{g}\|_{L^2}^{p} + \fr{p}{8} \|q\|_{L^2}^p.
\ee
Integrating in time \eqref{qpeq} from $0$ to $t$ and taking the supremum over $[0,T]$, we obtain 
\beg{align}
&\sup\limits_{0 \le t \le T} \|q\|_{L^2}^p + \fr{p}{2} \int_{0}^{T} \|q\|_{L^2}^{p-2} \|\l^{\fr{1}{2}} q\|_{L^2}^2 ds \nonumber
\\&\quad\le 2\|q_0\|_{L^2}^p  
+ C(p) \left(\|\Delta \Phi\|_{L^2}^p  + \|\tilde{g}\|_{L^2}^p \right) T 
+ 2\sup\limits_{0 \le t \le T} \left|\int_{0}^{t}  p \|q\|_{L^2}^{p-2} (\tilde{g}, q)_{L^2} dW \right|.
\end{align} We estimate 
\beg{align}
&\E \left\{\sup\limits_{0 \le t \le T} \left|\int_{0}^{t} 2p\|q\|_{L^2}^{p-2}  ( \tilde{g}, q)_{L^2} dW \right| \right\}
\le C(p) \E \left\{\left(\int_{0}^{T}  \|q\|_{L^2}^{2p-4} (\tilde{g}, q)_{L^2}^2 dt \right)^{\fr{1}{2}}\right\} \nonumber
\\&\quad\le C(p) \E \left\{\left(\int_{0}^{T} \|q\|_{L^2}^{2p-2} \|\tilde{g}\|_{L^2}^2 dt \right)^{\fr{1}{2}} \right\} 
\le \E \left\{\left(\sup\limits_{0 \le t \le T} \|q\|_{L^2}^{p-1} \right) \left(C(p) \int_{0}^{T} \|\tilde{g}\|_{L^2}^2 dt \right)^{\fr{1}{2}} \right\} \nonumber
\\&\quad\le \left( 1 - \fr{1}{p} \right) \E \left\{\sup\limits_{0 \le t \le T} \|q\|_{L^2}^p \right\} + C(p) \|\tilde{g}\|_{L^2}^p T^{\fr{p}{2}}
\end{align} and we obtain \eqref{SS2}.

\textbf{Step 3.} We show that the velocity $u$ obeys 
\beg{align} \la{SS3}
&\E \left\{\sup\limits_{0 \le t \le T} \|u^{\epsilon}\|_{L^2}^2 
+ \int_{0}^{T} \|\na u^{\epsilon}\|_{L^2}^2 dt \right\} 
\le  C(\|u_0\|_{L^2}, \|q_0\|_{L^2}, f, \Phi, \tilde{g}, g)e^{4T}.
\end{align}

\textbf{Proof of Step 3.} We apply It\^o's lemma pointwise in $x$ and we obtain the energy equality  
\beg{align}
\d \|u\|_{L^2}^2
&= -2(-\Delta u, u)_{L^2} dt 
-2 (u \cdot \na u, u)_{L^2} dt
-2 (q Rq, u)_{L^2}dt
-2 (q \na \Phi, u)_{L^2} dt \nonumber
\\&\quad\quad+2 (f, u)_{L^2} dt
+  \| g\|_{L^2}^2 dt
+ 2  ( g, u)_{L^2} dW,
\end{align} which implies
\beg{align}
&\d \|u\|_{L^2}^2 + 2\|\na u\|_{L^2}^2 dt  \nonumber
\\&\quad= -2 (q Rq + q \na  \Phi -f, u)_{L^2}dt
+  \| g\|_{L^2}^2 dt + 2  ( g, u)_{L^2} dW,
\end{align} where we used the cancellation
\be 
(u \cdot \na u, u)_{L^2} = 0
\ee due to the divergence-free condition satisfied by $u$. 
By Ladyzhenskaya's interpolation inequality 
\be 
\|u\|_{L^4} \le C\|u\|_{L^2} + C\|u\|_{L^2}^{\fr{1}{2}} \|\na u\|_{L^2}^{\fr{1}{2}},
\ee and the boundedness of the Riesz transforms in $L^4$, we estimate
\beg{align}
|(q Rq, u)_{L^2}| 
&\le \|q\|_{L^2} \|Rq\|_{L^4} \|u\|_{L^4} 
\le C\|q\|_{L^2} \|q\|_{L^4}\left(\|u\|_{L^2} + \|u\|_{L^2}^{\fr{1}{2}} \|\na u\|_{L^2}^{\fr{1}{2}} \right) \nonumber 
\\&\le C\|q\|_{L^2}^2 \|q\|_{L^4}^2 + \fr{1}{2} \|u\|_{L^2}^2 + \fr{1}{2}\|\na u\|_{L^2}^2.
\end{align} We also estimate 
\be 
|(q \na \Phi, u)_{L^2}| \le \fr{1}{2} \|u\|_{L^2}^2 + \fr{1}{2} \|\na \Phi\|_{L^4}^2 \|q\|_{L^4}^2 
\ee and 
\be 
|(f, u)_{L^2}| \le \fr{1}{2} \|u\|_{L^2}^2 + \fr{1}{2} \|f\|_{L^2}^2
\ee using H\"older's inequality followed by Young's inequality. 
We obtain the differential inequality 
\beg{align} 
\d \|u\|_{L^2}^2 + \|\na u\|_{L^2}^2 dt
&\le 3\|u\|_{L^2}^2 dt 
+ \|f\|_{L^2}^2 dt 
+ C\|q\|_{L^2}^2 \|q\|_{L^4}^2 dt \nonumber
\\&\quad\quad+ C \|\na \Phi\|_{L^4}^2\|q\|_{L^4}^2 dt
+ \|g\|_{L^2}^2 dt
+ 2 ( g, u)_{L^2} dW,
\end{align} hence
\beg{align}
&\d \left\{e^{-3t} \|u\|_{L^2}^2 \right\}(s) 
= -3e^{-3s} \|u\|_{L^2}^2 ds + e^{-3s} \d \|u (s)\|_{L^2}^2  \nonumber
\\&\quad\le -e^{-3s}\|\na u\|_{L^2}^2  ds + e^{-3s}  \left\{  \|f\|_{L^2}^2 ds + C\|q\|_{L^2}^2 \|q\|_{L^4}^2 ds + C \|\na \Phi\|_{L^4}^2\|q\|_{L^4}^2 ds \right\} \nonumber
\\&\quad\quad+ e^{-3s}\|g\|_{L^2}^2 ds + 2e^{-3s}  ( g, u)_{L^2} dW
\end{align} for all $s \in [0,t]$. Integrating in time from $0$ to $t$, we obtain 
\beg{align} \la{S3}
&e^{-3t} \|u(t)\|_{L^2}^2 + \int_{0}^{t} e^{-3s} \|\na u(s)\|_{L^2}^2 ds 
\le \|u_0\|_{L^2}^2  
+ \left(\|f\|_{L^2}^2 + \|g\|_{L^2}^2\right)t   \nonumber
\\&\quad\quad+ C \int_{0}^{t} \|q(s)\|_{L^2}^2 \|q(s)\|_{L^4}^2 ds 
+ C \int_{0}^{t} \|\na \Phi\|_{L^4}^2\|q(s)\|_{L^4}^2 ds \nonumber
\\&\quad\quad+ 2 \int_{0}^{t} e^{-3s}  ( g, u)_{L^2} dW(s).
\end{align}
We take the supremum in time over $[0,T]$ and apply $\E$. Using the continuous Sobolev embedding
\be 
H^{\fr{1}{2}} (\TT^2) \subset L^4(\TT^2)
\ee and \eqref{SS2} with $p=4$, we have 
\be 
\E \left\{\int_{0}^{T} \|q(s)\|_{L^2}^2 \|q(s)\|_{L^4}^2 ds  \right\} 
\le C\|q_0\|_{L^2}^4 + C \left(\|\Delta \Phi\|_{L^2}^4 + \|\tilde{g}\|_{L^2}^{2} \right)T + C\|\tilde{g}\|_{L^2}^2 T^2
\ee for all $t \in [0,T]$. From \eqref{SS2} with $p=2$, we have 
\be 
\E \left\{\int_{0}^{T} \|\na \Phi\|_{L^4}^2\|q(s)\|_{L^4}^2 ds \right\}
\le C\|\na \Phi\|_{L^4}^2 \left( \|q_0\|_{L^2}^2 + \|\l^{\fr{3}{2}} \Phi\|_{L^2}^2T + \|\tilde{g}\|_{L^2}^2 T \right)
\ee for all $t \in [0,T]$. We estimate 
\beg{align}
&\E \left\{\sup\limits_{0 \le t \le T}  \left|\int_{0}^{t}  2e^{-3s} ( g, u)_{L^2} dW \right| \right\}
\le  \E \left\{\sup\limits_{0 \le t \le T} \left(e^{-\fr{3}{2}t} \|u(t)\|_{L^2}\right) \left(\int_{0}^{T} Ce^{-3t} \|g\|_{L^2}^2 dt \right)^{\fr{1}{2}} \right\} \nonumber
\\&\quad\le \fr{1}{2} \E \left\{ \sup\limits_{0 \le t \le T} \left(e^{-3t}  \|u(t)\|_{L^2}^2 \right) \right\} + C\|g\|_{L^2}^2 
\end{align} and we obtain \eqref{SS3}.

\textbf{Step 4.} We prove that \eqref{SS5} holds for $p=4$.

\textbf{Proof of Step 4.} By It\^o's lemma, we have
\beg{align}
\d |q|^4 
&= -4 q^3 u \cdot \na q dt
- 4 q^3 \l q dt
+ 4 q^3 \Delta \Phi dt \nonumber
\\&\quad\quad+ 6  q^2 \tilde{g}^2 dt
+ 4q^3   \tilde{g} d W.  
\end{align}
Integrating in the space over $\TT^2$, we obtain the energy equality
\beg{align}  \la{1}
\d \|q\|_{L^4}^4 
&= -4 (u \cdot \na q, q^3)_{L^2} dt
- 4 (\l q, q^3)_{L^2} dt
+ 4 (\Delta \Phi, q^3)_{L^2} dt \nonumber
\\&+ 6( ( \tilde{g})^2, q^2)_{L^2} dt
+ 4  ( \tilde{g}, q^3)_{L^2} dW.
\end{align} We note that 
\be \la{2}
(u \cdot \na q, q^3)_{L^2} = 0
\ee due to the divergence-free condition for $u$. 
By the nonlinear Poincar\'e inequality for the fractional Laplacian in $L^4$ applied to the mean zero function $q$, we have 
\be  \la{3}
\int_{\TT^2} q^3 \l q dx \ge c \|q\|_{L^4}^4.
\ee Using H\"older's inequality with exponents $4, 4/3$ and Young's inequality with exponents $4, 4/3$, we get 
\be  \la{4}
4|(\Delta \Phi, q^3)_{L^2}|
\le 4\|\Delta \Phi\|_{L^{4}} \|q^3\|_{L^{4/3}}
= 4\|\Delta \Phi\|_{L^4} \|q\|_{L^4}^3
\le c\|q\|_{L^4}^4 + C\|\Delta \Phi\|_{L^4}^4.
\ee We also bound
\be  \la{5}
6|( \left( \tilde{g})^2, q^2 \right)_{L^2}|
\le 6\|q\|_{L^4}^2  \|\tilde{g}\|_{L^4}^2
\le c\|q\|_{L^4}^4 + C \|\tilde{g}\|_{L^4}^4,
\ee using H\"older and Young inequalities. 
Putting \eqref{1}--\eqref{5} together, we obtain the differential inequality 
\be \la{S4}
\d \|q\|_{L^4}^4 + c \|q\|_{L^4}^4 dt
\le C\|\Delta \Phi\|_{L^4}^4 dt + C \|\tilde{g}\|_{L^4}^4 dt
+ 4  ( \tilde{g}, q^3)_{L^2} dW. 
\ee Consequently,
\be 
\|q(t)\|_{L^4}^4 
+ c\int_{0}^{t} \|q\|_{L^4}^4 ds 
\le 2\|q_0\|_{L^4}^4  
+ C\|\Delta \Phi\|_{L^4}^4 t
+ C \|\tilde{g}\|_{L^4}^4 t
+ 4 \int_{0}^{t}  ( \tilde{g}, q^3)_{L^2} dW
\ee for all $t \in [0,T]$. We take the supremum over $[0,T]$ and then we apply $\E$. We estimate 
\beg{align} \la{martingaleterm2}
&\E \left\{\sup\limits_{0 \le t \le T} \left|8\int_{0}^{t}  ( \tilde{g}, q^3)_{L^2} d W \right| \right\}
\le C\E \left\{\left(\int_{0}^{T}  (\tilde{g}, q^3)_{L^2}^2 dt \right)^{\fr{1}{2}} \right\} \nonumber
\\&\quad\le C\E \left\{\left(\int_{0}^{T}  \|\tilde{g} \|_{L^4}^2 \|q^3\|_{L^{4/3}}^2 dt \right)^{\fr{1}{2}} \right\}
\le \E \left\{\sup\limits_{0 \le t \le T} \|q\|_{L^4}^3 \left(C\int_{0}^{T}  \|\tilde{g}\|_{L^4}^2 dt \right)^{\fr{1}{2}} \right\} \nonumber
\\&\quad\le \fr{3}{4} \E \left\{\sup\limits_{0 \le t \le T} \|q\|_{L^4}^4 \right\} + C  \|\tilde{g} \|_{L^4}^4  T^2
\end{align}
and we obtain \eqref{SS5} for $p=4$.

\textbf{Step 5.} We prove \eqref{SS5} for any $p \ge 8$. 

\textbf{Proof of Step 5.} The stochastic energy equality
\beg{align}  
\d \|q\|_{L^4}^p 
&= - p\|q\|_{L^4}^{p-4} (\l q, q^3)_{L^2} dt
+ p\|q\|_{L^4}^{p-4}  (\Delta \Phi, q^3)_{L^2} dt \nonumber
\\&\quad\quad+ \fr{3}{2} p \|q\|_{L^4}^{p-4}( \tilde{g}^2, q^2)_{L^2} dt 
+ 2p \left(\fr{p}{4} -1 \right)\|q\|_{L^4}^{p-8}  ( \tilde{g}, q^3)_{L^2}^2 dt \nonumber
\\&\quad\quad+ p \|q\|_{L^4}^{p-4}  ( \tilde{g}, q^3)_{L^2} dW
\end{align} holds for any $p \ge 8$.
By H\"older's inequality with exponents $4/3, 4$ and Young's inequality with exponents $p/(p-2), p/2$, we have 
\beg{align} 
&2p \left(\fr{p}{4} -1 \right)\|q\|_{L^4}^{p-8}  ( \tilde{g}, q^3)_{L^2}^2
\le 2p \left(\fr{p}{4} -1 \right)\|q\|_{L^4}^{p-8}  \|q^3\|_{L^{4/3}}^2  \|\tilde{g} \|_{L^4}^2 \nonumber
\\&\quad= 2p \left(\fr{p}{4} -1 \right) \|q\|_{L^4}^{p-8} \|q\|_{L^4}^6  \|\tilde{g}\|_{L^4}^2 
\le \fr{cp}{8} \|q\|_{L^4}^{p} + C\left( \|\tilde{g} \|_{L^4}^2 \right)^{\fr{p}{2}}.
\end{align}
We obtain 
\be \la{S5}
\d \|q\|_{L^4}^p  + \fr{cp}{2} \|q\|_{L^4}^p  dt
\le C\|\Delta \Phi\|_{L^4}^p dt + C\left( \|\tilde{g} \|_{L^4}^2 \right)^{\fr{p}{2}}  dt
+ p \|q\|_{L^4}^{p-4}  ( \tilde{g}, q^3)_{L^2} dW.
\ee Integrating \eqref{S5} in time from $0$ to $t$, taking the supremum over $[0,T]$, applying $\E$, and estimating
\beg{align} 
&\E \left\{\sup\limits_{0 \le t \le T} 2p\left|\int_{0}^{t} \|q\|_{L^4}^{p-4}  ( \tilde{g}, q^3)_{L^2} dW \right| \right\}\nonumber\\ 
&\quad\le \left(1 - \fr{1}{p} \right) \E \left\{\sup\limits_{0 \le t \le T} \|q\|_{L^4}^p \right\} + C(p) \|\tilde{g}\|_{L^4}^p T^{\fr{p}{2}} 
\end{align} we obtain \eqref{SS5}.

\textbf{Step 6.} We show that \eqref{SS6} holds. 

\textbf{Proof of Step 6.} We derive the stochastic energy equality
\beg{align}
\d (\|u\|_{L^2}^2)^{\fr{p}{2}}
&= -p \|u\|_{L^2}^{p-2} \|\na u\|_{L^2}^2 dt
+ p \|u\|_{L^2}^{p-2} (-q Rq - q \na \Phi + f, u)_{L^2}dt \nonumber
\\&\quad\quad+ \fr{p}{2}\|u\|_{L^2}^{p-2}  \|g\|_{L^2}^2 dt
+ p \left(\fr{p}{2} -1 \right) \|u\|_{L^2}^{p-4}  |( g, u)_{L^2}|^2 dt \nonumber
\\&\quad\quad+ p\|u\|_{L^2}^{p-2}  ( g, u)_{L^2} dW.
\end{align}
By Young's inequality with exponents $p/(p-2)$ and $p/2$, 
\be 
\fr{p}{2} \|u\|_{L^2}^{p-2}  \| g \|_{L^2}^2
\le  \fr{1}{5} \|u\|_{L^2}^p + C(p) \|g\|_{L^2}^p 
\ee and 
\beg{align} 
p\left(\fr{p}{2} -1 \right)\|u\|_{L^2}^{p-4}  |( g, u)_{L^2}|^2
&\le p\left(\fr{p}{2} -1 \right)\|u\|_{L^2}^{p-4} \|u\|_{L^2}^2 \|g\|_{L^2}^2 \nonumber
\\&\le \fr{1}{5}  \|u\|_{L^2}^p + C(p) \|g\|_{L^2}^p.
\end{align} Similarly, using Young's inequality with exponents $p/(p-1)$ and $p$, 
\be 
p\|u\|_{L^2}^{p-2} |(f, u)_{L^2}| \le p\|u\|_{L^2}^{p-2} \|u\|_{L^2}\|f\|_{L^2}
\le C(p)\|f\|_{L^2}^p + \fr{1}{5}  \|u\|_{L^2}^p
\ee and 
\beg{align} 
p\|u\|_{L^2}^{p-2}|(q \na  \Phi, u)_{L^2}| 
&\le p\|u\|_{L^2}^{p-2} \|u\|_{L^2} \|q\|_{L^2} \|\na \Phi\|_{L^{\infty}} \nonumber
\\&\le C(p) \|\na \Phi\|_{L^{\infty}}^p \|q\|_{L^2}^p + \fr{1}{5} \|u\|_{L^2}^p.
\end{align} 
By Ladyzhenskaya's interpolation inequality and the boundedness of the Riesz transforms in $L^4(\TT^2)$, we have 
\beg{align} 
&p\|u\|_{L^2}^{p-2} |(-q Rq, u)_{L^2}|
\le C(p)\|u\|_{L^2}^{p-2} \|u\|_{L^4}\|q\|_{L^2}\|q\|_{L^4} \nonumber
\\&\quad\le C(p)\|u\|_{L^2}^{p-2} \left(\|u\|_{L^2} + \|u\|_{L^2}^{\fr{1}{2}}\|\na u\|_{L^2}^{\fr{1}{2}} \right) \|q\|_{L^2}\|q\|_{L^4} \nonumber
\\&\quad\le \|u\|_{L^2}^{p} + \fr{p}{2} \|u\|_{L^2}^{p-2}\|\na u\|_{L^2}^2 + C(p) \|q\|_{L^2}^p \|q\|_{L^4}^p \nonumber
\\&\quad\le \fr{1}{5}  \|u\|_{L^2}^{p} + \fr{p}{2} \|u\|_{L^2}^{p-2}\|\na u\|_{L^2}^2 + C(p) \|q\|_{L^2}^{2p} +  C(p) \|q\|_{L^4}^{2p}.
\end{align} This yields the differential inequality
\beg{align} 
&\d \|u\|_{L^2}^{p} + \fr{p}{2} \|u\|_{L^2}^{p-2} \|\na u\|_{L^2}^2 dt
\le \|u\|_{L^2}^p dt 
+ C(p) \|g\|_{L^2}^p dt 
+ C(p) \|f\|_{L^2}^p dt \nonumber
\\&\quad\quad+ C(p) \|\na \Phi\|_{L^{\infty}}^p \|q\|_{L^2}^p dt 
+ C(p) \|q\|_{L^2}^{2p} dt + C(p) \|q\|_{L^4}^{2p} dt
+ p\|u\|_{L^2}^{p-2}  ( g, u)_{L^2} dW
\end{align} and thus
\beg{align}
&\d \left\{e^{-t} \|u\|_{L^2}^p \right\} (s) + e^{-s} \|u\|_{L^2}^{p-2} \|\na u\|_{L^2}^2 ds \nonumber
\\&\quad\le e^{-s} \left\{C(p) \|g\|_{L^2}^{p}ds
+ C(p) \|f\|_{L^2}^p ds
+ C(p) \|\na \Phi\|_{L^{\infty}}^p \|q\|_{L^2}^p ds
+ C(p) \|q\|_{L^2}^{2p} ds + C(p) \|q\|_{L^4}^{2p} ds \right\} \nonumber
\\&\quad\quad+ pe^{-s} \|u\|_{L^2}^{p-2}  ( g, u)_{L^2} dW.
\end{align}
We integrate in time from $0$ to $t$, take the supremum over $[0,T]$, and apply $E$. We obtain 
\beg{align} \la{mod1}
&\E \left\{\sup\limits_{0 \le t \le T} \left(e^{-t} \|u(t)\|_{L^2}^p \right) \right\} 
+ \E \left\{\int_{0}^{T} e^{-t} \|u\|_{L^2}^{p-2} \|\na u\|_{L^2}^2 dt  \right\}\nonumber
\\&\quad\le C(p) \left(\|g\|_{L^2}^p + \|f\|_{L^2}^p \right)
+ C(p) \|\na \Phi\|_{L^{\infty}}^p \E \left\{ \int_{0}^{T} \|q\|_{L^2}^p dt \right\}
+ C(p) \E \left\{\int_{0}^{T} \|q\|_{L^2}^{2p} dt \right\} \nonumber
\\&\quad\quad+ C(p) \E \left\{\int_{0}^{T} \|q\|_{L^4}^{2p} dt  \right\}
+ \sup\limits_{0 \le t \le T} \left|\int_{0}^{t} 2pe^{-s} \|u\|_{L^2}^{p-2}  ( g, u)_{L^2} dW \right|.
\end{align}
We estimate  
\beg{align} \la{mod2}
&\E \left\{\sup\limits_{0 \le t \le T}  \left|\int_{0}^{t} 2pe^{-s} \|u\|_{L^2}^{p-2}  (q_l, u)_{L^2} dW(s) \right| \right\} \nonumber
\\&\quad\le \left(1 - \fr{1}{p} \right) \E \left\{\sup\limits_{0 \le t \le T} \left(e^{-t}\|u(t)\|_{L^2}^p  \right) \right\}
+ C(p) \|g\|_{L^2}^p T^{\fr{p}{2}}.
\end{align}
Putting \eqref{mod1} and \eqref{mod2} together, and using \eqref{SS2} and \eqref{SS5}, we obtain \eqref{SS6}.

\textbf{Step 7.} We prove that \eqref{SS7} holds. 

\textbf{Proof of Step 7.} We write the equation satisfied by $\na u$, apply It\^o's lemma, and integrate in the space variable. We obtain the energy equality 
\beg{align}
&\d \|\na u\|_{L^2}^2  + 2\|\Delta u\|_{L^2}^2
= 2(u \cdot \na u, \Delta u)_{L^2} dt
+2 (qRq, \Delta u)_{L^2} dt \nonumber
\\&\quad\quad+ 2 (q \na \Phi, \Delta u)_{L^2} dt
-2 (f, \Delta u)_{L^2} dt
+ \| \na g \|_{L^2}^2 dt
- 2 (g, \Delta u)_{L^2} dW.
\end{align} The nonlinear term for the velocity vanishes, that is  
\be 
(u \cdot \na u, \Delta u)_{L^2} = 0,
\ee and using H\"older's inequality, we obtain 
\beg{align} \la{expna}
&\d \|\na u\|_{L^2}^2 + 2 \|\Delta u\|_{L^2}^2 dt
\le C\|q\|_{L^4}^2 \|\Delta u\|_{L^2} dt
+ 2\|\na \Phi\|_{L^{\infty}} \|q\|_{L^2} \|\Delta u\|_{L^2} dt \nonumber
\\&\quad\quad+ 2\|f\|_{L^2} \|\Delta u\|_{L^2} dt
+ \|\na g \|_{L^2}^2 dt
- 2 (g, \Delta u)_{L^2} dW.
\end{align}
An application of Young's inequality yields the differential inequality 
\beg{align} \la{naudiff}
&\d \|\na u\|_{L^2}^2 +  \|\Delta u\|_{L^2}^2 dt
\le C\|q\|_{L^4}^4  dt
+ C\|\na \Phi\|_{L^{\infty}}^2 \|q\|_{L^2}^2  dt \nonumber
\\&\quad\quad+ C\|f\|_{L^2}^2 dt
+ \|\na g\|_{L^2}^2 dt
- 2 (g, \Delta u)_{L^2} dW.
\end{align}
We integrate \eqref{naudiff} in time from $0$ to $t$, take the supremum in time, and then apply $\E$. We obtain
\beg{align} \la{naudiffmod}
&\E\left\{\sup\limits_{0 \le t \le T} \|\na u\|_{L^2}^2 \right\} 
+ \E \left\{\int_{0}^{T} \|\Delta u\|_{L^2}^2 dt \right\}
\le 2\|\na u_0\|_{L^2}^2
+ C\left(\|\na g\|_{L^2}^2 + \|f\|_{L^2}^2 \right)T \nonumber
\\&\quad+ C\E\left\{\int_{0}^{T} \|q\|_{L^4}^4 dt \right\}
+ C \|\na \Phi\|_{L^{\infty}}^2 \E \left\{\int_{0}^{T}  \|q\|_{L^2}^2  dt \right\}
+ \sup\limits_{0 \le t \le T} \left|\int_{0}^{t} 4 (g, \Delta u)_{L^2} dW \right|.
\end{align} We estimate the martingale term 
\beg{align} \la{martingaleterm}
&\E \left\{\sup\limits_{0 \le t \le T} \left|4 \int_{0}^{t} (g, \Delta u)_{L^2} dW \right| \right\} 
\le \E \left\{4\left(\int_{0}^{T}  (g, \Delta u)_{L^2}^2 dt \right)^{\fr{1}{2}} \right\} \nonumber
\\&\quad\le \E \left\{4\left(\int_{0}^{T} \|\na g\|_{L^2}^2 \|\na u\|_{L^2}^2 dt \right)^{\fr{1}{2}} \right\} 
\le \E \left\{4\sup\limits_{0 \le t \le T} \|\na u\|_{L^2} \left(\int_{0}^{T} \|\na g\|_{L^2}^2 dt \right)^{\fr{1}{2}} \right\} \nonumber
\\&\quad\le \fr{1}{2} \sup\limits_{0 \le t \le T} \E \left\{\|\na u\|_{L^2}^2 \right\} + C\|\na g\|_{L^2}^2 T.
\end{align}
Putting \eqref{naudiffmod} and \eqref{martingaleterm} together, and using \eqref{SS2} with $p=2$ and \eqref{SS5} with $p=4$, we get \eqref{SS7}.

\section{Pathwise Uniform Bounds for the Solutions} \la{B}

In this section, we prove Proposition \ref{pathneeded}. We let $(q,u)$ be the solution to \eqref{invsys} corresponding to the initial data $(q_0, u_0)$. Let \be  \la{noisedet1}
\tilde{\phi} (x,t,w) = \int_{0}^{t} e^{(t-s)\l} \tilde{g}(x) dW
\ee and 
\be  \la{noisedet2}
\phi (x,t,w) = \int_{0}^{t} e^{-(t-s)\Delta} g(x) dW.
\ee We set 
\be 
Q = q - \tilde{\phi}
\ee and 
\be 
U = u- \phi
\ee and we note that $(Q,U)$ obeys the deterministic system 
\be 
\begin{cases} \la{invsys3}
\pa_t Q + (U + \phi) \cdot \na (Q + \tilde{\phi})  + \l Q   = -  \l\tilde{\phi} 
\\ \pa_t U + (U + \phi) \cdot \na (U + \phi)  - \Delta U  + \na P 
= - (Q + \tilde{\phi}) R(Q + \tilde{\phi})   + f  + \Delta \phi
\\\na \cdot U = 0
\end{cases} 
\ee where we used the divergence-free condition imposed on $g$.

\textbf{Step 1. Bounds for the velocity in $L_{loc}^2(0,\infty;H^1(\TT^2))$.}
We take the $L^2$ inner product of the $Q$ equation with $Q$, and we obtain 
\be 
\fr{1}{2} \fr{d}{dt} \|Q\|_{L^2}^2 + \|\l^{\fr{1}{2}} Q\|_{L^2}^2
= - \int_{\TT^2} ((U + \phi) \cdot \na (Q + \tilde{\phi})) Q dx
- \int_{\TT^2} \l \tilde{\phi} Q dx.
\ee
We estimate the nonlinear term 
\beg{align} 
&\left|\int_{\TT^2} ((U + \phi) \cdot \na (Q + \tilde{\phi})) Q dx \right|
= \left|\int_{\TT^2} ((U + \phi) \cdot \na  \tilde{\phi})) Q dx  \right| \nonumber
\\&\le \|Q\|_{L^4} \|U\|_{L^4}\|\na \tilde{\phi} \|_{L^2} + \|Q\|_{L^4} \|\phi\|_{L^4}\|\na \tilde{\phi} \|_{L^2} \nonumber
\\&\le C\|\l^{\fr{1}{2}} Q\|_{L^2} \|U\|_{L^2}^{\fr{1}{2}} \|\na U\|_{L^2}^{\fr{1}{2}} \|\na \tilde{\phi} \|_{L^2} + C\|\l^{\fr{1}{2}} Q\|_{L^2} \|\phi\|_{L^4}\|\na \tilde{\phi} \|_{L^2} \nonumber
\\&\le \fr{1}{8} \|\l^{\fr{1}{2}} Q\|_{L^2}^2 + \fr{1}{8} \|\na U\|_{L^2}^2 + C\|\na \tilde{\phi} \|_{L^2}^4 \|U \|_{L^2}^2 + C\|\phi\|_{L^4}^2\|\na \tilde{\phi} \|_{L^2}^2
\end{align} using H\"older's inequality, Ladyzhenskaya's interpolation inequality applied to the mean zero function $U$, the continuous Sobolev embedding $H^{\fr{1}{2}} \subset L^4$, and Young's inequality. 
This yields the differential inequality
\be  \la{step11}
\fr{1}{2} \fr{d}{dt} \|Q\|_{L^2}^2 +  \fr{3}{4}\|\l^{\fr{1}{2}} Q\|_{L^2}^2
 \le \fr{1}{8} \|\na U\|_{L^2}^2 
 + C\|\na \tilde{\phi} \|_{L^2}^4 \|U \|_{L^2}^2 
 + C\|\phi\|_{L^4}^2\|\na \tilde{\phi} \|_{L^2}^2
 + C\|\l^{\fr{1}{2}} \tilde{\phi}\|_{L^2}^2.
\ee
Now we take the $L^2$ inner product of the $Q$ equation with $\l^{-1}Q$ and we get
\be 
\fr{1}{2} \fr{d}{dt} \|\l^{-\fr{1}{2}} Q\|_{L^2}^2 + \|Q\|_{L^2}^2
= - \int\limits_{\TT^2} \l\tilde{\phi} \l^{-1} Q dx
- \int\limits_{\TT^2} ((U + \phi) \cdot \na (Q + \tilde{\phi})) \l^{-1} Q dx
\ee
Integrating by parts and using the divergence-free condition obeyed by $U + \phi$, we can rewrite the nonlinear term as
\beg{align}
&- \int\limits_{\TT^2} ((U + \phi) \cdot \na (Q + \tilde{\phi})) \l^{-1} Q dx \nonumber
\\&= \int\limits_{\TT^2} (Q + \tilde{\phi}) R(Q + \tilde{\phi}) \cdot (U + \phi) dx
- \int\limits_{\TT^2} (Q + \tilde{\phi}) R \tilde{\phi} \cdot (U + \phi) dx \nonumber
\\&= \int\limits_{\TT^2} (Q + \tilde{\phi}) R(Q + \tilde{\phi}) \cdot U  dx
+ \int\limits_{\TT^2} (Q + \tilde{\phi}) \left[R(Q + \tilde{\phi}) \cdot \phi - R \tilde{\phi} \cdot (U + \phi) \right] dx  \nonumber
\\&= \int\limits_{\TT^2} (Q + \tilde{\phi}) R(Q + \tilde{\phi}) \cdot U  dx
+ \int\limits_{\TT^2} (Q + \tilde{\phi}) \left[RQ  \cdot \phi - R \tilde{\phi} \cdot U  \right] dx 
\end{align} and we estimate
\beg{align}
&\left|\int\limits_{\TT^2} (Q + \tilde{\phi}) \left[RQ  \cdot \phi - R \tilde{\phi} \cdot U  \right] dx   \right| \nonumber
\\&\le\left| \int\limits_{\TT^2} Q RQ \cdot \phi dx \right|
+ \left| \int\limits_{\TT^2} \tilde{\phi} RQ \cdot \phi dx \right|
+ \left| \int\limits_{\TT^2} Q R \tilde{\phi} \cdot U dx \right|
+ \left| \int\limits_{\TT^2} \tilde{\phi}R \tilde{\phi} \cdot U  dx \right| \nonumber
\\&\le C\|\phi\|_{L^4} \|Q\|_{L^2} \|\l^{\fr{1}{2}} Q \|_{L^2}
+ C\|\phi\|_{L^4} \|\tilde{\phi}\|_{L^4} \|Q\|_{L^2} 
+ C\|Q\|_{L^4} \|\tilde{\phi}\|_{L^4} \|U\|_{L^2}
+ C\|\tilde{\phi}\|_{L^4} \|\tilde{\phi} \|_{L^4} \|U\|_{L^2}  \nonumber
\\&\le \fr{1}{4} \|\l^{\fr{1}{2}} Q\|_{L^2}^2 
+ \fr{1}{8} \|Q\|_{L^2}^2 
+ C\|\phi\|_{L^4}^2 \|Q\|_{L^2}^2 
+ C(\|\phi\|_{L^4}^2 + \|\tilde{\phi}\|_{L^4}^2) \|\tilde{\phi} \|_{L^4}^2
+ C(1 + \|\tilde{\phi}\|_{L^4}^2 )\|U\|_{L^2}^2
\end{align} where we have used the boundedness of the Riesz transforms on $L^p(\TT^2)$ for $p \in (1, \infty)$. 
We obtain  
\beg{align} \la{step12}
&\fr{1}{2} \fr{d}{dt} \|\l^{-\fr{1}{2}} Q\|_{L^2}^2
+ \fr{3}{4} \|Q\|_{L^2}^2
\le C\| \tilde{\phi}\|_{L^2}^2
+ \int\limits_{\TT^2} (Q + \tilde{\phi}) R(Q + \tilde{\phi}) \cdot U  dx \nonumber
\\&+ \fr{1}{4} \|\l^{\fr{1}{2}} Q\|_{L^2}^2 
+ C\|\phi\|_{L^4}^2 \|Q\|_{L^2}^2 
+ C(\|\phi\|_{L^4}^2 + \|\tilde{\phi}\|_{L^4}^2) \|\tilde{\phi} \|_{L^4}^2
+ C(1 + \|\tilde{\phi}\|_{L^4}^2 )\|U\|_{L^2}^2
\end{align}
Finally, we take the $L^2$ inner product of the equation obeyed by $U$ with $U$ and we obtain 
\beg{align} 
\fr{1}{2} \fr{d}{dt} \|U\|_{L^2}^2 
+ \|\na U\|_{L^2}^2
&= - \int\limits_{\TT^2} (U + \phi) \cdot \na (U + \phi) \cdot U dx
- \int\limits_{\TT^2} (Q + \tilde{\phi}) R (Q + \tilde{\phi}) \cdot U dx \nonumber
\\& \int\limits_{\TT^2} \Delta \phi \cdot U dx
+ \int\limits_{\TT^2} f \cdot U dx
\end{align}
We integrate by parts the nonlinear term. Using the fact that $U+ \phi$ is divergence-free, we have  
\beg{align} 
&\left|\int\limits_{\TT^2} (U + \phi) \cdot \na (U + \phi) \cdot U dx  \right|
= \left| \int_{\TT^2} ((U + \phi) \cdot \na \phi ) \cdot U dx \right| \nonumber
\\&\le \|U\|_{L^2} \|U\|_{L^4}\|\na \phi\|_{L^4}
+ \|U\|_{L^2} \|\phi\|_{L^4} \|\na \phi\|_{L^4} \nonumber
\\&\le C(\|U\|_{L^2}^{\fr{1}{2}}\|\na U\|_{L^2}^{\fr{1}{2}} \|\na \phi\|_{L^4}  + \|\phi\|_{L^4}\|\na \phi\|_{L^4})\| U\|_{L^2} \nonumber
\\&\le \fr{1}{16} \|\na U\|_{L^2}^2 + C (\|\na \phi\|_{L^4}^{\fr{4}{3}} +1) \|U\|_{L^2}^2 + C \|\phi\|_{L^4}^2 \|\na \phi\|_{L^4}^2.
\end{align} 
This yields the differential inequality
\beg{align}  \la{step13}
&\fr{1}{2}\fr{d}{dt} \|U\|_{L^2}^2 + \fr{3}{4} \|\na U\|_{L^2}^2 
\le - \int\limits_{\TT^2} (Q + \tilde{\phi}) R (Q + \tilde{\phi}) \cdot U dx 
+ \fr{1}{4} \|Q\|_{L^2}^2  \nonumber
\\&+  C(\|\na \phi\|_{L^4}^{\fr{4}{3}} + 1) \|U\|_{L^2}^2 
+ C\|\tilde{\phi}\|_{L^2}^2 
+ C\|f\|_{L^2}^2
+ C\|\na \phi\|_{L^2}^2
+ C\|\phi\|_{L^4}^2 \|\na \phi\|_{L^4}^2
\end{align}
We add \eqref{step11}, \eqref{step12} and \eqref{step13}. Setting 
\be 
X (t,w) = \|Q\|_{L^2}^2 + \|\l^{-\fr{1}{2}} Q\|_{L^2}^2 + \|U\|_{L^2}^2, 
\ee  we get
\be 
 \fr{d}{dt} X +  \|\na U\|_{L^2}^2 \le CA(t) X(t) + CB(t)
\ee where $A(t)$ and $B(t)$ are some positive constants depending on $\phi, \tilde{\phi}$ and $f$. This implies
\be 
\fr{d}{dt} \left[e^{-C\int_{0}^{t} A(s) ds} X (t) \right]
+ e^{-C\int_{0}^{t} A(s) ds} \|\na U\|_{L^2}^2
\le C B(t).
\ee Integrating in time from $0$ to $t$, we obtain the bound
\be \la{lasts1}
X(t) + \int_{0}^{t} \|\na U\|_{L^2}^2 \le \left[C \int_{0}^{t} B(s) ds + 2\|q_0\|_{L^2}^2 + \|u_0\|_{L^2}^2 \right] e^{C\int_{0}^{t} A(s) ds}
\ee  for all $t \ge 0$.

\textbf{Step 2. Bounds for the charge density in $L_{loc}^{\infty}(0, \infty; L^4(\TT^2))$.}
We take the $L^2$ inner product of the $Q$ equation with $(Q)^3$. Using the Poincar\'e inequality for the fractional Laplacian, we get the deterministic differential inequality 
\be 
\fr{1}{4} \fr{d}{dt} \|Q\|_{L^4}^4 + c\|Q\|_{L^4}^4
\le - \int\limits_{\TT^2} \l \tilde{\phi} (Q)^3 dx
- \int\limits_{\TT^2} (U + \phi) \cdot \na (Q + \tilde{\phi}) (Q)^3 dx.
\ee 
In view of the continuous Sobolev embedding of $H^1(\TT^2)$ in $L^8(\TT^2)$, we bound the nonlinear term 
\beg{align}
&\left|\int\limits_{\TT^2} (U + \phi) \cdot \na (Q + \tilde{\phi}) (Q)^3 dx \right|
= \left| \int\limits_{\TT^2} ((U + \phi) \cdot \na  \tilde{\phi}) (Q)^3 dx \right| \nonumber
\\&\le \|U\|_{L^8} \|Q\|_{L^4}^3 \|\na \tilde{\phi}\|_{L^{8}}
+ \|Q\|_{L^4}^3 \|\phi \cdot \na \tilde{\phi}\|_{L^4} \nonumber
\\&\le C\|\na U\|_{L^2} \|Q\|_{L^4}^3 \|\na \tilde{\phi}\|_{L^{8}}
+ \|Q\|_{L^4}^3 \|\phi \cdot \na \tilde{\phi}\|_{L^4},
\end{align} hence
\be 
\fr{1}{4} \fr{d}{dt} \|Q\|_{L^4}^4 + c\|Q\|_{L^4}^4
\le \left[\|\l \tilde{\phi}\|_{L^4} +C \|\na U\|_{L^2} \|\na \tilde{\phi}\|_{L^{8}} +    \|\phi \cdot \na \tilde{\phi}\|_{L^4}  \right] \|Q\|_{L^4}^3
\ee which yields
\be 
\fr{d}{dt} \|Q\|_{L^4} + c\|Q\|_{L^4}
\le \|\l \tilde{\phi}\|_{L^4} + \fr{1}{2} \|\na U\|_{L^2}^2 + C\|\na \tilde{\phi}\|_{L^{8}}^2 +    \|\phi \cdot \na \tilde{\phi}\|_{L^4}.
\ee
Integrating in time from $0$ to $t$ and using the boundedness of $\na U$ in $L_{loc}^2(0, \infty; L^2(\TT^2))$ derived in Step 1, we obtain uniform in  bounds for  the $L^4$ norm of $Q$. 

\textbf{Step 3. Bounds for the velocity in $L_{loc}^2(0,\infty; H^2(\TT^2))$.}
Taking the $L^2$ inner product of the equation obeyed by $U$ with $-\Delta U$, we get
\beg{align}
\fr{1}{2} \fr{d}{dt} \|\na U\|_{L^2}^2
+ \|\Delta U\|_{L^2}^2
&= \int\limits_{\TT^2} (U + \phi) \cdot \na (U + \phi) \cdot \Delta U dx
+ \int\limits_{\TT^2} (Q + \tilde{\phi}) R  (Q + \tilde{\phi}) \cdot \Delta U dx \nonumber
\\&- \int\limits_{\TT^2} f \cdot \Delta U dx
- \int\limits_{\TT^2} \Delta {\phi} \cdot \Delta U dx.
\end{align}
Since the trace of $M^TM^2$ vanishes for any two-by-two traceless matrix $M$, we have
\beg{align}
&\left|\int\limits_{\TT^2} (U + \phi) \cdot \na (U + \phi) \cdot \Delta U dx \right|
= \left|\int\limits_{\TT^2} (U + \phi) \cdot \na (U + \phi) \cdot \Delta \phi dx \right| \nonumber
\\&= \left|\int\limits_{\TT^2} (U  \cdot \na U) \cdot \Delta \phi dx  + \int\limits_{\TT^2} (U  \cdot \na  \phi) \cdot \Delta \phi dx  + \int\limits_{\TT^2} ( \phi \cdot \na U) \cdot \Delta \phi dx \right| \nonumber
\\&\le \|U\|_{L^4}\|\na U\|_{L^4} \|\Delta \phi\|_{L^2}
+ \|U\|_{L^4} \|\na \phi\|_{L^{4}} \|\Delta \phi\|_{L^2}
+ \|\phi\|_{L^2}\|\na U\|_{L^4} \|\Delta \phi\|_{L^2} \nonumber
\\&\le C\|\na U\|_{L^2}^{\fr{3}{2}} \|\Delta U\|_{L^2}^{\fr{1}{2}} \|\Delta \phi\|_{L^2}
+ \|\Delta U\|_{L^2} \|\na \phi\|_{L^{4}} \|\Delta \phi\|_{L^2}
+ \|\phi\|_{L^2}\|\Delta U\|_{L^2} \|\Delta \phi\|_{L^2} \nonumber
\\&\le \fr{1}{4} \|\Delta U\|_{L^2}^2 + (\|\na U\|_{L^2}^2 + \|\na \phi\|_{L^{4}}^2 + \|\phi\|_{L^2}^2) \|\Delta \phi\|_{L^2}^2.
\end{align} 
We obtain 
\beg{align}
 \fr{d}{dt} \|\na U\|_{L^2}^2 
+  \|\Delta U\|_{L^2}^2
&\le  (\|\na U\|_{L^2}^2 + \|\na \phi\|_{L^{4}}^2 + \|\phi\|_{L^2}^2) \|\Delta \phi\|_{L^2}^2  \nonumber
\\&+ C\|Q + \tilde{\phi}\|_{L^4}^4
+ C\|f\|_{L^2}^2
+ C\|\Delta {\phi}\|_{L^2}^2.
\end{align}
We integrate in time from $0$ to $t$ and we use the bounds derived in Step 1 and Step 2 to obtain uniform bounds for $\|\na U\|_{L^2}$ and $\int_{0}^{t}\|\Delta U\|_{L^2}^2 ds$.  

\vspace{0.5cm}

\textbf{Acknowledgment.} The research of N.E.G.H. was partially supported by NSF-DMS-1816551, NSF-DMS-2108790. The research of M.I. was partially supported by NSF-DMS-2204614.

\end{document}